\documentclass[10pt,reqno]{amsart}

\usepackage[T1]{fontenc}
\usepackage[utf8]{inputenc}
\usepackage[margin=0.65in]{geometry}
\usepackage{lmodern}     

\usepackage{amsmath}
\usepackage{amssymb}
\usepackage{amsthm}
\usepackage{mathtools}
\usepackage{mathrsfs}
\usepackage{bm}

\usepackage{graphicx}
\usepackage{booktabs}
\usepackage{enumitem}
\usepackage{microtype}
 \usepackage{tikz}
\usetikzlibrary{arrows.meta}
\AtBeginDocument{\hypersetup{pdfauthor={Kailun Chen},
pdftitle={Half-space polynuclear growth: KPZ universality}}}

\usepackage[hidelinks]{hyperref}
\usepackage[nameinlink,noabbrev]{cleveref}
\crefname{section}{Section}{Sections}
\Crefname{section}{Section}{Sections}
\crefname{subsection}{Section}{Sections}
\Crefname{subsection}{Section}{Sections}
\crefname{subsubsection}{Section}{Sections}
\Crefname{subsubsection}{Section}{Sections}

\allowdisplaybreaks
\numberwithin{equation}{section}
\theoremstyle{plain}
\newtheorem{theorem}{Theorem}[section]
\newtheorem{proposition}[theorem]{Proposition}
\newtheorem{lemma}[theorem]{Lemma}
\newtheorem{corollary}[theorem]{Corollary}

\theoremstyle{definition}
\newtheorem{definition}[theorem]{Definition}

\theoremstyle{remark}
\newtheorem{remark}[theorem]{Remark}

\DeclareMathOperator{\Pf}{Pf}

\title[Half-space PNG: KPZ universality]{Half-space polynuclear growth: KPZ universality}

\author{Kailun Chen}
\thanks{%
Department of Mathematics, National University of Singapore,
2 Science Drive 2, Singapore 117543.
\newline
Email: \href{mailto:kailunxsx@gmail.com}
{\texttt{kailunxsx@gmail.com}}.
\newline
Homepage: \url{https://sites.google.com/view/kailunchenshomepage/home}.
}

\makeatletter
\def\@setaddresses{\par\smallskip\begingroup\footnotesize
 \def\author##1{}%
 \def\address##1##2{\par\noindent\textsc{##2}\enspace}%
 \def\email##1##2{\href{mailto:##2}{\texttt{##2}}}%
 \def\urladdr##1##2{}%
 \addresses\par\endgroup}
\makeatother

\begin{document}
\setlength{\abovedisplayskip}{2pt plus 2pt minus 1pt}
\setlength{\belowdisplayskip}{2pt plus 2pt minus 1pt}
\setlength{\abovedisplayshortskip}{2pt plus 1pt}
\setlength{\belowdisplayshortskip}{2pt plus 1pt minus 1pt}

\begin{abstract}
We study polynuclear growth (PNG) on the half-line with Poisson nucleations in the bulk and at the boundary.
For every nonnegative boundary nucleation rate, we construct stationary measures and jointly stationary couplings under
a common graphical environment, with profiles pinned at the origin. We realize these couplings as a measurable stationary
PNG horizon and prove convergence to the half-space stationary horizon under Kardar–Parisi–Zhang (KPZ) scaling. The
horizon convergence includes its base profile and positive representing measure, together with joint profile convergence at
every prescribed finite parameter list. We also prove convergence of the height evolution from general upper-semicontinuous
initial profiles to the half-space KPZ fixed point. The associated passage field converges locally uniformly to the half-space
directed landscape, jointly with evolutions from finitely many initial profiles driven by the same environment. The
convergence results cover the critical boundary window and every fixed subcritical boundary rate, including zero. All fixed
subcritical rates produce the same limiting laws.
\end{abstract}

\maketitle
\enlargethispage{2pt}
\tableofcontents

\section{Introduction}\label{sec:introduction}

\subsection{Overview}\label{subsec:overview}
Many one-dimensional random growth and interacting particle models
share a characteristic large-scale behavior: at large time $t$, height
fluctuations are of order $t^{1/3}$ and spatial correlations extend
over distances of order $t^{2/3}$. This $1{:}2{:}3$ scaling is
characteristic of the Kardar--Parisi--Zhang (KPZ) universality class
\cite{kardar1986dynamic}. Three related objects describe the expected
universal limits. The KPZ fixed point describes height evolution from
general initial profiles. The directed landscape is a random directed
metric that records passage values between space--time points and
couples these evolutions. The stationary horizon is the jointly
stationary coupling of profiles with different slopes for the KPZ
fixed point, after recentering at the origin.

Early progress in KPZ universality concerned fluctuation laws for
special initial conditions. Baik, Deift, and Johansson proved that
rescaled lengths of longest increasing subsequences of uniform random
permutations converge to the Tracy--Widom law of the Gaussian unitary
ensemble (GUE) \cite{baik1999distribution}; Johansson established
corresponding one-point laws for last-passage percolation and
step-initial TASEP \cite{johansson2000shape}. To describe spatial
correlations, Pr\"ahofer and Spohn introduced multilayer PNG and
identified the Airy process in the PNG droplet \cite{prahofer2002scale}.
Its one-point marginals recover the GUE Tracy--Widom law. Johansson
proved functional convergence for discrete PNG
\cite{johansson2003discrete}. The multilayer construction also led to
the Airy line ensemble, whose top curve is the Airy process. Corwin
and Hammond established its continuous, nonintersecting realization
and, after a parabolic shift, its Brownian Gibbs property
\cite{corwin2014brownian}. For general initial profiles, growth at an
endpoint depends on competition among paths from different starting
locations. The Airy sheet couples passage values as both starting and
ending positions vary. Corwin, Quastel, and Remenik proposed a
variational formula for universal height evolution, maximizing initial
height plus an Airy-sheet passage value
\cite{corwin2015renormalization}. Building on earlier determinantal
representations \cite{sasamoto2005spatial,borodin2007fluctuation},
Matetski, Quastel, and Remenik constructed the KPZ fixed point as a
Feller Markov process by taking the scaling limit of TASEP transition
probabilities from general initial profiles \cite{matetski2021kpz}.
Dauvergne, Ortmann, and Vir\'ag characterized the Airy-sheet law using
the Airy line ensemble and constructed the directed landscape as the
locally uniform limit of Brownian last-passage percolation, together
with convergence of maximizing paths to geodesics in the directed
landscape \cite{dauvergne2022directed}. The directed landscape records
both endpoint positions and times; passage values compose by
maximization over intermediate positions, and disjoint time intervals
carry independent passage fields. Using one-sided reflected Brownian
motions, Nica, Quastel, and Remenik showed that maximizing initial
height plus a passage value in the directed landscape gives the KPZ
fixed point \cite{nica2020reflected}. The directed
landscape thus couples evolutions from different initial profiles. Subsequent results established fixed-point convergence for
the O'Connell--Yor polymer \cite{virag2020heat}; convergence to the
directed landscape and the KPZ fixed point for classical last-passage
models and PNG \cite{dauvergne2021scaling} and the KPZ equation
\cite{wu2026kpzlandscape}; fixed-point convergence for ASEP and the
stochastic six-vertex model, with directed landscape convergence for
their colored versions \cite{aggarwal2024kpz,aggarwal2024colored}; and
directed landscape convergence for the log-gamma polymer
\cite{zhang2025loggamma}. Dauvergne and Zhang characterized the
directed landscape coupling by fixed-point laws, independent time
increments, monotonicity, and commutation with vertical shifts
\cite{dauvergne2024characterization}, providing a framework for deriving
directed landscape convergence from fixed-point convergence and
control of coupled evolutions.

The stationary horizon describes the joint stationary behavior of
profiles with different slopes. Its origins lie in Busemann functions,
obtained from limits of passage-time differences in a prescribed
direction. Fan and Sepp\"al\"ainen gave a queueing description of the
joint Busemann profiles of exponential last-passage percolation
\cite{fan2020joint}. Building on this description, Busani introduced
the stationary horizon as their diffusive scaling limit, with
convergence of the entire process indexed by drift
\cite{busani2024diffusive}. Each profile is a two-sided Brownian motion
with its prescribed drift. The same coupling independently appeared
in the Brownian last-passage Busemann process studied by
Sepp\"al\"ainen and Sorensen \cite{seppalainen2023global}.
Busani, Sepp\"al\"ainen, and Sorensen then proved that the stationary
horizon is jointly invariant under KPZ fixed-point evolutions driven
by one directed landscape, with each evolved profile recentered at
the origin. Under prescribed asymptotic slope conditions, it is the
unique jointly invariant law and an attractor. They also identified
it with the directed landscape's Busemann profiles, connecting
stationary growth with semi-infinite geodesics
\cite{busani2024stationary}. Their subsequent convergence criterion
uses directed landscape convergence, joint invariance, and control of
maximizers to identify the stationary horizon as the limit. It yields
joint convergence at finitely many slopes for Brownian, exponential,
geometric, and Poisson last-passage percolation, the Poisson lines
model, and the Sepp\"al\"ainen--Johansson model
\cite{busani2024scaling}. Convergence at finitely many drift parameters
also holds for the KPZ horizon, which describes jointly invariant
profiles for the KPZ equation, in the low-temperature limit
\cite{groathouse2025jointly}, and for stationary ASEP and stochastic
six-vertex profiles \cite{aggarwal2024colored}.
The stationary profiles encoded by the TASEP speed process also
converge to the stationary horizon, with convergence of the entire
process indexed by drift \cite{busani2026speed}.

In half-space, boundary interactions produce parameter-dependent
analogues of these objects. Zhang obtained half-space TASEP transition
Pfaffians and their fixed-time KPZ scaling limits \cite{zhang2024tasep}.
Dauvergne and Zhang constructed the half-space stationary horizon from
jointly stationary last-passage profiles and obtained the half-space
directed landscape as the scaling limit of exponential last-passage
percolation and a family of Poisson-avoiding metrics generalizing
colored TASEP \cite{dauvergne2026directed}. Their variational
construction realizes Zhang's transition laws as a continuous Markov
process \cite{dauvergne2026directed}, the half-space KPZ fixed
point. Their directed landscape characterization is the universality
criterion used here.

We study half-space PNG, a growth process on the half-line driven by
bulk and boundary Poisson nucleations \cite{sasamoto2004fluctuations}.
We construct its stationary measures, jointly stationary couplings, and
PNG horizon. We prove convergence of the PNG horizon, height evolution,
and common passage field to the half-space stationary horizon, KPZ
fixed point, and directed landscape, respectively, in the critical
boundary window and at fixed subcritical boundary rates.
We prove the half-space KPZ fixed-point limit from the transition Pfaffian
of Chen and Zhang \cite{chenZhangPNG} and Zhang's limiting formula
\cite{zhang2024tasep}, together with tightness and a Markov argument.
We obtain the directed landscape limit through the characterization
of Dauvergne and Zhang \cite{dauvergne2026directed}.

\subsection{Model and setting}\label{subsec:model-and-setting}

We work on the half-line \(\mathbb R_{+}:=[0,\infty)\). We use the normalization of
full-space PNG from \cite{matetski2025polynuclear}, in which bulk
nucleations have intensity \(2\), and add a source at the origin as in
\cite{sasamoto2004fluctuations,chenZhangPNG}. 

\begin{definition}[Half-line height profiles]
\label{def:half-line-height-profiles}
Let $\mathrm{UC}_+$ consist of locally bounded-above,
upper-semicontinuous functions $h:[0,\infty)\to\mathbb Z\cup\{-\infty\}$,
with local Hausdorff convergence of their closed hypographs
$\operatorname{hypo}(h)=\{(x,r):x\geq0,\ r\leq h(x)\}$.
\end{definition}
Upper semicontinuity selects the upper value at a step and includes
height-one nucleation spikes. The value $-\infty$ allows localized
initial data. Write $u_+=\max(u,0)$.

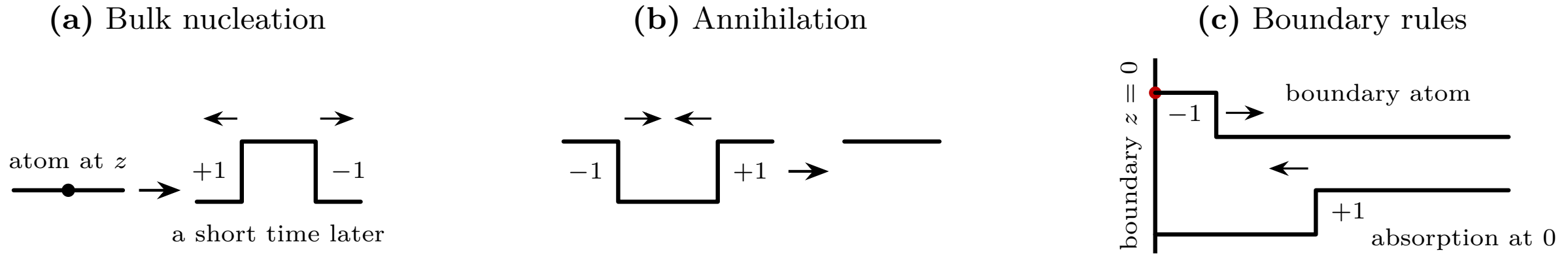
\begin{figure}[t]
\centering
\begin{tikzpicture}[
    x=0.65cm,
    y=0.58cm,
    profile/.style={very thick,line cap=round,line join=round},
    motion/.style={-{Stealth[length=2mm]},semithick},
    evolve/.style={-{Stealth[length=2.4mm]},thick},
    event/.style={circle,fill=black,inner sep=1.35pt}
]

\node[font=\small] at (2.8,4.2)
    {\textbf{(a)} Bulk nucleation};
\draw[profile] (0.10,1.25)--(1.80,1.25);
\node[event] at (0.95,1.25) {};
\node[font=\scriptsize,above] at (0.95,1.42)
    {atom at \(z\)};
\draw[evolve] (2.05,1.25)--(2.70,1.25);
\draw[profile]
    (2.95,1.05)--(3.65,1.05)
    --(3.65,2.10)--(4.80,2.10)
    --(4.80,1.05)--(5.50,1.05);
\node[font=\scriptsize,left] at (3.61,1.58) {\(+1\)};
\node[font=\scriptsize,right] at (4.84,1.58) {\(-1\)};
\draw[motion] (3.58,2.50)--(3.05,2.50);
\draw[motion] (4.87,2.50)--(5.40,2.50);
\node[font=\scriptsize,below] at (4.22,0.86)
    {a short time later};

\begin{scope}[xshift=1.5cm]
\node[font=\small] at (9.25,4.2)
    {\textbf{(b)} Annihilation};
\draw[profile]
    (6.35,2.10)--(7.20,2.10)
    --(7.20,1.05)--(8.75,1.05)
    --(8.75,2.10)--(9.60,2.10);
\node[font=\scriptsize,left] at (7.16,1.57) {\(-1\)};
\node[font=\scriptsize,right] at (8.79,1.57) {\(+1\)};
\draw[motion] (7.30,2.50)--(7.88,2.50);
\draw[motion] (8.65,2.50)--(8.07,2.50);
\draw[evolve] (9.85,1.58)--(10.48,1.58);
\draw[profile] (10.72,2.10)--(12.20,2.10);
\end{scope}

\begin{scope}[xshift=3.0cm]
\node[font=\small] at (16.0,4.2)
    {\textbf{(c)} Boundary rules};
\draw[very thick] (13.25,0.15)--(13.25,3.55);
\node[font=\scriptsize,rotate=90] at (12.87,1.85)
    {boundary \(z=0\)};

\fill[red!75!black] (13.25,2.95) circle (1.8pt);
\draw[profile]
    (13.25,2.95)--(14.20,2.95)
    --(14.20,2.18)--(18.75,2.18);
\node[font=\scriptsize,right] at (13.24,2.57) {\(-1\)};
\draw[motion] (14.35,2.60)--(14.95,2.60);
\node[font=\scriptsize,anchor=west] at (15.08,2.92)
    {boundary atom};

\draw[profile]
    (13.25,0.48)--(15.75,0.48)
    --(15.75,1.25)--(18.75,1.25);
\node[font=\scriptsize,right] at (15.79,0.87) {\(+1\)};
\draw[motion] (15.63,1.63)--(15.02,1.63);
\node[font=\scriptsize,anchor=west] at (16.42,0.40)
    {absorption at \(0\)};
\end{scope}

\end{tikzpicture}
\caption{Local rules of half-space PNG.
\textup{(a)} A bulk atom creates a coincident up-step/down-step pair,
which separates into an expanding island.
\textup{(b)} A down-step followed by an up-step annihilates on contact.
\textup{(c)} A boundary atom raises the height at the origin and emits
a down-step, whereas an incoming up-step is absorbed at the origin.}
\label{fig:half-space-png-local-rules}
\end{figure}

\begin{definition}[Half-space PNG]\label{def:half-space-png}
Fix a boundary parameter \(c\geq0\).  The half-space PNG process with
boundary rate \(c\) is the evolution on \(\mathrm{UC}_{+}\) determined
by the following rules.

\begin{enumerate}[label=\textup{(\roman*)},leftmargin=2.6em]
\item
In the absence of new nucleations, a profile \(h\) evolves for time
\(s\geq0\) according to
\begin{equation}
 (\mathsf D_s h)(x)
    :=
    \sup_{y\in[(x-s)_{+},\,x+s]}h(y),
    \qquad x\in\mathbb R_{+}.
    \label{eq:half-space-png-deterministic-evolution}
\end{equation}
For locally finite unit-step profiles, an up-step, where the height
rises by one from left to right, moves to the left with speed one,
whereas a down-step moves to the right with speed one. A down-step
followed by an up-step annihilates when the two steps meet. An up-step
that reaches the origin is absorbed, so that the upper terrace has
expanded to the boundary.

\item
Bulk nucleations are the atoms of a Poisson point process
\(\Pi_{\mathrm{bulk}}\) on
\(\mathbb R_{>0}\times\mathbb R_{>0}\) with intensity
\(2\,\mathrm dt\,\mathrm dz\). At an atom \((t,z)\), the creation
map \(C_z\) sets \((C_z h)(z)=h(z)+1\) and leaves the profile
unchanged elsewhere.
Where $h(z)$ is finite, this adds a height-one upper spike. On a flat
terrace it creates a coincident up-step/down-step pair; immediately
afterwards the two steps separate and form an expanding island.

\item
Independently of \(\Pi_{\mathrm{bulk}}\), boundary nucleations are the
atoms of a Poisson process \(\Pi_{\partial}\) of rate \(c\) on
\(\mathbb R_{>0}\). At each boundary atom, the creation map \(C_0\)
sets \((C_0h)(0)=h(0)+1\) and leaves the profile unchanged for \(x>0\).
Where $h(0)$ is finite, the raised boundary height propagates into
the half-line under the deterministic evolution; on a flat terrace it
emits a right-moving down-step. When an incoming up-step reaches the origin, it is absorbed, raising the height at \(0\) by one.
\end{enumerate}
We use the convention \(-\infty+1=-\infty\): an atom at a location
whose current height is $-\infty$ leaves the profile unchanged.
Figure~\ref{fig:half-space-png-local-rules} illustrates these rules
on finite terraces.
\end{definition}

The maps $\mathsf D_s$, $C_z$, and $C_0$ preserve $\mathrm{UC}_+$,
and $\mathsf D_s\mathsf D_r=\mathsf D_{s+r}$. Indeed, for each integer
$a$, the superlevel set of $\mathsf D_s h$ is the closed
$s$-neighbourhood of $\{h\geq a\}$ in $\mathbb R_+$, while a creation
map either leaves a superlevel set unchanged or adds its creation
location. To construct the evolution on $[0,R]$ up to time $T$, only
$h_0|_{[0,R+T]}$, bulk atoms in $(0,T]\times(0,R+T]$, and boundary atoms
in $(0,T]$ can contribute. Their total number $N_{R,T}$ is almost
surely finite, and the heights in this window are bounded above by
$\sup_{[0,R+T]}h_0+N_{R,T}$. Ordering these atoms in time and
alternating deterministic evolution with the creation maps therefore
defines a pathwise evolution in $\mathrm{UC}_+$. Finite propagation
makes these constructions agree on overlapping windows. We denote the
resulting height by \(h^{(c)}(t,x;h_0)\), \(t,x\geq0\),
for \(h_0\in\mathrm{UC}_{+}\). Evolving several initial profiles with
one environment \(\Pi=(\Pi_{\mathrm{bulk}},\Pi_{\partial})\) defines the
\emph{graphical}, or \emph{basic}, coupling. We write
$h^{(c)}(t,x;h_0;\Pi)$ when specifying the environment.

The value at \((t,x)\) depends on the initial profile only through
the interval $I_{t,x}$ and on the bulk atoms in the backward light cone
$\mathcal C_{t,x}$, where
\begin{equation*}
 I_{t,x}:=[(x-t)_+,x+t],\qquad
 \mathcal C_{t,x}:=\{(s,z):0<s\leq t,\ z>0,\ |z-x|\leq t-s\}.
\end{equation*}
Boundary atoms can contribute only for $0<s\leq t-x$, an interval
of length $(t-x)_+$.

\begin{definition}[Narrow wedge and point-to-point field]
\label{def:half-space-png-point-to-point-field}
For \(y\in\mathbb R_{+}\), the narrow wedge based at \(y\) is the
profile \(\mathfrak d_y\in\mathrm{UC}_{+}\) with
\(\mathfrak d_y(y)=0\) and \(\mathfrak d_y(z)=-\infty\) for \(z\neq y\).
The point-to-point half-space PNG field is
\(\mathcal A_t^{(c)}(x,y):=h^{(c)}(t,x;\mathfrak d_y)\).
We write $\mathcal A_t^{(c)}(x,y;\Pi)$ when specifying the graphical
environment.
\end{definition}

Equivalently, \(\mathcal A_t^{(c)}(x,y)\) is the largest number of
bulk and boundary atoms that can be collected by a path from
\((y,0)\) to \((x,t)\) which remains in the half-space and whose
spatial coordinate is \(1\)-Lipschitz in time. Collected atoms have
time coordinates in \((0,t]\). If
\(|x-y|>t\), there is no such path and we set
\(\mathcal A_t^{(c)}(x,y)=-\infty\). A feasible path that collects
no atoms has score zero.

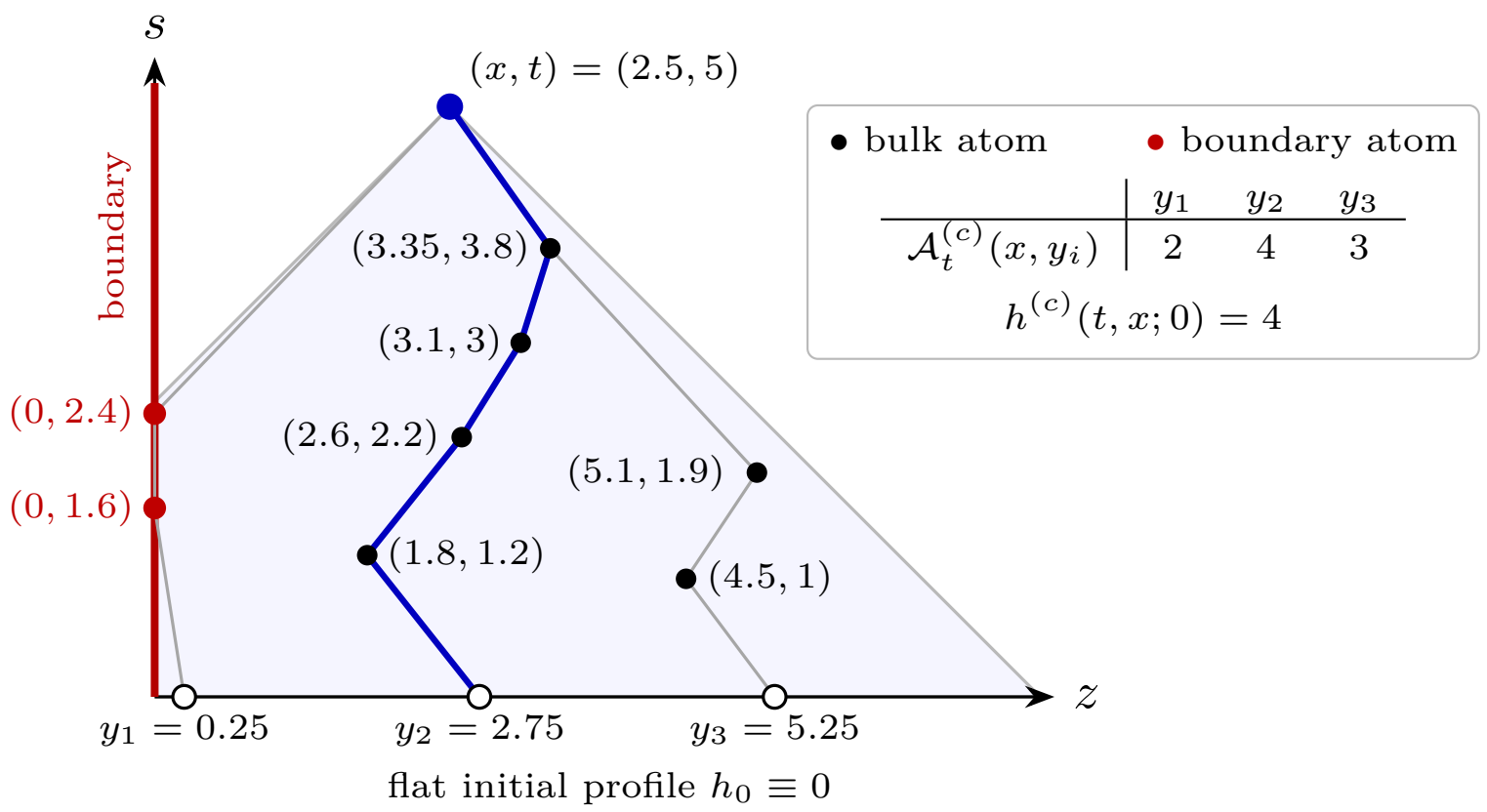
\begin{figure}[t]
\centering
\begin{tikzpicture}[
    x=0.80cm,
    y=0.80cm,
    axis/.style={-{Stealth[length=2mm]},semithick},
    cone/.style={gray!55,semithick},
    pointpath/.style={gray!70,semithick},
    maximizing/.style={blue!75!black,very thick},
    bulk/.style={circle,fill=black,inner sep=1.40pt},
    bdry/.style={circle,fill=red!75!black,inner sep=1.55pt},
    start/.style={
        circle,
        draw=black,
        fill=white,
        inner sep=1.55pt,
        line width=0.6pt
    },
    endpoint/.style={
        circle,
        fill=blue!75!black,
        inner sep=1.80pt
    },
    score/.style={
        draw=gray!55,
        fill=white,
        rounded corners=2pt,
        inner sep=4pt,
        align=center,
        font=\scriptsize
    }
]

\fill[blue!4]
    (0,0)--(7.5,0)--(2.5,5)--(0,2.5)--cycle;
\draw[cone]
    (0,0)--(7.5,0)--(2.5,5)--(0,2.5);

\draw[axis] (0,0)--(7.62,0) node[right] {\(z\)};
\draw[axis] (0,0)--(0,5.42) node[above] {\(s\)};
\draw[red!70!black,line width=1.5pt]
    (0,0)--(0,5.20);
\node[
    font=\scriptsize,
    rotate=90,
    text=red!70!black
] at (-0.34,3.90) {boundary};

\draw[pointpath]
    (0.25,0)--(0,1.60)--(0,2.40)--(2.50,5.00);

\draw[pointpath]
    (5.25,0)--(4.50,1.00)--(5.10,1.90)
    --(3.35,3.80)--(2.50,5.00);

\draw[maximizing]
    (2.75,0)--(1.80,1.20)--(2.60,2.20)
    --(3.10,3.00)--(3.35,3.80)--(2.50,5.00);

\node[bulk] at (1.80,1.20) {};
\node[font=\scriptsize,anchor=east] at (3.47,1.20) {\((1.8,1.2)\)};
\node[bulk] at (2.60,2.20) {};
\node[font=\scriptsize,anchor=east] at (2.57,2.20) {\((2.6,2.2)\)};
\node[bulk] at (3.10,3.00) {};
\node[font=\scriptsize,anchor=east] at (3.10,3.00) {\((3.1,3)\)};
\node[bulk] at (3.35,3.80) {};
\node[font=\scriptsize,anchor=west] at (1.50,3.80) {\((3.35,3.8)\)};
\node[bulk] at (4.50,1.00) {};
\node[font=\scriptsize,anchor=east] at (5.90,1.00) {\((4.5,1)\)};
\node[bulk] at (5.10,1.90) {};
\node[font=\scriptsize,anchor=west] at (3.33,1.90) {\((5.1,1.9)\)};

\node[bdry] at (0,1.60) {};
\node[
    font=\scriptsize,
    anchor=west,
    text=red!75!black
] at (-1.4,1.60) {\((0,1.6)\)};
\node[bdry] at (0,2.40) {};
\node[
    font=\scriptsize,
    anchor=west,
    text=red!75!black
] at (-1.4,2.40) {\((0,2.4)\)};

\node[start] at (0.25,0) {};
\node[start] at (2.75,0) {};
\node[start] at (5.25,0) {};
\node[endpoint] at (2.50,5.00) {};

\node[font=\scriptsize,above right]
    at (2.50,5.00)
    {\((x,t)=(2.5,5)\)};

\node[font=\scriptsize,below] at (0.25,0) {\(y_1=0.25\)};
\node[font=\scriptsize,below] at (2.75,0) {\(y_2=2.75\)};
\node[font=\scriptsize,below] at (5.25,0) {\(y_3=5.25\)};

\node[font=\scriptsize,anchor=north]
    at (3.85,-0.48)
    {flat initial profile \(h_0\equiv0\)};

\node[score,anchor=north west] at (5.52,5.02) {%
    \(\bullet\) bulk atom
    \qquad
    \textcolor{red!75!black}{\(\bullet\)} boundary atom\\[3pt]
    \(\begin{array}{c|ccc}
        &y_1&y_2&y_3\\ \hline
        \mathcal A_t^{(c)}(x,y_i)&2&4&3
      \end{array}\)\\[3pt]
    \(\displaystyle
      h^{(c)}(t,x;0)=4\)};

\end{tikzpicture}

\caption{Last passage from the flat profile $h_0\equiv0$.
All atoms in the backward cone are shown in coordinates $(z,s)$.
Each displayed path maximizes from its starting point; the blue path
also maximizes over the starting location.}
\label{fig:half-space-png-last-passage}
\end{figure}

Each of the maps \(\mathsf D_s\), \(C_z\), and \(C_0\) preserves
pointwise suprema whenever the supremum belongs to $\mathrm{UC}_+$,
and commutes with integer vertical shifts.
The profile identically equal to $-\infty$ is absorbing. Applying these
identities to the pointwise decomposition
\(h_0=\sup_{y\geq0}\{h_0(y)+\mathfrak d_y\}\)
through the finitely many updates in a dependence cone gives the
pathwise variational formula
\begin{equation}
 h^{(c)}(t,x;h_0)
    =
    \sup_{y\in I_{t,x}}
    \bigl\{h_0(y)+\mathcal A_t^{(c)}(x,y)\bigr\}.
    \label{eq:half-space-png-max-plus-formula}
\end{equation}
Thus a general initial profile is propagated by adding to each
possible starting height \(h_0(y)\) the point-to-point passage value
from \(y\), and then selecting the largest contribution. Figure~\ref{fig:half-space-png-last-passage} illustrates \eqref{eq:half-space-png-max-plus-formula} for flat initial data.

\subsection{Stationary measures}
\label{subsec:intro-stationary-measures}

We formulate stationarity after pinning at the origin. Write
\begin{equation*}
 \mathrm{UC}_{+,0}:=\{h\in\mathrm{UC}_+:h(0)=0\},\qquad
 \operatorname{Rec}(h):=h-h(0),
\end{equation*}
where recentering is defined when $h(0)$ is finite. Equip
$\mathrm{UC}_{+,0}$ with its induced Borel $\sigma$-field.

\begin{definition}[Stationary measure]
\label{def:intro-half-space-png-stationary-measure}
Fix $c\geq0$.  A law $\nu$ on $\mathrm{UC}_{+,0}$ is stationary for
half-space PNG if, for $H\sim\nu$ independent of a future graphical
environment $\Pi$,
\begin{equation*}
 \operatorname{Rec}\bigl(h^{(c)}(t,\cdot;H;\Pi)\bigr)
 \overset{\mathrm d}{=}H,\qquad t\geq0.
\end{equation*}
\end{definition}

Set $r_c:=c\vee1$, $\mathcal G_c:=[r_c,\infty)$, and
$\mathbb N_0:=\{0,1,\ldots\}$. Our geometric convention is
$\mathbb P(\operatorname{Geo}(q)=n)=(1-q)q^n$ for $q\in[0,1)$ and
$n\in\mathbb N_0$. A signed Poisson walk with up/down rates $(a,b)$ is $N^+-N^-$,
where $N^+,N^-$ are independent Poisson processes with rates $a,b$,
respectively, started at zero.

For a spectral parameter $g\in\mathcal G_c$, let $S_1^{(g)},S_2^{(g)}$ be independent
signed Poisson walks with respective up/down rates $(g,g^{-1})$ and
$(g^{-1},g)$.  If $g>c$, independently take
$Y^{(c,g)}\sim\operatorname{Geo}(c/g)$ and set
\begin{equation}
 \widetilde H^{(c,g)}(x)
 :=S_2^{(g)}(x)
 +\left(
   \sup_{0\leq y\leq x}
   \{S_1^{(g)}(y)-S_2^{(g)}(y)\}
   -Y^{(c,g)}
  \right)_{+}.
 \label{eq:intro-half-space-png-stationary-profile-raw}
\end{equation}
For locally bounded $f$, its upper-semicontinuous regularization is
\begin{equation*}
 (\operatorname{usc}f)(x):=\lim_{\varepsilon\downarrow0}
 \sup\{f(y):y\geq0,\ |y-x|<\varepsilon\}.
\end{equation*}
Define the profiles and their laws by
\begin{align}
 H^{(c,g)}&:=\operatorname{usc}\widetilde H^{(c,g)},
 &&g>c,
 \label{eq:intro-half-space-png-stationary-profile}\\
 H^{(c,c)}&:=\operatorname{usc}S_2^{(c)},
 &&c\geq1,
 \label{eq:intro-half-space-png-boundary-dominated-profile}
\end{align}
and $\nu_{c,g}:=\operatorname{Law}(H^{(c,g)})$.

\begin{theorem}[A stationary family for half-space PNG]
\label{thm:intro-half-space-png-stationary-family}
For every $c\geq0$ and $g\in\mathcal G_c$, the law $\nu_{c,g}$ is
stationary for half-space PNG with boundary rate $c$.
\end{theorem}

Section~\ref{subsec:stationary-measures} proves this by generator
balance. 

\begin{definition}[Jointly stationary measure]
\label{def:intro-half-space-png-joint-stationarity}
Fix $c\geq0$ and $k\in\mathbb N$. A law $\boldsymbol\nu$ on
$(\mathrm{UC}_{+,0})^k$ is \emph{jointly stationary} if, for
$(H_1,\ldots,H_k)\sim\boldsymbol\nu$
independent of one common future environment $\Pi$,
\begin{equation*}
 \left(
  \operatorname{Rec}\bigl(h^{(c)}(t,\cdot;H_i;\Pi)\bigr)
 \right)_{i=1}^{k}
 \overset{\mathrm d}{=}(H_i)_{i=1}^{k},\qquad t\geq0.
\end{equation*}
\end{definition}

For $g_1>\cdots>g_k\geq r_c$, write
$\boldsymbol g=(g_1,\ldots,g_k)$.
When $g_k>c$, let $(R_i^{(c;\boldsymbol g)})_{i=1}^k$ be the
profiles of Definition~\ref{def:half-space-png-joint-candidate-law};
Figure~\ref{fig:half-space-png-k-two-candidate-law} illustrates their
shared line and bridge environment for $k=2$.
At $g_k=c\geq1$, omit the singular terminal bridge and use the
shifted source in \eqref{eq:half-space-png-singular-shifted-source}.
For an arbitrary list $g_1,\ldots,g_k\in\mathcal G_c$, let
$a_1>\cdots>a_\ell$ be its distinct values, set
$\boldsymbol a=(a_1,\ldots,a_\ell)$, and define $\kappa(i)$ by
$g_i=a_{\kappa(i)}$. Set
\begin{equation*}
 \nu_{c;g_1,\ldots,g_k}:=\operatorname{Law}
 \bigl((R_{\kappa(i)}^{(c;\boldsymbol a)})_{i=1}^k\bigr).
\end{equation*}

\begin{theorem}[A jointly stationary family for half-space PNG]
\label{thm:intro-half-space-png-joint-stationary-family}
For every $c\geq0$, $k\in\mathbb N$, and
$g_1,\ldots,g_k\in\mathcal G_c$, the law $\nu_{c;g_1,\ldots,g_k}$ is
jointly stationary for half-space PNG with boundary rate $c$.  Its $i$-th marginal is $\nu_{c,g_i}$, and the family is consistent under deleting, permuting, and repeating parameters.
\end{theorem}

\begin{remark}[Conjecture]
\label{rem:expected-classification}
We expect that the stationary laws
in Theorem~1.5 are all the extremal stationary measures of
half-space PNG, and that the jointly stationary laws in Theorem~1.7
are all the extremal jointly stationary measures under the
basic coupling. Proving these classifications requires further work,
which we do not pursue here. For the corresponding classifications
in half-space geometric last-passage percolation, see
\cite{das2026invariant}.
\end{remark}

Section~\ref{subsec:half-space-png-horizon} constructs a jointly
measurable field whose finite-parameter distributions are the jointly
stationary laws $\nu_{c;g_1,\ldots,g_k}$.
Uniqueness below concerns regular representations of this family
in the sense of Definition~\ref{def:half-space-png-horizon}.

\begin{theorem}[Half-space PNG horizon]
\label{thm:intro-half-space-png-horizon}
For every $c\geq0$, there exists a regular representation
$(\mathfrak b_c^{\mathrm{PNG}},\mu_c^{\mathrm{PNG}})$, unique in law,
with the following properties.
\begin{enumerate}[label=\textup{(\roman*)},leftmargin=*]

\item \textup{Construction.}
The base profile $\mathfrak b_c^{\mathrm{PNG}}$ is a pinned, locally
finite, integer-valued c\`adl\`ag step function.
The measure $\mu_c^{\mathrm{PNG}}$ is a positive locally finite
Borel measure on $\mathcal G_c\times\mathbb R_+$, zero on
$\{r_c\}\times\mathbb R_+$ and $\mathcal G_c\times\{0\}$.
Define
\begin{align*}
 H_c^{\mathrm{PNG},\to}(g,x)
 :=\mathfrak b_c^{\mathrm{PNG}}(x)
   +\mu_c^{\mathrm{PNG}}\bigl((r_c,g]\times(0,x]\bigr),\qquad
 H_c^{\mathrm{PNG}}(g,\cdot)
 :=\operatorname{usc}H_c^{\mathrm{PNG},\to}(g,\cdot).
\end{align*}
We call $H_c^{\mathrm{PNG}}$ the \emph{half-space PNG horizon}.

\item \textup{Finite-parameter laws and regularity.}
The field $H_c^{\mathrm{PNG}}$ is jointly measurable.
For every deterministic finite list $g_1,\ldots,g_k\in\mathcal G_c$,
\begin{equation*}
 \operatorname{Law}\left(
 H_c^{\mathrm{PNG}}(g_1,\cdot),\ldots,
 H_c^{\mathrm{PNG}}(g_k,\cdot)
 \right)=\nu_{c;g_1,\ldots,g_k}.
\end{equation*}
Almost surely, every coordinate is a locally finite unit-step
profile in $\mathrm{UC}_{+,0}$, and
$H_c^{\mathrm{PNG},\to}(g_n,\cdot)$ converges locally uniformly to
$H_c^{\mathrm{PNG},\to}(g,\cdot)$ whenever $g_n\downarrow g$.

\item \textup{Stationarity.}
Let $\Pi$ be a future graphical environment independent of
$H_c^{\mathrm{PNG}}$. For each fixed $t\geq0$, define
\begin{equation*}
 H_{c,t}^{\mathrm{ev}}(g,\cdot)
 :=\operatorname{Rec}\bigl(
 h^{(c)}(t,\cdot;H_c^{\mathrm{PNG}}(g,\cdot);\Pi)\bigr).
\end{equation*}
The coordinates indexed by
$\mathbb Q_c:=\{r_c\}\cup(\mathbb Q\cap(r_c,\infty))$
determine a canonical regular reconstruction
$\widehat H_{c,t}^{\mathrm{PNG}}$.
For each fixed deterministic $g\in\mathcal G_c$,
\begin{equation}
 \widehat H_{c,t}^{\mathrm{PNG}}(g,\cdot)
 =H_{c,t}^{\mathrm{ev}}(g,\cdot)
 \quad\text{almost surely},\qquad
 \widehat H_{c,t}^{\mathrm{PNG}}
 \overset{\mathrm d}{=}H_c^{\mathrm{PNG}}
 \quad\text{in the }(\mathfrak b,\mu)\text{ representation}.
 \label{eq:intro-half-space-png-horizon-stationarity}
\end{equation}
\end{enumerate}
\end{theorem}

\subsection{KPZ universality}
\label{subsec:intro-kpz-universality}

We prove $1{:}2{:}3$ limits for the stationary horizon, the height
evolution from general initial data, and the common passage field.
Use the boundary convention of \cite{dauvergne2026directed} and set
\begin{equation}
 c_\varepsilon(\rho)=e^{\rho\varepsilon^{1/2}}, \quad
 \lambda_\rho=(2\rho)_+\quad(\rho\in\mathbb R),\qquad\qquad
 c_\varepsilon(-\infty)=c_\star\in[0,1),\quad
 \lambda_{-\infty}=0.
 \label{eq:intro-png-critical-boundary-scaling}
\end{equation}
Every fixed $c_\star\in[0,1)$ has the same subcritical limit.
The critical window includes rates above one when $\rho>0$;
a fixed rate $c>1$ is outside these scaling limits. Zhang's
parameter in \cite{zhang2024tasep} is $-\rho$; the affine boundary
parameter in \cite[Section~1.3]{chenZhangPNG} is
$-\rho+O(\varepsilon^{1/2})$.
For the stationary family, scale the spectral parameter $g$ by:
\begin{equation}
 g_\varepsilon(\lambda)=e^{\lambda\varepsilon^{1/2}/2},
 \qquad\lambda\geq\lambda_\rho,
 \label{eq:intro-png-slope-scaling}
\end{equation}
equivalently $g_\varepsilon(\lambda)\geq c_\varepsilon(\rho)\vee1$.
Write $\mathcal H^\rho$ for the half-space stationary horizon of
\cite[Theorem~2.19]{dauvergne2026directed}, restricted here to $x\geq0$,
and $\mathcal L^\rho$, $\mathfrak h^\rho$ for the half-space directed
landscape and KPZ fixed point.

For the passage and height fields below, local uniform convergence
is understood through continuous mesh interpolation, with uninterpolated and
interpolated fields asymptotically uniformly close on each compact set.
We write $C_{\mathrm{loc}}(D)$ for continuous real-valued functions on
$D$ with uniform convergence on compact subsets, and $C_c(D)$ for
continuous functions with compact support in $D$. For a scalar function
$f$ and a compact set $K$, write $\|f\|_K:=\sup_{x\in K}|f(x)|$.
Write $\Longrightarrow$ for convergence in distribution.

\subsubsection{Convergence to the half-space stationary horizon}
\label{subsubsec:intro-stationary-horizon-convergence}
For $x=\varepsilon(n+\theta)$, $n\in\mathbb N_0$, $0\leq\theta<1$,
define the linear interpolation by
$
 (\operatorname{Lin}_\varepsilon h)(x)
 =(1-\theta)h(n)+\theta h(n+1),
$
and set
\begin{equation*}
 \mathcal H_\varepsilon^{\mathrm{PNG},\rho}(\lambda,x)
 =\varepsilon^{1/2}(\operatorname{Lin}_\varepsilon
 H_{c_\varepsilon(\rho)}^{\mathrm{PNG},\to}
       (g_\varepsilon(\lambda),\cdot))(x).
\end{equation*}
Here $H_c^{\mathrm{PNG},\to}$ is the canonical right-continuous representative
from Theorem~\ref{thm:intro-half-space-png-horizon}.
The space $\mathscr D(\lambda_\rho)$ consists of pairs $(\mathfrak b,\mu)$,
where $\mathfrak b$ is continuous with $\mathfrak b(0)=0$ and $\mu$ is a positive locally
finite measure on $[\lambda_\rho,\infty)\times[0,\infty)$ with zero
mass on the two lower faces. Identify the pair with
$G(\lambda,x)=\mathfrak b(x)+\mu((\lambda_\rho,\lambda]\times(0,x])$.
Its topology is local uniform convergence of $\mathfrak b$ and vague convergence
of the bivariate Stieltjes measure $\mu$. This is
\cite[Definition~2.13]{dauvergne2026directed} restricted to $x\geq0$,
with zero mass assigned to the two lower faces.

\begin{theorem}[Stationary horizon limit]
\label{thm:intro-half-space-png-stationary-horizon-limit}
For $\rho\in\mathbb R\cup\{-\infty\}$,
$\mathcal H_\varepsilon^{\mathrm{PNG},\rho}$ is a Borel
$\mathscr D(\lambda_\rho)$-valued random element and
\begin{equation*}
 \mathcal H_\varepsilon^{\mathrm{PNG},\rho}
 \Longrightarrow
 \mathcal H^\rho|_{[\lambda_\rho,\infty)\times\mathbb R_+}
 \quad\text{in }\mathscr D(\lambda_\rho).
\end{equation*}
For every $\lambda_1>\cdots>\lambda_k\geq\lambda_\rho$, also
\begin{equation*}
 (\mathcal H_\varepsilon^{\mathrm{PNG},\rho}(\lambda_i,\cdot))_{i=1}^k
 \Longrightarrow(\mathcal H^\rho(\lambda_i,\cdot))_{i=1}^k
 \quad\text{in }C_{\mathrm{loc}}([0,\infty))^k.
\end{equation*}
\end{theorem}
The proof is in Section~\ref{subsec:stationary-horizon-convergence}.

\subsubsection{Convergence to the half-space KPZ fixed point}
\label{subsubsec:intro-kpz-fixed-point-convergence}
Let $\mathrm{UC}_+^{\mathbb R}$ consist of locally bounded above,
upper-semicontinuous functions $[0,\infty)\to\mathbb R\cup\{-\infty\}$
not identically $-\infty$. Local hypograph convergence means that
$x_\varepsilon\to x$ implies
$\limsup f_\varepsilon(x_\varepsilon)\leq f(x)$, and each finite
$f(x)$ has a recovery sequence attaining its value in the limit.
For $f$ with a linear upper bound and
$\rho\in\mathbb R\cup\{-\infty\}$, let $\mathfrak h^\rho(t,y;f)$
denote the half-space KPZ fixed point started from $f$, realized as a
continuous Markov process in
\cite[Theorem~4.11]{dauvergne2026directed}. Its transition laws are the Fredholm Pfaffians of
\cite{zhang2024tasep}, recalled in Section~\ref{subsubsec:critical-transition-limit}
for finite $\rho$, and their
subcritical limit for $\rho=-\infty$.
Let $\mathcal X_{\mathrm{adm}}\subset\mathrm{UC}_+$ be the state space of
\cite[Section~1.2]{chenZhangPNG}: on every $[0,R_{\mathrm{obs}}]$, a profile is constant
on the open intervals of some finite partition. It contains the
unit-step class $\mathcal X_{\mathrm{step}}$ of
Definition~\ref{def:half-space-png-step-coordinate-space}
and hence supports every stationary law $\nu_{c,g}$.
For deterministic $h_{0,\varepsilon}\in\mathcal X_{\mathrm{adm}}$, set
$f_\varepsilon(x)=\varepsilon^{1/2}h_{0,\varepsilon}(\varepsilon^{-1}x)$
and
\begin{equation*}
 \mathfrak h_\varepsilon^\rho(t,y)
 =\varepsilon^{1/2}\bigl[
 h^{(c_\varepsilon(\rho))}(\varepsilon^{-3/2}t,
       \varepsilon^{-1}y;h_{0,\varepsilon})
       -2\varepsilon^{-3/2}t\bigr].
\end{equation*}

\begin{theorem}[Half-space KPZ fixed-point limit]
\label{thm:intro-half-space-png-kpz-fixed-point-limit}
Fix $\rho\in\mathbb R\cup\{-\infty\}$. Suppose
$f_\varepsilon\to f\in\mathrm{UC}_+^{\mathbb R}$ locally in hypograph
and $f_\varepsilon(x)\leq C(1+x)$ for one finite $C\geq0$.
Then
\begin{equation*}
 \mathfrak h_\varepsilon^\rho\Longrightarrow
 \mathfrak h^\rho(\cdot,\cdot;f)
 \quad\text{locally uniformly on }(0,\infty)\times[0,\infty).
\end{equation*}
\end{theorem}
The proof is in Section~\ref{subsec:kpz-fixed-point-convergence}.

\subsubsection{Convergence to the half-space directed landscape}
\label{subsubsec:intro-directed-metric-convergence}
Take bulk atoms of intensity $2\,ds\,dz$ on
$\mathbb R\times(0,\infty)$ and independent rate-$c$ boundary atoms
on $\mathbb R\times\{0\}$, with union $\Pi^{(c)}$. Put
\begin{equation*}
 \mathbb H=[0,\infty)\times\mathbb R,\qquad
 \mathbb H_\uparrow^2=\{(x,s;y,t)\in\mathbb H^2:s<t\}.
\end{equation*}
\begin{definition}[Half-space PNG directed metric]
\label{def:half-space-png-random-directed-metric}
With $\vartheta_s\Pi=\{(r-s,z):(r,z)\in\Pi\}$ and the passage
field of Definition~\ref{def:half-space-png-point-to-point-field}, set
\begin{equation*}
 \mathcal L_{\mathrm{PNG}}^{(c)}(x,s;y,t)
 =\mathcal A_{t-s}^{(c)}(y,x;\vartheta_s\Pi^{(c)}).
\end{equation*}
It counts the greatest number of distinct atoms in $(s,t]$ collected
by a $1$-Lipschitz path $\gamma:[s,t]\to[0,\infty)$ with
$\gamma(s)=x$, $\gamma(t)=y$, and is $-\infty$ if $|x-y|>t-s$.
\end{definition}
For passage from time $s$ to time $t$, split an admissible path at its
position $z$ at time $r$. Any atom at time $r$ is counted in the first
segment and excluded from the second. Maximizing over the segments and their meeting position gives
\begin{equation}
 \mathcal L_{\mathrm{PNG}}^{(c)}(x,s;y,t)
 =\max_{z\geq0}\{
   \mathcal L_{\mathrm{PNG}}^{(c)}(x,s;z,r)
   +\mathcal L_{\mathrm{PNG}}^{(c)}(z,r;y,t)\},
 \qquad s<r<t.
 \label{eq:half-space-png-random-metric-composition}
\end{equation}
Define
\begin{equation*}
 L_\varepsilon^\rho(x,s;y,t)
 =\varepsilon^{1/2}\mathcal L_{\mathrm{PNG}}^{(c_\varepsilon(\rho))}
  (\varepsilon^{-1}x,\varepsilon^{-3/2}s;
   \varepsilon^{-1}y,\varepsilon^{-3/2}t)-2\varepsilon^{-1}(t-s).
\end{equation*}

The max-plus representation gives the exact identity
\begin{equation*}
 \mathfrak h_\varepsilon^\rho(t,y)
 =\sup_{x\geq0}\{f_\varepsilon(x)
       +L_\varepsilon^\rho(x,0;y,t)\}.
\end{equation*}

\begin{theorem}[Half-space directed landscape limit]
\label{thm:intro-half-space-png-directed-landscape-limit}
For every $\rho\in\mathbb R\cup\{-\infty\}$, there are continuous
mesh interpolations $\widetilde L_\varepsilon^\rho$ such that
\begin{equation*}
 \widetilde L_\varepsilon^\rho\Longrightarrow\mathcal L^\rho
 \quad\text{in }C_{\mathrm{loc}}(\mathbb H_\uparrow^2),\qquad
 \|\widetilde L_\varepsilon^\rho-L_\varepsilon^\rho\|_K
       \longrightarrow0\quad\text{in probability}
\end{equation*}
for every compact $K\subset\mathbb H_\uparrow^2$.
Thus the raw metrics converge locally uniformly in the convention above.
For $\rho=-\infty$, this holds for each fixed $c_\star\in[0,1)$.
\end{theorem}
The centering retains the parabola $-(x-y)^2/(t-s)$ in the limit.
Section~\ref{subsec:directed-landscape-convergence} proves this
theorem using the finite-wedge specialization of the transition
limits in Section~\ref{subsec:kpz-fixed-point-convergence}.
The half-space KPZ fixed point has the variational representation
\cite[Theorem~4.11]{dauvergne2026directed}:
\begin{equation}
 \mathfrak h^\rho(t,y;f)
 =\sup_{x\geq0}\{f(x)+\mathcal L^\rho(x,0;y,t)\}.
 \label{eq:intro-half-space-kpz-variational}
\end{equation}
Corollary~\ref{cor:png-general-initial-transfer} gives joint convergence
of the metric and finitely many height evolutions driven by the same
Poisson environment. Every limiting height field is given by this
variational formula with the same $\mathcal L^\rho$.

\subsection{Proof strategy and organization}\label{subsec:outline}
Section~\ref{sec:stationary-measures} encodes the one-profile law by a
reflected spatial walk and proves generator balance at the wall.
The joint construction combines geometric last-passage stationarity
and symmetries from \cite{dauvergne2026directed} with a finite line
rearrangement and a limit retaining one common future. The singular
endpoint uses a shifted source, and planar crossing gives the regular
PNG horizon. Section~\ref{subsec:stationary-horizon-convergence} proves
its limit by an invariance principle for the lines and bridges and
convergence of the representing Stieltjes measures.

Section~\ref{subsec:kpz-fixed-point-convergence} proves uniform
boundary and hitting estimates for the exact transition Pfaffian
\cite[Theorem~1.1]{chenZhangPNG}. Trace and Hilbert--Schmidt bounds
justify its limit to Zhang's transition laws \cite{zhang2024tasep}. Fixed subcritical rates
require a skew derivative form and stationary comparison at the wall,
where an individual kernel block has only a distributional limit.
Localization, profile tightness, Markov iteration, and short-time
bounds give the height-field limit.

Section~\ref{subsec:directed-landscape-convergence} uses finite-wedge
tests and \cite[Theorem~5.1 and Corollary~5.2]{dauvergne2026directed}
to identify the rational metric law, then proves four-variable
tightness. The exact variational formula and localization identify
the joint height limits as functions of that same metric.

\phantomsection
\subsection*{AI assistance}
The main ideas and overall research direction of this work originated
with the author. The author acknowledges substantial assistance from
OpenAI's ChatGPT in research discussions, proof development, and
manuscript preparation, and takes full responsibility for the
mathematical arguments and the final manuscript.

\subsection*{Acknowledgements}\label{subsec:acknowledgements}
The author thanks the organizers of the Virginia Integrable Probability
Summer School 2024 and the workshop ``All roads to the KPZ universality
class'' at the American Institute of Mathematics for providing
stimulating research environments. The author is grateful to Xincheng
Zhang and Lingfu Zhang for helpful discussions.

\section{Stationary measures}
\label{sec:stationary-measures}

We first prove stationarity of the one-parameter laws $\nu_{c,g}$
by a generator calculation in Section~\ref{subsec:stationary-measures}.
Section~\ref{subsec:joint-stationary-measures} constructs the laws
$\nu_{c;g_1,\ldots,g_k}$ and proves their joint stationarity under
a common graphical environment, together with parameter-deletion
consistency. Finally, Section~\ref{subsec:half-space-png-horizon}
assembles this family into a canonical measurable field, the
half-space PNG horizon.

\subsection{Stationarity at one parameter}
\label{subsec:stationary-measures}

We prove Theorem~\ref{thm:intro-half-space-png-stationary-family}
by checking generator balance for a marked reflected random-walk representation of
$\nu_{c,g}$.

\subsubsection{Step coordinates and a local stationarity criterion}
\label{subsubsec:step-coordinates}

For bounded measurable \(F\) on \(\mathrm{UC}_{+,0}\), set
\begin{equation*}
 (\mathsf S_t^{(c)}F)(h)
 :=\mathbb E_h\!\left[
 F\!\left(\operatorname{Rec}
   \bigl(h^{(c)}(t,\cdot;h)\bigr)\right)\right].
\end{equation*}
The graphical cocycle, independent time slabs, and vertical-shift
covariance give \(\mathsf S_0^{(c)}=I\) and
\(\mathsf S_{s+t}^{(c)}=\mathsf S_s^{(c)}\mathsf S_t^{(c)}\).

For the index set of a finite or infinite step configuration, write
\begin{equation*}
 I_0=\varnothing,\qquad
 I_N=\{1,\ldots,N\}\quad(1\leq N<\infty),\qquad
 I_\infty=\mathbb N.
\end{equation*}
\begin{definition}[Step-coordinate state space]
\label{def:half-space-png-step-coordinate-space}
A profile belongs to
\(\mathcal X_{\mathrm{step}}\subset\mathrm{UC}_{+}\) if it has a
representation \(h=(m;\xi)\), where \(m\in\mathbb Z\) and
\(\xi=((y_i,\sigma_i))_{i\in I_N}\), satisfying
\begin{equation*}
 \begin{gathered}
 y_i\geq0,\quad \sigma_i\in\{+1,-1\},\quad y_i\leq y_{i+1},
 \qquad \#\{i:y_i\leq R\}<\infty\quad(R<\infty),\\
 y_i=y_{i+1}>0\ \Longrightarrow\
 (\sigma_i,\sigma_{i+1})=(+1,-1),
 \qquad
 y_i=0\ \Longrightarrow\ i=1\text{ and }\sigma_i=-1,
\end{gathered}
\end{equation*}
and
\begin{equation*}
 h(x)=m+\sum_{y_i<x}\sigma_i
       +\#\{i:y_i=x,\ \sigma_i=+1\},\qquad x\geq0.
\end{equation*}
\end{definition}
Thus \(m=h(0)\), and reconstruction takes the upper value at each step,
including a coincident \((+1,-1)\)-spike. The representation is unique. Put
\(\mathcal X_{\mathrm{step},0}:=\mathcal X_{\mathrm{step}}\cap\mathrm{UC}_{+,0}\).
Reconstruction is a Borel injection from the standard Borel space of
admissible marked configurations; its image and inverse are therefore
Borel. Thus $\mathcal X_{\mathrm{step},0}$ is standard Borel.

\begin{definition}[Regular configuration]
\label{def:half-space-png-regular-configuration}
A configuration in $\mathcal X_{\mathrm{step}}$ is \emph{regular}
if its step locations are distinct and strictly positive.
\end{definition}

Write \(\xi_{<z}\) and \(\xi_{>z}\) for the subconfigurations with
locations below and above \(z\); the symbols \(\xi^-\) and \(\xi^+\)
denote the ordered configurations preceding and following a displayed
block of marks. Between events, $\dot y_i=-\sigma_i$. At locations
without an existing step, nucleations act by
\begin{equation*}
 C_z(m;\xi)=(m;\xi_{<z},(z,+1),(z,-1),\xi_{>z}),
 \qquad C_0(m;\xi)=(m+1;(0,-1),\xi).
\end{equation*}
An adjacent \((-1,+1)\)-pair is deleted on contact, and an incoming
up-step is absorbed at the wall:
\begin{equation*}
 (m;\xi^-,(q,-1),(q,+1),\xi^+)\longmapsto(m;\xi^-,\xi^+),
 \qquad
 (m;(y,+1),\xi)\xrightarrow[y\downarrow0]{}(m+1;\xi).
\end{equation*}
The incoming collision and absorption faces are identified with the
configurations after deletion. For each fixed initial state, these coordinate rules
almost surely reconstruct graphical PNG and
preserve $\mathcal X_{\mathrm{step}}$, and, after recentering,
$\mathcal X_{\mathrm{step},0}$. Indeed, the level sets
$A_a(h):=\{x\geq0:h(x)\geq a\}$ satisfy
\begin{equation*}
 A_a(\mathsf D_sh)
 =\{x\geq0:\operatorname{dist}_{\mathbb R_+}(x,A_a(h))\leq s\},
 \qquad a\in\mathbb Z.
\end{equation*}
Their endpoints have the stated velocities and deletion rules; isolated
components produce outgoing $(+1,-1)$ pairs. Poisson atoms almost surely
avoid existing steps, so the creation rules agree with one-point raisings.
Induction over the finitely many events in each compact dependence cone,
followed by countable exhaustion, proves agreement.

Lift observables on \(\mathcal X_{\mathrm{step},0}\) to vertically shift-invariant
functions on $\mathcal X_{\mathrm{step}}$. Let $\mathscr D_{\mathrm{pol}}$
consist of bounded observables $F$ admitting finite radii $0<R_0<R$
such that $F(h)$ depends only on $h|_{[0,R_0]}$. For
\(\mathbf y=(y_1,\ldots,y_n)\) and
\(\boldsymbol\sigma=(\sigma_1,\ldots,\sigma_n)\), let
\(F_n(m;\mathbf y,\boldsymbol\sigma)\) be its coordinate form when
there are \(n\) marks in \([0,R)\); locality makes the continuation
beyond \(R\) irrelevant. On each ordered simplex with fixed signs,
these forms are locally Lipschitz and
piecewise $C^1$, independent of $m$, and satisfy the boundary matching conditions
\begin{align*}
 \lim_{\varepsilon\downarrow0}
 F_n(m;\xi^-,(q,-1),(q+\varepsilon,+1),\xi^+)
 &=F_{n-2}(m;\xi^-,\xi^+), &&q>0,
\\
 \lim_{y\downarrow0}F_n(m;(y,+1),\xi)
 &=F_{n-1}(m+1;\xi).
\end{align*}
At bulk and boundary creation faces, $F(C_zh)$ and $F(C_0h)$ are
the traces from the outgoing configurations. Require, almost everywhere,
\begin{equation}
 \sum_{i:y_i<R}|\partial_{y_i}F(h)|
 \leq C(1+N_R(h))^p,
 \qquad N_R(h):=\#\{i:y_i\leq R\},
 \label{eq:stationary-polynomial-test-bound}
\end{equation}
for some $C<\infty$ and integer $p\geq0$. On cells of regular configurations, define
almost everywhere,
$
 \mathcal TF(h):=-\sum_i\sigma_i\partial_{y_i}F_n,
$
and the local generator 
\begin{equation*}
 (\mathscr G_cF)(h)=\mathcal TF(h)
 +2\int_0^\infty\!\bigl(F(C_zh)-F(h)\bigr)\,\mathrm dz
 +c\bigl(F(C_0h)-F(h)\bigr).
\end{equation*}
Set the bulk integrand to zero at existing step locations, a Lebesgue-null
set. The integral vanishes beyond $R_0$; the matching conditions account
for deterministic deletions.

\begin{lemma}[A local generator criterion]
\label{lem:stationary-echeverria}
Fix $c\geq0$. Let $\mu$ be a probability measure supported on regular
configurations in $\mathcal X_{\mathrm{step},0}$. Assume that its
finite-window coordinate distributions are absolutely continuous on
each fixed-count, fixed-sign ordered simplex, and that
\begin{equation*}
 \int N_R(h)^p\,\mu(dh)<\infty,
 \qquad 0<R<\infty,\quad p\in\mathbb N.
\end{equation*}
If
\begin{equation*}
 \int\mathscr G_cF\,d\mu=0,
 \qquad F\in\mathscr D_{\mathrm{pol}},
\end{equation*}
then $\mu\mathsf S_t^{(c)}=\mu$ for every $t\geq0$.
\end{lemma}

\begin{proof}
\emph{Step 1: A determining family.}
For $h=(0;\xi)$, write $\eta_h=\sum_i\sigma_i\delta_{y_i}$ and consider
\begin{equation*}
 F(h)=\Phi\bigl(\langle\eta_h,\varphi_1\rangle,\ldots,
                    \langle\eta_h,\varphi_j\rangle\bigr),
 \qquad
 \varphi_a\in C_c^\infty((0,\infty)),\quad
 \Phi\in C_b^\infty(\mathbb R^j).
\end{equation*}
These functions determine laws on regular pinned profiles: the smooth
spatial tests determine $\eta_h$, and pinning fixes the additive height
constant. They belong to $\mathscr D_{\mathrm{pol}}$, since opposite marks
cancel at collisions and the spatial tests vanish near the wall.
Their total coordinate gradient is bounded by $a_FN_R$, where
\begin{equation*}
 a_F:=\sum_{a=1}^j
 \|\partial_a\Phi\|_\infty\|\varphi_a'\|_\infty
\end{equation*}
and $R$ contains all spatial supports.

\smallskip
\emph{Step 2: Truncated backward observables.}
Fix $T<\infty$ and $R_{\mathrm{src}}>R+T+1$.
Let $N_t^{\mathrm{src}}$ count the bulk and wall sources in
$(0,t]\times(0,R_{\mathrm{src}})$ and $(0,t]$, respectively, and put
$\Lambda=2R_{\mathrm{src}}+c$. With $h_t=h^{(c)}(t,\cdot;h)$, define
\begin{equation*}
 U_k(t,h):=\mathbb E_h\!\left[
 F(\operatorname{Rec}h_t)\mathbf1_{\{N_t^{\mathrm{src}}=k\}}\right],
 \qquad V_M(t,h):=\sum_{k=0}^M U_k(t,h).
\end{equation*}
We verify directly that $V_M(t,\cdot)\in\mathscr D_{\mathrm{pol}}$
for $0\leq t\leq T$.

For fixed source types and $k\leq M$, the ordered source coordinates
range over a bounded box of dimension at most $2M$.
With $n=N_{R_{\mathrm{src}}}(h)$ initial marks, all candidate steps
follow affine rays. Spatial insertion order and comparisons of births,
meetings, wall arrivals, and the observation time partition this box
into history cells using at most $C_M(n+1)^4$ affine hyperplanes.
Inside each cell the terminal observable is smooth, with total
initial-coordinate gradient bounded by $C_F(n+2M)$.
Interfaces involving source coordinates have nonzero source normals
from a finite set depending only on $M$ and uniformly bounded
displacement speeds in the initial coordinates. Under an
initial-coordinate displacement, their swept volumes are bounded
by a polynomial in $n$ times that displacement, and the integrand
jumps by at most $2\|F\|_\infty$.
Across interfaces independent of the source coordinates, the observable
is continuous by collision cancellation and vanishing of the spatial
tests near the wall. Integrating over the source coordinates gives
local Lipschitz continuity and
\begin{equation*}
 \sum_{i:y_i<R_{\mathrm{src}}}
 |\partial_{y_i}V_M(t,h)|
 \leq C_{F,R_{\mathrm{src}},T,M}
       (1+N_{R_{\mathrm{src}}}(h))^{p_M}.
\end{equation*}
On regions of fixed combinatorial type, the history polytopes have
affine vertices; triangulation and differentiation under the integral
give piecewise $C^1$ regularity.

The boundary matching conditions follow by bounded convergence.
Indeed, for almost every fixed source list, a collapsing incoming pair,
or an up-step approaching zero, disappears before the first source and
the observation time. The subsequent evolutions then agree with those
started after deletion. Outgoing creation traces similarly agree with
the graphical creation states. At $t=0$ the matching holds for $F$.
Finally, finite propagation gives locality in $[0,R+T]$, proving the
claimed domain membership.

\smallskip
\emph{Step 3: The truncated backward equation.}
Write
\begin{equation*}
 (\mathscr N f)(h):=
 2\int_0^{R_{\mathrm{src}}}f(C_zh)\,dz+cf(C_0h).
\end{equation*}
Differentiating the ordered source integrals gives
\begin{equation*}
 \partial_tU_k=(\mathcal T-\Lambda)U_k+\mathscr N U_{k-1},
 \qquad U_{-1}=0.
\end{equation*}
By locality, $\mathscr G_c=\mathcal T-\Lambda+\mathscr N$ on these
observables. Summing over $k\leq M$ yields, coordinate-almost everywhere,
\begin{equation}
 \partial_tV_M=\mathscr G_cV_M-\mathscr N U_M,
 \qquad V_M(0)=F.
 \label{eq:stationary-truncated-generator-evolution}
\end{equation}
The absolute-continuity assumption makes this identity valid
$\mu$-almost everywhere. Moreover, $F$ is continuous through graphical
events, and between events its speed is at most
$a_F(N_{R_{\mathrm{src}}}(h)+2N_T^{\mathrm{src}})$.
The truncation indicator changes only at a source, so
\begin{equation*}
 |V_M(t,h)-V_M(s,h)|
 \leq\bigl[a_F(N_{R_{\mathrm{src}}}(h)+2\Lambda T)
       +\Lambda\|F\|_\infty\bigr]|t-s|.
\end{equation*}
This integrable bound and the polynomial gradient bound justify
time integration and Fubini against $\mu$.

\smallskip
\emph{Step 4: Generator balance and removal of the truncation.}
Since $V_M(t,\cdot)\in\mathscr D_{\mathrm{pol}}$, the balance hypothesis
and \eqref{eq:stationary-truncated-generator-evolution} give
\begin{equation*}
 \int V_M(t,\cdot)\,d\mu-\int F\,d\mu
 =-\int_0^t\int\mathscr N U_M(s,\cdot)\,d\mu\,ds.
\end{equation*}
The bound $|U_M(s,h)|\leq
\|F\|_\infty e^{-\Lambda s}(\Lambda s)^M/M!$ implies
\begin{equation*}
 \left|\int V_M(t,\cdot)\,d\mu-\int F\,d\mu\right|
 \leq\|F\|_\infty
 \mathbb P\bigl(\operatorname{Pois}(\Lambda t)>M\bigr).
\end{equation*}
Also, uniformly in $h$,
\begin{equation*}
 |\mathsf S_t^{(c)}F(h)-V_M(t,h)|
 \leq\|F\|_\infty
 \mathbb P\bigl(\operatorname{Pois}(\Lambda t)>M\bigr).
\end{equation*}
Letting $M\to\infty$ proves
$\int\mathsf S_t^{(c)}F\,d\mu=\int F\,d\mu$.
At each fixed deterministic time the evolved configuration is regular
almost surely: incoming collision and wall steps are deleted immediately,
and a Poisson source occurs at that time with probability zero.
The determining family therefore identifies the evolved law with $\mu$.
Since $T$ was arbitrary, the conclusion holds for every $t\geq0$.
\end{proof}

\subsubsection{A reflected random-walk representation}
\label{subsubsec:reflected-walk}

We encode $\nu_{c,g}$ by a reflected walk whose clock marks record
the visible steps of the height profile. This gives the finite-window
densities used in the generator calculation below.

\paragraph*{\textbf{The reflected random-walk}}
Fix $g\in\mathcal G_c$ with $c<g$, and abbreviate
$S_j=S_j^{(g)}$, $Y=Y^{(c,g)}$,
$\widetilde H=\widetilde H^{(c,g)}$, and $H=H^{(c,g)}$.
Set
\begin{equation*}
 W(x):=S_1(x)-S_2(x),\qquad
 M(x):=\sup_{0\leq y\leq x}W(y),\qquad
 Q(x):=\max\{Y,M(x)\}-W(x).
\end{equation*}
Thus $Q$ is the reflection of $Y-W$ at zero, and
\begin{equation*}
 Q(0)=Y,\qquad
 \widetilde H(x)=S_1(x)+Q(x)-Y,\qquad
 \pi_q:=\mathbb P(Q(0)=q)=(1-c/g)(c/g)^q.
\end{equation*}
At $c=0$, the last expression means $\boldsymbol\pi=\delta_0$.
We regard $\boldsymbol\pi=(\pi_q)_{q\geq0}$ as a row vector.
It specifies the initial law of this auxiliary spatial process.
Upper-semicontinuous regularization preserves the locations and signs
of the jumps, so the same step marks describe $H$.

Call a clock event \emph{visible} if it changes $\widetilde H$.
For the up- and down-jump clocks $N_{j,+},N_{j,-}$ of $S_j$,
the identity above gives
\begin{equation*}
 \begin{array}{c|c|c|c}
 \text{clock}&\text{rate}&Q:q\longmapsto q'&
 \Delta\widetilde H\\ \hline
 N_{1,+}&g&(q-1)_+&\mathbf1_{\{q=0\}}\\
 N_{1,-}&g^{-1}&q+1&0\\
 N_{2,+}&g^{-1}&q+1&1\\
 N_{2,-}&g&(q-1)_+&-\mathbf1_{\{q>0\}}.
 \end{array}
\end{equation*}
In particular, the down-step intensity at the spatial origin is
$g\mathbb P(Q(0)>0)=c$.
The clock marks must be retained: at $q=0$, both $N_{1,+}$ and
$N_{2,-}$ leave $Q$ unchanged, but only $N_{1,+}$ is visible.

\paragraph*{\textbf{Visible-step matrices}}
Let $I$ be the identity matrix on $\mathbb N_0$, and let $\mathbf1$
be the column vector of ones. The rate matrices for visible up-steps
and down-steps are, respectively,
\begin{equation*}
 (\mathsf Q_+)_{qp}
 =g\mathbf1_{\{q=p=0\}}+g^{-1}\mathbf1_{\{p=q+1\}},
 \qquad
 (\mathsf Q_-)_{qp}
 =g\mathbf1_{\{q\geq1,\ p=q-1\}}.
\end{equation*}
Put $r=g+g^{-1}$ and $\mathsf A=\mathsf Q_++\mathsf Q_-$.
The table shows that the invisible events have the same total
hidden-state rate matrix $\mathsf A$. Since the four clocks have
total rate $2r$, evolution with no visible step is described by
\begin{equation*}
 \mathsf Q_0:=\mathsf A-2rI,
 \qquad \mathsf E(u):=e^{u\mathsf Q_0}.
\end{equation*}
The entry $\mathsf E(u)_{qp}$ is the probability of reaching state
$p$ from $q$ over spatial length $u$ without a visible step.
Allowing visible steps as well gives the reflected-walk generator
\begin{equation*}
 \mathsf Q=\mathsf Q_0+\mathsf Q_++\mathsf Q_-
          =2(\mathsf A-rI).
\end{equation*}

These matrices act boundedly on $\ell^\infty(\mathbb N_0)$, paired
with row vectors in $\ell^1(\mathbb N_0)$.
Since $\mathsf A\mathbf1=r\mathbf1$,
\begin{equation*}
 \mathsf E(u)=e^{-2ru}\sum_{m\geq0}\frac{u^m\mathsf A^m}{m!},
 \qquad \mathsf E(u)\mathbf1=e^{-ru}\mathbf1.
\end{equation*}
Thus $\mathsf E(u)$ is positive with operator norm $e^{-ru}$,
and $\mathsf Q$ is a bounded conservative Markov generator.
The following identities will be used in the generator balance:
\begin{equation*}
 \mathsf Q_+\mathsf Q_-=I,\qquad
 \boldsymbol\pi\mathsf Q_-=c\boldsymbol\pi,\qquad
 \mathfrak C:=\mathsf Q_-\mathsf Q_+-I.
\end{equation*}

\paragraph*{\textbf{Finite-window jump densities}}
For $L>0$, let $\Delta_n(L)=\{0<y_1<\cdots<y_n<L\}$,
with the empty configuration when $n=0$, and write
$\mathsf Q_{\pm1}=\mathsf Q_\pm$.
For prescribed ordered jump locations $\mathbf y$ and signs
$\boldsymbol\sigma$, the density is
\begin{equation*}
 \mathcal J_L(\mathbf y,\boldsymbol\sigma)
 :=\boldsymbol\pi\mathsf E(y_1)\mathsf Q_{\sigma_1}
 \mathsf E(y_2-y_1)\cdots
 \mathsf Q_{\sigma_n}\mathsf E(L-y_n)\mathbf1.
\end{equation*}
Each $\mathsf E$ factor accounts for a gap with no visible step;
each $\mathsf Q_{\sigma_i}$ inserts a step of the specified sign;
the final multiplication by $\mathbf1$ sums over the terminal hidden
state. For $n=0$, the product is $\boldsymbol\pi\mathsf E(L)\mathbf1$.
This is the ordered marked Janossy density, meaning the joint density
of exactly these steps and no others in $(0,L)$.

Summing over signs replaces each $\mathsf Q_{\sigma_i}$ by $\mathsf A$.
The identities above give
\begin{equation*}
 \sum_{\boldsymbol\sigma\in\{\pm1\}^n}
 \mathcal J_L(\mathbf y,\boldsymbol\sigma)=r^ne^{-rL}.
\end{equation*}
Since $\Delta_n(L)$ has volume $L^n/n!$, these densities have total
mass one. They also show that the unmarked step locations form a
homogeneous Poisson process of rate $r$; in particular,
$N_L\sim\operatorname{Pois}(rL)$.

There is almost surely no step at a fixed endpoint $L$. Hence, for
bounded measurable $F$ depending only on $h|_{[0,L]}$,
\begin{equation*}
 \int F(h)\,\nu_{c,g}(dh)
 =\sum_{n\geq0}\sum_{\boldsymbol\sigma\in\{\pm1\}^n}
 \int_{\Delta_n(L)}
 \mathcal J_L(\mathbf y,\boldsymbol\sigma)
 F_n(\mathbf y,\boldsymbol\sigma)\,d\mathbf y.
\end{equation*}
Here $F_n$ is evaluated on the pinned upper-semicontinuous profile
with the indicated step marks. Thus $\nu_{c,g}$ has regular support,
absolutely continuous finite-window coordinate laws, and finite moments
of all orders for every local step count, as required by
Lemma~\ref{lem:stationary-echeverria}.

\begin{remark}[Relation to full-space PNG stationary measures]
\label{rem:half-space-png-mqr-invariant-measures}
For $a>0$, let $\mu_a^+$ be the restriction to $[0,\infty)$, pinned
at zero, of the full-space PNG stationary law with independent
Poisson up-steps and down-steps of rates $a$ and $a^{-1}$
\cite[Section~4.3, Proposition~4.2]{matetski2025polynuclear}.
For every $c>0$,
\begin{equation*}
 \mu_{1/c}^+=
 \begin{cases}
  \nu_{c,1/c},&0<c<1,\\
  \nu_{c,c},&c\geq1.
 \end{cases}
\end{equation*}
For $c\geq1$, this follows directly from
\eqref{eq:intro-half-space-png-boundary-dominated-profile}, whose
up/down rates are $(c^{-1},c)$.
For $0<c<1$, take $g=1/c$. Then
$\pi_q=(1-g^{-2})g^{-2q}$, and direct multiplication gives
$\boldsymbol\pi\mathsf Q_+=g\boldsymbol\pi$,
 $\boldsymbol\pi\mathsf Q_-=g^{-1}\boldsymbol\pi$,
$ \boldsymbol\pi\mathsf E(u)=e^{-ru}\boldsymbol\pi$.
Consequently,
$\mathcal J_L(\mathbf y,\boldsymbol\sigma)
 =e^{-rL}\prod_{i=1}^n g^{\sigma_i}$,
which is the density of independent Poisson up-steps and down-steps
with rates $(g,g^{-1})=(c^{-1},c)$.
For general $g$, the signs need not be independent, although the
unmarked step locations remain Poisson.
\end{remark}

\subsubsection{Generator balance}
\label{subsubsec:generator-balance}

Fix \(F\in\mathscr D_{\mathrm{pol}}\), with locality radii
\(0<R_0<R\), and take the integration cutoff \(R_0<L<R\).
Locality gives
\begin{equation*}
 2\int_0^\infty(F(C_zh)-F(h))\,\mathrm dz
 =2\int_0^L F(C_zh)\,\mathrm dz-2L F(h).
\end{equation*}
Abbreviate
\begin{equation*}
 \sum{}^{\!*}G
 :=\sum_{n\geq0}\sum_{\boldsymbol\sigma\in\{\pm1\}^n}
 \int_{\Delta_n(L)}G_n(\mathbf y,\boldsymbol\sigma)\,\mathrm d\mathbf y.
\end{equation*}
Write \(\mathcal J_L^{(i)}[\mathsf M]\) for the matrix product with
\(\mathsf Q_{\sigma_i}\) replaced by \(\mathsf M\). For \(z\in(y_k,y_{k+1})\),
with \(y_0=0,y_{n+1}=L\), write \(\mathcal J_{L,z}[\mathsf M]\) for the
matrix product with its gap factor replaced by \(\mathsf E(z-y_k)\mathsf M\mathsf E(y_{k+1}-z)\), using
one-sided endpoint conventions. The mean generator comprises
\begin{align*}
 \mathcal T_F&=-\sum{}^{\!*}\mathcal J_L
   \sum_{i=1}^n\sigma_i\partial_{y_i}F_n,
 &\mathcal T_{\mathrm{bulk}}^{+}&=2\sum{}^{\!*}\mathcal J_L\int_0^L F(C_zh)\,\mathrm dz,
 \notag\\
 \mathcal T_{\mathrm{bulk}}^{-}&=-2L\sum{}^{\!*}\mathcal J_LF_n,
 &\mathcal C_+&=c\sum{}^{\!*}\mathcal J_LF(C_0h),
 \qquad \mathcal C_-=-c\sum{}^{\!*}\mathcal J_LF_n.
\end{align*}

Set \(\mathfrak Z=\mathsf Q_+-\mathsf Q_-\), and write \([\mathsf U,\mathsf V]:=\mathsf U\mathsf V-\mathsf V\mathsf U\) for the
commutator. Direct multiplication also gives
\begin{equation*}
 [\mathsf Q_0,\mathfrak Z]=2\mathfrak C,\qquad
 [\mathfrak Z,\mathsf Q_\sigma]=-\mathfrak C,\qquad
 \sigma\,\partial_y\!\left(\mathsf E(y-u)\mathsf Q_\sigma \mathsf E(v-y)\right)
 =\mathsf E(y-u)\mathfrak C \mathsf E(v-y).
\end{equation*}
Integration by parts yields
\begin{equation}
 -\int_{\Delta_n(L)}\!\mathcal J_L\sigma_i\partial_{y_i}F_n
 =\int_{\Delta_n(L)}\!\mathcal J_L^{(i)}[\mathfrak C]F_n
  +\int_{\partial_i^-}\!\sigma_i\mathcal J_LF_n
  -\int_{\partial_i^+}\!\sigma_i\mathcal J_LF_n,
\label{eq:stationary-integration-by-parts}
\end{equation}
where \(\partial_i^-\) is the face \(y_i=y_{i-1}\) and
\(\partial_i^+\) is the face \(y_i=y_{i+1}\), with
\(y_0=0\) and \(y_{n+1}=L\). Face integrals use the remaining
\(n-1\) ordered coordinates.
Internal-face coefficients \(\sigma_{i+1}-\sigma_i\) are zero for
equal signs, \(-2\) for \((+1,-1)\), and \(+2\) for
\((-1,+1)\).

The creation face has observable \(F(C_zh)\) and factor
\(\mathsf Q_+\mathsf Q_-=I\), so cancels \(\mathcal T_{\mathrm{bulk}}^{+}\). At the annihilation
face, the matching condition deletes the pair. Its \(\mathsf Q_-\mathsf Q_+\) factor and the
\(-2I\) insertion representing \(\mathcal T_{\mathrm{bulk}}^{-}\) leave
\(2\sum{}^{\!*}\int_0^L\mathcal J_{L,z}[\mathfrak C]F_n\,\mathrm dz\).

Within gaps,
\(\partial_z\mathcal J_{L,z}[\mathfrak Z]
=2\mathcal J_{L,z}[\mathfrak C]\).
At \(y_i\), its left trace minus right trace is
\(\mathcal J_L^{(i)}[[\mathfrak Z,\mathsf Q_{\sigma_i}]]
=-\mathcal J_L^{(i)}[\mathfrak C]\). Telescoping therefore gives
\begin{equation*}
 2\int_0^L\mathcal J_{L,z}[\mathfrak C]\,\mathrm dz
 =\mathcal J_{L,L}[\mathfrak Z]-\mathcal J_{L,0}[\mathfrak Z]
  -\sum_{i=1}^n\mathcal J_L^{(i)}[\mathfrak C].
\end{equation*}
After multiplication by \(F_n\) and summation, the last term cancels
the density derivative in \eqref{eq:stationary-integration-by-parts}.
Since \(F\) does not depend on steps near \(L\), the boundary contributions there equal
\(-\mathcal J_{L,L}[\mathsf Q_+]F+\mathcal J_{L,L}[\mathsf Q_-]F
=-\mathcal J_{L,L}[\mathfrak Z]F\), canceling the right endpoint.

At zero, shift invariance and absorption matching give
\begin{equation*}
 F(0;(0,-1),\xi)=F(C_0(0;\xi)),\qquad
 \lim_{y\downarrow0}F(0;(y,+1),\xi)=F(0;\xi).
\end{equation*}
For the remaining column vector
\begin{equation*}
 \boldsymbol v_{\mathrm{rem}}:=\mathsf E(y_1)\mathsf Q_{\sigma_1}\mathsf E(y_2-y_1)\cdots
 \mathsf Q_{\sigma_n}\mathsf E(L-y_n)\mathbf1,
\end{equation*}
with \(\boldsymbol v_{\mathrm{rem}}=\mathsf E(L)\mathbf1\) when \(n=0\), the boundary
contributions at zero combine as
\begin{align*}
 &-\boldsymbol\pi \mathsf Q_-\boldsymbol v_{\mathrm{rem}} F(C_0h)
   +\boldsymbol\pi \mathsf Q_+\boldsymbol v_{\mathrm{rem}} F(h)
   -\boldsymbol\pi(\mathsf Q_+-\mathsf Q_-)\boldsymbol v_{\mathrm{rem}} F(h)\notag\\
 &\hspace{25mm}
 =\boldsymbol\pi \mathsf Q_-\boldsymbol v_{\mathrm{rem}}\bigl(F(h)-F(C_0h)\bigr)
 =c\boldsymbol\pi\boldsymbol v_{\mathrm{rem}}\bigl(F(h)-F(C_0h)\bigr).
\end{align*}
Summing over the remaining products cancels \(\mathcal C_++\mathcal C_-\).

By locality, deleting the tail beyond $L$ bounds the observable derivatives
by $C(1+N_L)^p$, integrable under the Poisson step counts.
If $e_q$ denotes the $q$th coordinate vector, then
$\mathfrak C=(g^2e_1-e_0)e_0^{\mathsf T}$.
Positivity of $\mathsf E$, its row sums, and
$\mathsf A\mathbf1=r\mathbf1$ give, for every nonnegative bounded matrix
$\mathsf D$,
\begin{align*}
 \sum{}^{\!*}\int_0^L\mathcal J_{L,z}[\mathsf D] \,\mathrm dz
 &\leq L\|\mathsf D\|_{\infty\to\infty},\\
 \sum_{n\geq1}\sum_{\boldsymbol\sigma}
 \int_{\Delta_n(L)}\sum_{i=1}^n
       |\mathcal J_L^{(i)}[\mathfrak C]|\,\mathrm d\mathbf y
 &\leq2L\|\lvert\mathfrak C\rvert\|_{\infty\to\infty}.
\end{align*}
Here $\lvert\mathfrak C\rvert$ is taken entrywise; the factor two accounts
for the two signs of the replaced mark. Boundedness of $\mathsf Q_\sigma$ and
$\mathsf Q_\sigma\mathsf Q_\tau$ similarly controls the reindexed face terms.
Thus Fubini, integration by parts, and face reindexing are justified
absolutely. Continuity cancels boundary terms between the $C^1$ regions. Hence
\begin{equation}
 \int(\mathscr G_cF)(h)\,\nu_{c,g}(\mathrm dh)=0,
 \qquad F\in\mathscr D_{\mathrm{pol}}.
 \label{eq:stationary-generator-balance-final}
\end{equation}

For the endpoint \(g=c\geq1\), the right-continuous profile is \(S_2^{(c)}\),
whose regularization has the same signed steps and scalar Janossy data
\(\mathsf Q_+=c^{-1}\), \(\mathsf Q_-=c\), \(\mathsf Q_0=-(c+c^{-1})\), and
\(\boldsymbol\pi=\mathbf1=1\).
Here \(\mathsf Q_+\mathsf Q_-=1\), \(\boldsymbol\pi \mathsf Q_-=c\boldsymbol\pi\), and \(\mathfrak C=0\);
the same face calculation proves
\eqref{eq:stationary-generator-balance-final}.

\begin{proof}[Proof of Theorem
\ref{thm:intro-half-space-png-stationary-family}]
For every \(g\in\mathcal G_c\), \(\nu_{c,g}\) has regular support
(Definition~\ref{def:half-space-png-regular-configuration}),
absolutely continuous finite-window Janossy densities, and Poisson local
step counts. The balance and Lemma \ref{lem:stationary-echeverria} give
\(\nu_{c,g}\mathsf S_t^{(c)}=\nu_{c,g}\).

\end{proof}

\subsection{Jointly stationary measures}
\label{subsec:joint-stationary-measures}

We prove Theorem~\ref{thm:intro-half-space-png-joint-stationary-family}.
For \(\boldsymbol h=(h_1,\ldots,h_k)\in(\mathrm{UC}_{+,0})^k\) and bounded measurable \(F\), write
\begin{align*}
 \Theta_t^{(c),[k]}(\boldsymbol h;\Pi)
 :=\left(\operatorname{Rec}
      \bigl(h^{(c)}(t,\cdot;h_i;\Pi)\bigr)\right)_{i=1}^k,\qquad
 (\overline{\mathsf P}_t^{(c),[k]}F)(\boldsymbol h)
 :=\mathbb E_\Pi F\bigl(\Theta_t^{(c),[k]}(\boldsymbol h;\Pi)\bigr).
\end{align*}
The future \(\Pi\) is independent of the input and shared by all
components. Joint stationarity means \(\nu\overline{\mathsf P}_t^{(c),[k]}=\nu\).
A common geometric update gives joint invariance, and parameter deletion
identifies the marginals. Its finite line representation then allows the
initial profiles and this shared update environment to converge together
to PNG.

\subsubsection{The Poisson-line construction}
\label{subsubsec:half-space-png-finite-parameter-line-law}

We give the precise Poisson-line construction of the jointly stationary
measures $\nu_{c;g_1,\ldots,g_k}$ from
Section~\ref{subsec:intro-stationary-measures}.

\paragraph*{\textbf{Bridge-inclusive c\`adl\`ag last passage}}
Let \(E=(E_p,\ldots,E_q)\) be c\`adl\`ag functions.  Write
\begin{equation*}
 \mathrm dE_j((r,s])=E_j(s)-E_j(r),\qquad
 \mathrm dE_j(\{r\})=E_j(r)-E_j(r-)
\end{equation*}
for increments and jumps;
these agree with signed Stieltjes measures when the lines have locally
finite variation.

\begin{definition}[Last passage with bridge atoms]
\label{def:half-space-png-cadlag-lpp}
Fix a finite set \(\mathcal M_{\mathrm{br}}\subset(-\infty,0)\) of marked
coordinates.  For \(u=(a,p)\), \(v=(d,q)\), \(a\leq d\), define
paths, weights, and passage values by
\begin{gather*}
a=x_{p-1}\leq x_p\leq\cdots\leq x_q=d.
\\
\operatorname{wt}_{E;\mathcal M_{\mathrm{br}}}
 (x_{p-1},\ldots,x_q)
 :=
 \sum_{j=p}^{q}
 \left\{
  \mathrm dE_j((x_{j-1},x_j])
  +\mathbf 1_{\{x_{j-1}\in\mathcal M_{\mathrm{br}}\}}
   \mathrm dE_j(\{x_{j-1}\})
 \right\}.
\\
\mathcal L_{E;\mathcal M_{\mathrm{br}}}(u;v)
 :=
 \sup_{a=x_{p-1}\leq\cdots\leq x_q=d}
 \operatorname{wt}_{E;\mathcal M_{\mathrm{br}}}
 (x_{p-1},\ldots,x_q).
\end{gather*}
The value is \(-\infty\) if \(a>d\) or \(p>q\).  When the marked set
is fixed by the environment, we abbreviate \(\mathcal L_{E;\mathcal M_{\mathrm{br}}}\) to \(\mathcal L_E\).
\end{definition}

An ordinary switch atom belongs to the departing line.  At marked
coordinates, the singleton also collects the bridge on each entered line,
including the source line, exactly once.

Every path from \((a,p)\) to \((d,q)\) has absolute weight at most
\begin{equation*}
 2\sum_{j=p}^q\max\{\sup_{[a,d]}|E_j|,|E_j(a-)|\}<\infty,
\end{equation*}
so the passage value is finite for \(a\leq d\), \(p\leq q\).
For finite step environments, switches at atoms or in complementary
gaps give finitely many weight patterns, so the supremum is attained.

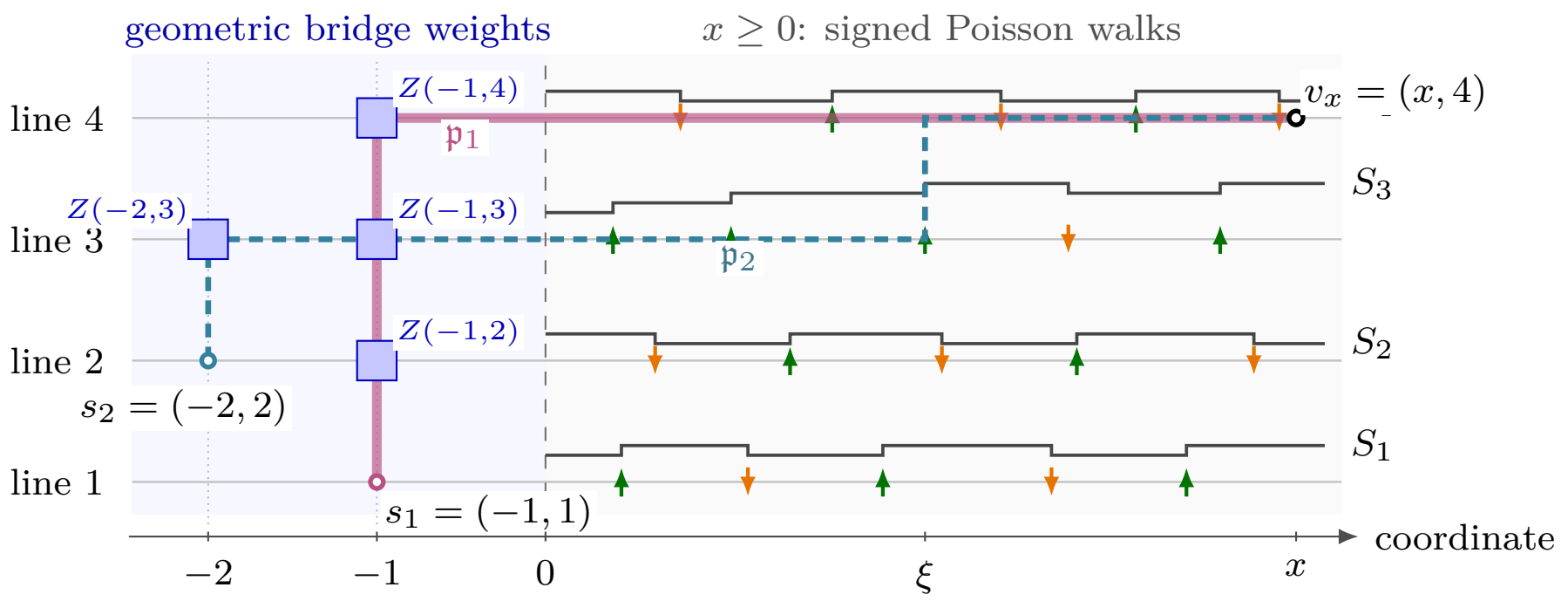
\begin{figure}[t]
\centering
\begin{tikzpicture}[
 x=1.28cm,
 y=0.92cm,
 every node/.style={font=\footnotesize},
 >={Latex[length=1.5mm,width=1.2mm]},
 carrier/.style={gray!48,line width=0.45pt},
 walk/.style={black!72,line width=0.7pt},
 bridge/.style={rectangle,draw=blue!75!black,fill=blue!22,
                minimum width=3mm,minimum height=3mm,inner sep=0pt},
 plusatom/.style={green!45!black,line width=0.8pt,
                 -{Latex[length=1.5mm,width=1.2mm]}},
 minusatom/.style={orange!90!black,line width=0.8pt,
                  -{Latex[length=1.5mm,width=1.2mm]}},
 pathone/.style={magenta!72!black,line width=2.1pt,opacity=0.66},
 pathtwo/.style={cyan!55!black,line width=1.25pt,densely dashed}
]
\fill[blue!3] (-2.45,0.73) rectangle (0,4.52);
\fill[black!2] (0,0.73) rectangle (4.72,4.52);
\node[blue!65!black] at (-1.23,4.72) {geometric bridge weights};
\node[black!70] at (2.35,4.72)
 {\(x\geq0\): signed Poisson walks};

\foreach \j in {1,2,3,4}{
 \draw[carrier] (-2.45,\j)--(4.72,\j);
 \node[anchor=east] at (-2.52,\j) {line \(\j\)};
}
\draw[gray!52,densely dotted] (-2,0.65)--(-2,4.48);
\draw[gray!52,densely dotted] (-1,0.65)--(-1,4.48);
\draw[black!55,dashed] (0,0.65)--(0,4.48);

\draw[walk]
 (0,1.22)--(0.45,1.22)--(0.45,1.30)--(1.20,1.30)
 --(1.20,1.22)--(2.00,1.22)--(2.00,1.30)--(3.00,1.30)
 --(3.00,1.22)--(3.80,1.22)--(3.80,1.30)--(4.62,1.30);
\node[anchor=west] at (4.68,1.30) {\(S_1\)};
\draw[walk]
 (0,2.22)--(0.65,2.22)--(0.65,2.14)--(1.45,2.14)
 --(1.45,2.22)--(2.35,2.22)--(2.35,2.14)--(3.15,2.14)
 --(3.15,2.22)--(4.20,2.22)--(4.20,2.14)--(4.62,2.14);
\node[anchor=west] at (4.68,2.14) {\(S_2\)};
\draw[walk]
 (0,3.22)--(0.40,3.22)--(0.40,3.30)--(1.10,3.30)
 --(1.10,3.38)--(2.25,3.38)--(2.25,3.46)--(3.10,3.46)
 --(3.10,3.38)--(4.00,3.38)--(4.00,3.46)--(4.62,3.46);
\node[anchor=west] at (4.68,3.46) {\(S_3\)};
\draw[walk]
 (0,4.22)--(0.80,4.22)--(0.80,4.14)--(1.70,4.14)
 --(1.70,4.22)--(2.70,4.22)--(2.70,4.14)--(3.50,4.14)
 --(3.50,4.22)--(4.35,4.22)--(4.35,4.14)--(4.62,4.14);
\node[anchor=west] at (4.68,4.14) {\(S_4\)};

\foreach \xx in {0.45,2.00,3.80}
 \draw[plusatom] (\xx,0.88)--(\xx,1.12);
\foreach \xx in {1.20,3.00}
 \draw[minusatom] (\xx,1.12)--(\xx,0.88);
\foreach \xx in {1.45,3.15}
 \draw[plusatom] (\xx,1.88)--(\xx,2.12);
\foreach \xx in {0.65,2.35,4.20}
 \draw[minusatom] (\xx,2.12)--(\xx,1.88);
\foreach \xx in {0.40,1.10,2.25,4.00}
 \draw[plusatom] (\xx,2.88)--(\xx,3.12);
\draw[minusatom] (3.10,3.12)--(3.10,2.88);
\foreach \xx in {1.70,3.50}
 \draw[plusatom] (\xx,3.88)--(\xx,4.12);
\foreach \xx in {0.80,2.70,4.35}
 \draw[minusatom] (\xx,4.12)--(\xx,3.88);

\draw[pathone] (-1,1)--(-1,4)--(4.45,4);
\draw[pathtwo] (-2,2)--(-2,3)--(2.25,3)--(2.25,4)--(4.45,4);
\node[magenta!72!black,fill=white,inner sep=1pt]
 at (-0.48,3.84) {\(\mathfrak p_1\)};
\node[cyan!45!black,fill=white,inner sep=1pt]
 at (1.15,2.84) {\(\mathfrak p_2\)};

\node[bridge] at (-1,2) {};
\node[anchor=south west,blue!75!black,fill=white,inner sep=0.8pt]
 at (-0.91,2.08) {\(\scriptstyle Z(-1,2)\)};
\node[bridge] at (-2,3) {};
\node[anchor=south east,blue!75!black,fill=white,inner sep=0.8pt]
 at (-2.09,3.08) {\(\scriptstyle Z(-2,3)\)};
\node[bridge] at (-1,3) {};
\node[anchor=south west,blue!75!black,fill=white,inner sep=0.8pt]
 at (-0.91,3.08) {\(\scriptstyle Z(-1,3)\)};
\node[bridge] at (-1,4) {};
\node[anchor=south west,blue!75!black,fill=white,inner sep=0.8pt]
 at (-0.91,4.08) {\(\scriptstyle Z(-1,4)\)};

\filldraw[fill=white,draw=cyan!55!black,line width=1pt]
 (-2,2) circle[radius=1.4pt];
\node[anchor=north east,fill=white,inner sep=0.8pt]
 at (-1.5,1.8) {\(s_2=(-2,2)\)};
\filldraw[fill=white,draw=magenta!72!black,line width=1pt]
 (-1,1) circle[radius=1.4pt];
\node[anchor=north west,fill=white,inner sep=0.8pt]
 at (-0.98,0.93) {\(s_1=(-1,1)\)};
\filldraw[fill=white,draw=black,line width=1pt]
 (4.45,4) circle[radius=1.5pt];
\node[anchor=south west,fill=white,inner sep=0.8pt]
 at (4.47,4.02) {\(v_x=(x,4)\)};

\draw[->,black!70] (-2.47,0.55)--(4.82,0.55);
\foreach \xx/\lab in {-2/{-2},-1/{-1},0/{0},2.25/{\xi},4.45/{x}}{
 \draw[black!70] (\xx,0.51)--(\xx,0.59);
 \node[anchor=north] at (\xx,0.49) {\(\lab\)};
}
\node[anchor=west] at (4.82,0.55) {coordinate};
\end{tikzpicture}
\caption{The shared \(k=2\) environment.  Squares are bridges;
the highlighted paths run from the two sources to \((x,4)\). The signed Poisson walks are schematic.}
\label{fig:half-space-png-k-two-candidate-law}
\end{figure}

\paragraph*{\textbf{Candidate jointly stationary measures}}
We specialize Definition~\ref{def:half-space-png-cadlag-lpp} to a
shared environment of signed Poisson lines and geometric bridge atoms.
Figure~\ref{fig:half-space-png-k-two-candidate-law} illustrates the
construction for $k=2$.

\begin{definition}[Poisson-line candidate law]
\label{def:half-space-png-joint-candidate-law}
Fix $c\geq0$ and $\boldsymbol g=(g_1,\ldots,g_k)$ with
$g_1>\cdots>g_k\geq r_c$ and $g_k>c$.
Write $\bar j:=2k+1-j$ and
\begin{equation*}
 (\gamma_1,\ldots,\gamma_{2k})
 :=(g_1,\ldots,g_k,g_k^{-1},\ldots,g_1^{-1}).
\end{equation*}
The environment and its passage profiles are defined as follows.
\begin{enumerate}[label=\textup{(\roman*)},leftmargin=*]

\item \textup{Poisson lines.}
Let $S_1^{\boldsymbol g},\ldots,S_{2k}^{\boldsymbol g}$ be independent
signed Poisson walks, started at zero, with line-$j$ up/down rates
$(\gamma_j,\gamma_j^{-1})$.

\item \textup{Geometric bridges.}
Independently of the walks, take independent variables $Z(-i,j)$,
$1\leq i\leq k$ and $i<j\leq\bar i$, with
\begin{equation}
 Z(-i,j)\sim
 \begin{cases}
  \operatorname{Geo}(\gamma_j/g_i),&i<j<\bar i,\\
  \operatorname{Geo}(c/g_i),&j=\bar i.
 \end{cases}
 \label{eq:half-space-png-bridge-laws-spectral}
\end{equation}
For each line $j$, let $A_j^{(c;\boldsymbol g)}$ be the
right-continuous step function on $(-\infty,0]$, normalized by
$A_j^{(c;\boldsymbol g)}(0)=0$, whose Stieltjes measure is
\begin{equation*}
 \mathrm dA_j^{(c;\boldsymbol g)}
 =\sum_{\substack{1\leq i\leq k\\i<j\leq\bar i}}
 Z(-i,j)\,\delta_{-i},
\end{equation*}
where $\delta_{-i}$ is the unit point mass at $-i$.
Concatenate the bridges and Poisson lines at zero:
\begin{equation*}
 E_j^{(c;\boldsymbol g)}(x):=
 \begin{cases}
  A_j^{(c;\boldsymbol g)}(x),&x\leq0,\\
  S_j^{\boldsymbol g}(x),&x\geq0.
 \end{cases}
\end{equation*}

\item \textup{Passage profiles and their joint law.}
Use the passage functional of
Definition~\ref{def:half-space-png-cadlag-lpp} with marked bridge
columns $\mathcal M_{\mathrm{br}}=\{-k,\ldots,-1\}$ and sources
$s_i=(-i,i)$. In this shared environment, define
\begin{align}
 G_i^{(c;\boldsymbol g)}(x)
 &:=\mathcal L_{E^{(c;\boldsymbol g)}}(s_i;(x,2k)),
 \label{eq:half-space-png-unpinned-joint-profile}\\
 \widetilde R_i^{(c;\boldsymbol g)}(x)
 &:=G_i^{(c;\boldsymbol g)}(x)-G_i^{(c;\boldsymbol g)}(0),
 \qquad
 R_i^{(c;\boldsymbol g)}
 :=\operatorname{usc}\widetilde R_i^{(c;\boldsymbol g)},
 \label{eq:half-space-png-joint-profile}
\end{align}
for $1\leq i\leq k$ and $x\geq0$.
The \emph{candidate jointly stationary law} is
\begin{equation*}
 \nu_{c;g_1,\ldots,g_k}
 :=\operatorname{Law}\bigl(
 R_1^{(c;\boldsymbol g)},\ldots,R_k^{(c;\boldsymbol g)}\bigr).
\end{equation*}
All coordinates use the same Poisson lines and bridge variables.
\end{enumerate}
\end{definition}

All geometric parameters in the definition lie in $[0,1)$.
In particular, when $c<1$ and $g_k=1$, the central terminal bridge
has law $\operatorname{Geo}(c)$ and remains well defined.
Each $G_i$ is a measurable, locally finite c\`adl\`ag step function.
Almost surely, the Poisson lines have no atom in some interval
immediately to the right of zero, so upper-semicontinuous regularization
preserves pinning. Hence $R_i^{(c;\boldsymbol g)}\in\mathrm{UC}_{+,0}$.

\smallskip
For illustration, take $k=2$ with $g_2>c$;
Figure~\ref{fig:half-space-png-k-two-candidate-law} shows the shared environment.

The bridge functions and weight distributions in this example are
\begin{equation*}
 \begin{aligned}
 \mathrm dA_1&=0,\\
 \mathrm dA_2&=Z(-1,2)\,\delta_{-1},\\
 \mathrm dA_3&=Z(-2,3)\,\delta_{-2}
                +Z(-1,3)\,\delta_{-1},\\
 \mathrm dA_4&=Z(-1,4)\,\delta_{-1},
\end{aligned}
\qquad
\begin{aligned}
 Z(-1,2)&\sim\operatorname{Geo}(g_2/g_1),\\
 Z(-2,3)&\sim\operatorname{Geo}(c/g_2),\\
 Z(-1,3)&\sim\operatorname{Geo}((g_1g_2)^{-1}),\\
 Z(-1,4)&\sim\operatorname{Geo}(c/g_1).
\end{aligned}
\end{equation*}
The two profiles use the same four lines and bridges, with entrances \((-1,1)\) and \((-2,2)\), respectively.

\paragraph*{\textbf{The one-parameter consistency check}}
\begin{lemma}[The one-parameter marginal]
\label{lem:half-space-png-one-parameter-check}
For \(k=1\) and \(g>c\), the profile \(R_1^{(c;g)}\) equals
\(H^{(c,g)}\) almost surely under the shared-walk coupling, where \(H^{(c,g)}\) is defined in
\eqref{eq:intro-half-space-png-stationary-profile}.
\end{lemma}

\begin{proof}
The bridge is \(Y:=Z(-1,2)\sim\operatorname{Geo}(c/g)\). Switching at \(-1\), in \((-1,0]\), or at \(y\in(0,x]\) gives \(Y+S_2^{(g)}(x)\), \(S_2^{(g)}(x)\), or \(S_1^{(g)}(y)+S_2^{(g)}(x)-S_2^{(g)}(y)\), respectively.  Including \(y=0\) in the last supremum gives
\begin{equation*}
 \begin{aligned}
  G_1^{(c;g)}(x)
  &=S_2^{(g)}(x)+\max\!\left\{Y,
    \sup_{0\leq y\leq x}(S_1^{(g)}(y)-S_2^{(g)}(y))\right\},\\
  \widetilde R_1^{(c;g)}(x)
  &=S_2^{(g)}(x)+\left(
    \sup_{0\leq y\leq x}\{S_1^{(g)}(y)-S_2^{(g)}(y)\}-Y
    \right)_{+}.
 \end{aligned}
\end{equation*}
Here \(G_1(0)=Y\); the second identity is therefore exactly
\eqref{eq:intro-half-space-png-stationary-profile-raw}, and applying \(\operatorname{usc}\) proves the claim.
\end{proof}

This includes \(c=0\), when \(Y=0\), and \(c<1\), \(g=1\), when the
walks have rates \((1,1)\) and \(Y\sim\operatorname{Geo}(c)\).

\paragraph*{\textbf{The singular endpoint}}
To complete the construction at \(g_k=c\geq1\), observe that the terminal bridge \(Z(-k,k+1)\) would have
parameter one.  For \(\boldsymbol g^{\partial}=(g_1,\ldots,g_{k-1},c)\), set \(g_k=c\) in the remaining laws, omit this bridge, and replace
the last source by
\begin{equation}
 s_k^{\partial}:=(-k+1,k+1).
\label{eq:half-space-png-singular-shifted-source}
\end{equation}
All retained geometric parameters are below one.  Put \(s_i^\partial=(-i,i)\) for \(i<k\), and define the profiles by
\eqref{eq:half-space-png-unpinned-joint-profile}--%
\eqref{eq:half-space-png-joint-profile} using the resulting environment \(E^\partial\) and sources \(s_i^\partial\).

\begin{figure}[t]
\centering
\begin{tikzpicture}[
 x=1.30cm,y=0.69cm,
 every node/.style={font=\footnotesize},
 >={Latex[length=1.5mm,width=1.2mm]},
 carrier/.style={gray!48,line width=0.45pt},
 bridge/.style={rectangle,draw=blue!75!black,fill=blue!22,
                minimum width=3mm,minimum height=3mm,inner sep=0pt},
 terminal/.style={rectangle,draw=orange!85!black,fill=orange!25,
                  minimum width=3.6mm,minimum height=3.6mm,inner sep=0pt},
 continuation/.style={cyan!55!black,line width=1.35pt},
 forced/.style={orange!85!black,line width=1.6pt}
]
\foreach \panel/\offset in {a/0,b/8.1}{
 \begin{scope}[xshift=\offset cm]
  \fill[blue!3] (-2.38,0.75) rectangle (0,4.48);
  \fill[black!2] (0,0.75) rectangle (1.75,4.48);
  \foreach \j in {1,2,3,4}{
   \draw[carrier] (-2.38,\j)--(1.75,\j);
   \node[anchor=east] at (-2.45,\j) {line \(\j\)};
  }
  \draw[gray!52,densely dotted] (-2,0.65)--(-2,4.42);
  \draw[gray!52,densely dotted] (-1,0.65)--(-1,4.42);
  \draw[black!55,dashed] (0,0.65)--(0,4.42);
  \draw[continuation] (-1,3)--(0.55,3)--(0.55,4)--(1.45,4);
  \node[bridge] at (-1,2) {};
  \node[bridge] at (-1,3) {};
  \node[bridge] at (-1,4) {};
  \filldraw[fill=white,draw=magenta!72!black,line width=1pt]
   (-1,1) circle[radius=1.5pt];
  \node[anchor=north west,fill=white,inner sep=0.8pt]
   at (-0.94,0.98) {\(s_1\)};
  \filldraw[fill=white,draw=black,line width=1pt]
   (1.45,4) circle[radius=1.5pt];
  \node[anchor=south,fill=white,inner sep=1pt]
   at (1.40,4.12) {\((x,4)\)};
  \draw[->,black!70] (-2.40,0.5)--(1.83,0.5);
  \foreach \xx/\lab in {-2/{-2},-1/{-1},0/{0},1.45/{x}}{
   \draw[black!70] (\xx,0.46)--(\xx,0.54);
   \node[anchor=north] at (\xx,0.43) {\(\lab\)};
  }
 \end{scope}
}
\node at (-0.30,4.88) {(a) \(g_2^{(n)}\downarrow c\)};
\draw[forced] (-2,2)--(-2,3);
\draw[continuation] (-2,3)--(-1,3);
\node[terminal] at (-2,3) {};
\node[anchor=south west,orange!85!black,fill=white,inner sep=1pt]
 at (-1.93,3.16) {\(Y_n\to\infty\)};
\filldraw[fill=white,draw=cyan!55!black,line width=1pt]
 (-2,2) circle[radius=1.5pt];
\node[anchor=north west,fill=white,inner sep=1pt]
 at (-1.92,1.92) {\(s_2=(-2,2)\)};
\begin{scope}[xshift=8.1cm]
 \node at (-0.30,4.88) {(b) \(g_2=c\geq1\)};
 \node[anchor=south,blue!75!black,fill=white,inner sep=1pt]
  at (-0.98,3.22) {\(Z(-1,3)\)};
 \draw[cyan!55!black,line width=1pt]
  (-1,3) circle[radius=3pt];
 \node[anchor=north west,fill=white,inner sep=1pt]
  at (-0.88,2.80) {\(s_2^\partial=(-1,3)\)};
\end{scope}
\end{tikzpicture}
\caption{The singular endpoint for \(k=2\).
Removing the diverging bridge \(Y_n=Z_n(-2,3)\) shifts
\(s_2=(-2,2)\) to \(s_2^\partial=(-1,3)\).
The bridge \(Z(-1,3)\) remains included at the new source.
Highlighted paths illustrate the shift; the first source is unchanged.}
\label{fig:half-space-png-singular-endpoint}
\end{figure}
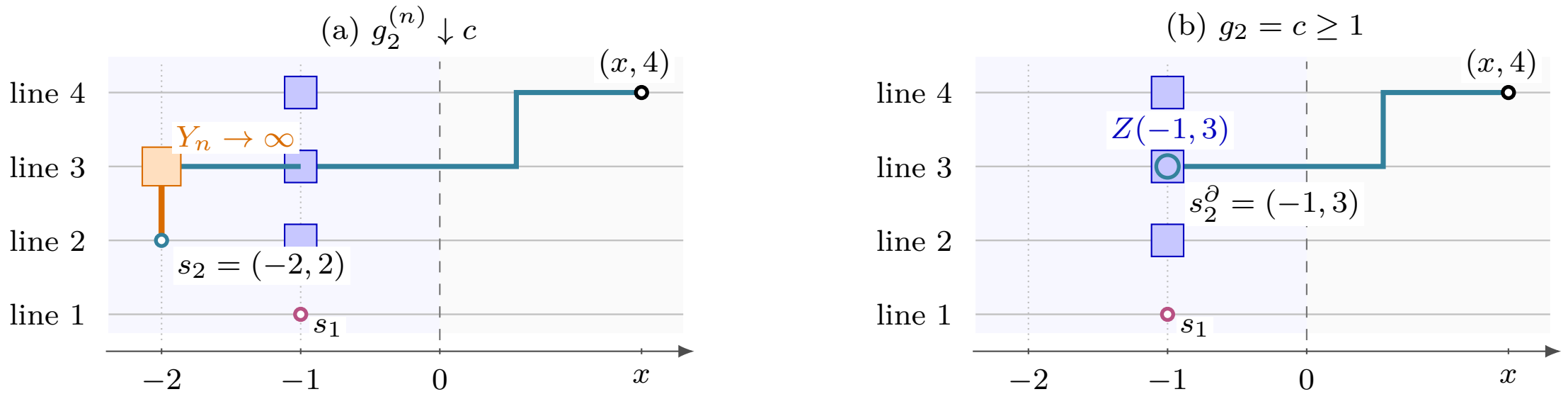

Figure~\ref{fig:half-space-png-singular-endpoint} illustrates the
source shift for \(k=2\).  The following lemma justifies the construction.

\begin{lemma}[A diverging terminal bridge]
\label{lem:half-space-png-forced-terminal-bridge}
Let \(c\geq1\), fix the first \(k-1\) spectral parameters, and let \(g_k^{(n)}\downarrow c\) through strictly ordered admissible vectors \(\boldsymbol g^{(n)}=(g_1,\ldots,g_{k-1},g_k^{(n)})\). Writing \(\boldsymbol g^\partial=(g_1,\ldots,g_{k-1},c)\), one can
couple the profiles so that, for every \(L<\infty\),
\begin{equation*}
 \mathbb P\left(R_i^{(c;\boldsymbol g^{(n)})}
 =R_i^{(c;\boldsymbol g^\partial)}\text{ on }[0,L]
 \text{ for every }i\leq k\right)\longrightarrow1.
\end{equation*}
\end{lemma}

\begin{proof}
Remove \(Y_n=Z_n(-k,k+1)\sim\operatorname{Geo}(c/g_k^{(n)})\) to obtain \(E^{(n),\circ}\).  Its strip \([-k,-k+1)\) is empty, so
\begin{equation*}
 \mathcal L_{E^{(n),\circ}}((-k,k+1);(x,2k))
=\mathcal L_{E^{(n),\circ}}(s_k^\partial;(x,2k)).
\end{equation*}
Indeed, moving switches to \(-k+1\) can only gain nonnegative bridges;
conversely extend a shifted path backwards on line \(k+1\).  The
interval atom at \(-k+1\) equals the shifted source singleton.

The total variation of \(E^{(n),\circ}\) on lines \(k,\ldots,2k\) and \([-k,L+1]\) gives a tight bound \(\mathfrak w_n\) for all nonterminal
weights.  Since \(Y_n\to\infty\) in probability, \(\mathbb P(Y_n>2\mathfrak w_n)\to1\).  On this event every maximizer for the
last source and an exit in \([0,L+1]\) uses \(Y_n\), and
\begin{equation*}
 G_k^{(c;\boldsymbol g^{(n)})}(x)
=Y_n+\mathcal L_{E^{(n),\circ}}(s_k^\partial;(x,2k)).
\end{equation*}
Pinning cancels \(Y_n\).  Couple the remaining bridges maximally and
the Poisson clocks by common rate cutoffs; their configurations agree
with the endpoint environment on this compact window with probability
tending to one.  The other sources cannot access \(Y_n\), and working on \([0,L+1]\) preserves equality on \([0,L]\) after upper-semicontinuous regularization.
\end{proof}

For \(k=1\), the shifted source is \(s_1^{\partial}=(0,2)\). There is then no line choice and no bridge at the source, so
\begin{equation}
 R_1^{(c;c)}=\operatorname{usc}S_2^{(c)},
\label{eq:half-space-png-singular-k-one}
\end{equation}
in agreement with
\eqref{eq:intro-half-space-png-boundary-dominated-profile}.

For arbitrary parameter lists, use the repetition and permutation
convention of Section~\ref{subsec:intro-stationary-measures}.

\subsubsection{Geometric stationarity and consistency}
\label{subsubsec:half-space-png-equal-time-geometric-approximation}

We construct a geometric approximation with three properties:
consistent spatial-window laws, joint invariance under a shared
discrete PNG update, and consistency under deletion of spectral
parameters. The next subsections identify its Poisson-line limit and
transfer these properties to PNG. Assume
\begin{equation}
 g_1>\cdots>g_k>c\vee1,\qquad c>0.
\label{eq:half-space-png-strict-interior-parameters}
\end{equation}
The physical mesh $\delta$ remains fixed throughout this subsection.
The auxiliary parameter $\eta$ regularizes reciprocal pairs; we take
$\eta\downarrow0$ after cancelling source weights. The limit
$\delta\downarrow0$ is taken later.  Write \(\bar j=2k+1-j\) and set
\begin{equation*}
 \gamma_i^{(\eta)}=g_i,\qquad
 \gamma_{\bar i}^{(\eta)}=(1-\eta)g_i^{-1},\qquad 1\leq i\leq k.
\end{equation*}
Choose \(0<\delta<1\), \(c\delta<1\), \(\delta g_1<1\), and \(0<\eta<1-\delta g_1\).  For a cut base \(a\in\mathbb Z\), define
\begin{equation*}
 \alpha_p^{\delta,\eta,a}=
 \begin{cases}
  \gamma_{p-a}^{(\eta)},&a+1\leq p\leq a+2k,\\
  \delta,&\text{otherwise}.
 \end{cases}
\end{equation*}
Let \(X^{\delta,\eta,a}\) be the extended half-space geometric
environment of \cite[Definition~2.2]{dauvergne2026directed}, with
boundary parameter \(c\) and coordinate parameters \(\alpha^{\delta,\eta,a}\).  Thus the variables are independent over
reflection orbits, \(X(p,q)=X(q,p)\), and
\begin{equation*}
 X(p,q)\sim
 \begin{cases}
  \operatorname{Geo}(\alpha_p^{\delta,\eta,a}
                            \alpha_q^{\delta,\eta,a}),&p\neq q,\\
  \operatorname{Geo}(c\alpha_p^{\delta,\eta,a}),&p=q,
 \end{cases}
\end{equation*}
for parameters in \((0,1)\), and \(+\infty\) otherwise.
For \(u,v\in\mathbb Z^2\), let \(\mathcal P(u,v)\) consist of the vertex
paths from \(u\) to \(v\) with steps \((1,0)\) and \((0,1)\), subject to
\((r,r)\in\pi\Rightarrow(r-1,r)\notin\pi\) for every \(r\in\mathbb Z\).
Define the extended half-space last-passage value by
\begin{equation*}
 L_X(u;v):=\max_{\pi\in\mathcal P(u,v)}
             \sum_{w\in\pi}X(w),
 \qquad \max\varnothing:=-\infty .
\end{equation*}
Both endpoint weights are included.

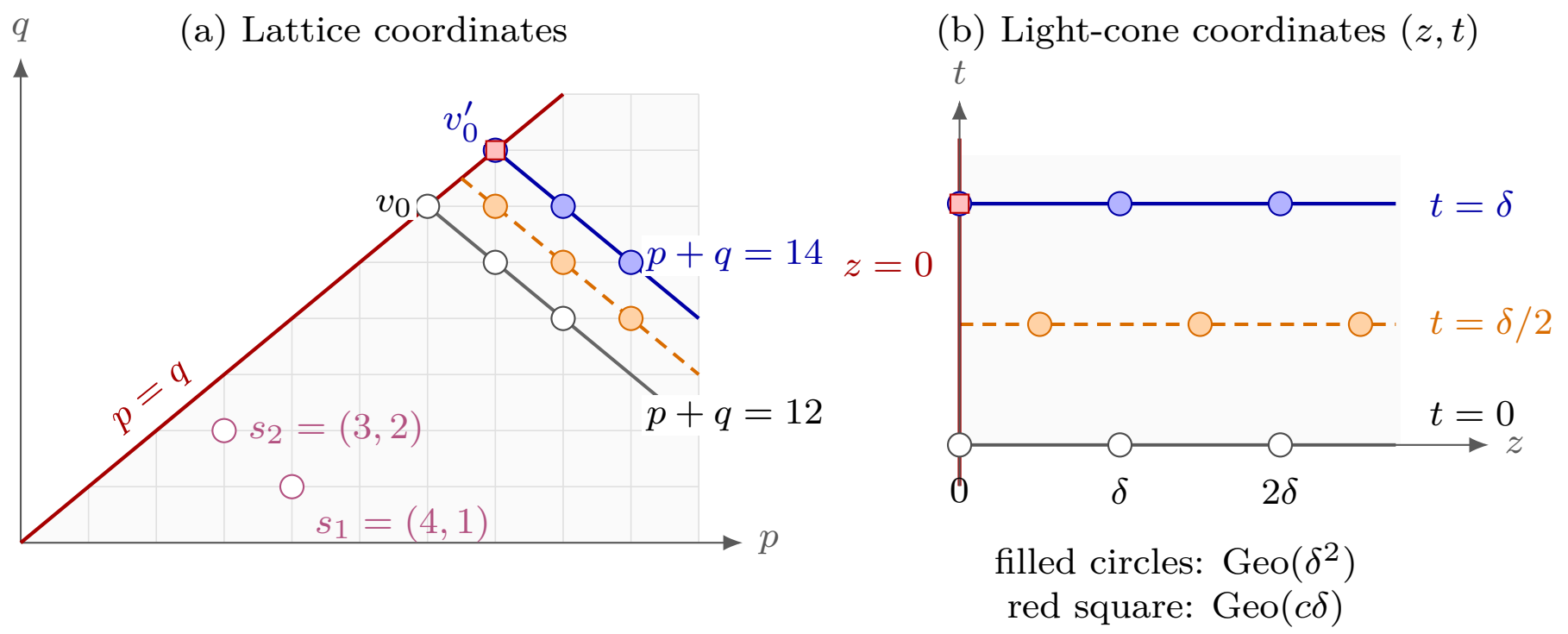
\begin{figure}[t]
\centering
\begin{tikzpicture}[
 every node/.style={font=\footnotesize},
 >={Latex[length=1.5mm,width=1.2mm]},
 cutinput/.style={black!60,line width=0.75pt},
 cutodd/.style={orange!85!black,densely dashed,line width=0.75pt},
 cuteven/.style={blue!65!black,line width=0.8pt},
 oddweight/.style={circle,draw=orange!85!black,fill=orange!35,
                   inner sep=1.8pt},
 evenweight/.style={circle,draw=blue!65!black,fill=blue!30,
                    inner sep=1.8pt},
 boundaryweight/.style={rectangle,draw=red!75!black,fill=red!25,
                        inner sep=2pt},
 inputpoint/.style={circle,draw=black!70,fill=white,inner sep=1.8pt},
 sourcepoint/.style={circle,draw=magenta!70!black,fill=white,inner sep=1.8pt}
]
\begin{scope}[x=0.52cm,y=0.43cm]
 \fill[black!2] (0,0)--(10,0)--(10,8)--(8,8)--cycle;
 \foreach \p in {0,...,10}{
  \pgfmathsetmacro{\qtop}{min(\p,8)}
  \draw[black!12,line width=0.3pt] (\p,0)--(\p,\qtop);
 }
 \foreach \q in {0,...,8}
  \draw[black!12,line width=0.3pt] (\q,\q)--(10,\q);
 \draw[->,black!65] (0,0)--(10.65,0) node[right] {$p$};
 \draw[->,black!65] (0,0)--(0,8.65) node[above] {$q$};
 \draw[red!65!black,line width=0.8pt] (0,0)--(8,8);
 \node[red!65!black,rotate=39,above] at (2.2,2.2) {$p=q$};
 \node at (5.2,9.1) {(a) Lattice coordinates};
 \draw[cutinput] (6,6)--(10,2);
 \draw[cutodd] (6.5,6.5)--(10,3);
 \draw[cuteven] (7,7)--(10,4);
 \node[anchor=north west,fill=white,inner sep=1pt] at (9.15,2.65)
  {$p+q=12$};
 \node[anchor=south west,blue!65!black,fill=white,inner sep=1pt]
  at (9.15,4.75) {$p+q=14$};
 \foreach \m in {0,1,2}{
  \node[inputpoint] at ({6+\m},{6-\m}) {};
  \node[evenweight] at ({7+\m},{7-\m}) {};
  \node[oddweight] at ({7+\m},{6-\m}) {};
 }
 \node[boundaryweight] at (7,7) {};
 \node[anchor=east,fill=white,inner sep=0.7pt] at (5.85,6.0) {$v_0$};
 \node[anchor=south east,blue!65!black,fill=white,inner sep=0.7pt]
  at (6.85,7.08) {$v'_0$};
 \node[sourcepoint] at (4,1) {};
 \node[anchor=north west,magenta!70!black] at (4.1,0.95) {$s_1=(4,1)$};
 \node[sourcepoint] at (3,2) {};
 \node[anchor=west,magenta!70!black] at (3.12,2.03) {$s_2=(3,2)$};
\end{scope}
\begin{scope}[xshift=7.2cm,yshift=0.75cm,x=1.23cm,y=1.85cm]
 \node at (1.55,1.71) {(b) Light-cone coordinates $(z,t)$};
 \fill[black!2] (0,0) rectangle (2.75,1.2);
 \draw[cutinput] (0,0)--(2.72,0);
 \draw[cutodd] (0,0.5)--(2.72,0.5);
 \draw[cuteven] (0,1)--(2.72,1);
 \draw[red!65!black,line width=0.8pt] (0,-0.17)--(0,1.27);
 \draw[->,black!65] (0,-0.17)--(0,1.43) node[above] {$t$};
 \draw[->,black!65] (0,0)--(3.30,0) node[right] {$z$};
 \node[anchor=east,red!65!black] at (-0.06,0.75) {$z=0$};
 \foreach \m in {0,1,2}{
  \node[inputpoint] at (\m,0) {};
  \node[evenweight] at (\m,1) {};
  \node[oddweight] at ({\m+0.5},0.5) {};
 }
 \node[boundaryweight] at (0,1) {};
 \foreach \m/\lab in {0/0,1/{\delta},2/{2\delta}}
  \node[anchor=north] at (\m,-0.07) {$\lab$};
 \node[anchor=west] at (2.83,0.13) {$t=0$};
 \node[anchor=west,orange!85!black] at (2.83,0.5) {$t=\delta/2$};
 \node[anchor=west,blue!65!black] at (2.83,1) {$t=\delta$};
 \node[anchor=north,align=center] at (1.35,-0.34)
  {filled circles: $\operatorname{Geo}(\delta^2)$\\
   red square: $\operatorname{Geo}(c\delta)$};
\end{scope}
\end{tikzpicture}
\caption{Antidiagonal cuts and the two-layer update, with
$k=2$, $a=0$, $M=2$, and $N=6$.
The vertices $v_m=(6+m,6-m)$ and $v'_m=v_m+(1,1)$,
$m=0,1,2$, map to $(t,z)=(0,m\delta)$ and $(\delta,m\delta)$
under $\Phi_{\delta,N}$; panel (b) plots space horizontally and time
vertically. Orange and blue vertices belong to the two fresh layers;
the red square is the new diagonal vertex.}
\label{fig:half-space-png-cut-embedding}
\end{figure}

For \(1\leq i\leq k\) and \(M\in\mathbb N_0\), put
\begin{equation*}
 s_i^a=(a+\bar i,a+i),\qquad N_M=a+2k+M,\qquad
 v_m^M=(N_M+m,N_M-m),\quad 0\leq m\leq M.
\end{equation*}
Define the pinned antidiagonal field
\begin{equation}
 H_{i,M}^{\delta,\eta,a}(m)
 =L_{X^{\delta,\eta,a}}(s_i^a;v_m^M)
  -L_{X^{\delta,\eta,a}}(s_i^a;v_0^M).
\label{eq:half-space-png-geometric-pinned-antidiagonal-field}
\end{equation}
For the sources and terminals in
\eqref{eq:half-space-png-geometric-pinned-antidiagonal-field}, the path
restriction keeps every admissible path in \(\{(p,q):p\geq q\}\).
Simultaneous translation of the parameters, sources, and terminals
by $(-a,-a)$ shows that these window laws do not depend on $a$.
Define the light-cone map by
\begin{equation}
 \Phi_{\delta,N}(p,q)
 =\left(\delta\left(\frac{p+q}{2}-N\right),
         \frac{\delta}{2}(p-q)\right).
\label{eq:half-space-png-null-coordinate-embedding}
\end{equation}
With $N=N_M$, this map sends \(v_m^M\) to \((0,\delta m)\) and \(v_m^M+(n,n)\) to \((n\delta,\delta m)\),
so the antidiagonal field is read at one physical time.
Figure~\ref{fig:half-space-png-cut-embedding} shows the cuts and the two
fresh lattice layers in both coordinate systems.

\begin{lemma}[Accessible cells and cancellation of source weights]
\label{lem:half-space-png-accessible-cells-and-forced-limit}
For a path starting at $s_i^a$, the source is the only accessible
cell whose geometric parameter tends to one as $\eta\downarrow0$.
Its weight occurs in every path and cancels in
\eqref{eq:half-space-png-geometric-pinned-antidiagonal-field}.
The remaining weights can be coupled so that, for every fixed $M$,
jointly in $i\leq k$ and $0\leq m\leq M$,
\begin{equation*}
 H_{i,M}^{\delta,\eta,a}(m)
 \xrightarrow[\eta\downarrow0]{\mathrm{a.s.}}
 H_{i,M}^{\delta,a}(m).
\end{equation*}
\end{lemma}

\begin{proof}
The source parameter is
$\gamma_{\bar i}^{(\eta)}\gamma_i^{(\eta)}=1-\eta$.
An accessible mixed special cell with labels $\bar r,\ell$ satisfies
$r\leq i\leq\ell$, so its parameter is
$ (1-\eta)\frac{g_\ell}{g_r}\leq1-\eta$,
with equality only for $r=\ell=i$. Every other mixed parameter has
a limit strictly below one. Two reflected labels give a parameter
at most $(1-\eta)^2/g_k^2$; a background coordinate gives at most
$\delta g_1$ or $\delta^2$; an accessible special diagonal parameter
is at most $c(1-\eta)/g_i$. All these limits are below one.
Cells with two unreflected special labels are inaccessible.
Thus every accessible weight is finite at positive $\eta$, and all
but the source weight have proper limits. Common-uniform coupling
of the finitely many remaining weights and finite maximization
prove the assertion after the source contribution is cancelled.
\end{proof}

This cancellation parallels
\cite[Proposition~2.8]{dauvergne2026directed}.
All parameter swaps below are performed at positive $\eta$.

\paragraph*{\textbf{Parameter swaps and moving sources}}
Let $X^{\langle\ell\rangle}$ have the coordinate parameters at
$\ell,\ell+1$ interchanged. The shifted version of
\cite[Lemma~2.6]{dauvergne2026directed} gives
\begin{equation}
 \{L_X((p,q);(r,s)):p\leq r,\ q\leq s,
        \ p,q\neq\ell+1,\ r,s\neq\ell\}
 \overset{\mathrm d}{=}
 \{L_{X^{\langle\ell\rangle}}((p,q);(r,s)):p\leq r,\ q\leq s,
        \ p,q\neq\ell+1,\ r,s\neq\ell\}.
 \label{eq:half-space-png-shifted-full-field-parameter-swap}
\end{equation}
This identity applies jointly to all retained sources and terminals.
We first apply it to unpinned passage values, then take their pinned
differences.

The swaps also require transporting paired sources with their labels.
The following observation explains these source shifts. In a
post-swap environment, suppose that the new source $w$ is one
admissible north or east step from its old source $u$ and can reach
every retained terminal, that
$Y^\eta:=X(w)\sim\operatorname{Geo}(1-\eta)$, and that every other
weight accessible from $u$ to the retained terminals has a proper
limit. Their finite sum $\mathfrak w^\eta$ is tight, whereas
$Y^\eta\to\infty$ in probability. On
$\{Y^\eta>\mathfrak w^\eta\}$, every maximizer visits $w$:
a path through $w$ has weight at least $Y^\eta$, while a path
avoiding it has weight at most $\mathfrak w^\eta$.
Endpoint inclusion and adjacency then give
\begin{equation*}
 L_X(u;v)=X(u)+L_X(w;v)
\end{equation*}
for every retained terminal $v$. The extra term cancels under
pinning. Thus exact equality in law under a swap can be combined
with a source shift whose pinned fields agree with probability
tending to one as $\eta\downarrow0$.

\paragraph*{\textbf{Cut consistency and the common update}}
One discrete time step adds two independent lattice layers. Let
$(\xi_m^{\mathrm o})_{m\geq0}$, $(\xi_m^{\mathrm e})_{m\geq1}$,
and $d$ be mutually independent and independent of the input, with
$\xi_m^{\mathrm o},\xi_m^{\mathrm e}\sim\operatorname{Geo}(\delta^2)$
and $d\sim\operatorname{Geo}(c\delta)$. For pinned $h$, set
\begin{align}
 K_h\!\left(m+\tfrac12\right)
 &=\xi_m^{\mathrm o}+\max\{h(m),h(m+1)\}, &&m\geq0,
 \label{eq:half-space-png-discrete-png-half-step}\\
 \mathfrak r_h&:=d+K_h(\tfrac12),\qquad
 (\Psi_\delta^{(c)}h)(0):=0,\notag\\
 (\Psi_\delta^{(c)}h)(m)
 &:=\xi_m^{\mathrm e}
   +\max\{K_h(m-\tfrac12),K_h(m+\tfrac12)\}-\mathfrak r_h,
 &&m\geq1.
 \label{eq:half-space-png-discrete-png-update}
\end{align}
Use the same fresh array $\xi=(\xi^{\mathrm o},\xi^{\mathrm e},d)$
for every component, each with its own recentering constant.
Write $\Psi_{\delta,\xi}^{(c),[k]}$ for this shared update,
or $\Psi_\delta^{(c),[k]}$ when the noise is implicit.
Different time steps use independent fresh arrays.

\begin{lemma}[Cut consistency and stationarity]
\label{lem:half-space-png-equal-time-parameter-train}
For every $M\in\mathbb N_0$,
\begin{equation}
 \bigl(H_{i,M+1}^{\delta,a}(m):i\leq k,\ 0\leq m\leq M\bigr)
 \overset{\mathrm d}{=}
 \bigl(H_{i,M}^{\delta,a}(m):i\leq k,\ 0\leq m\leq M\bigr).
 \label{eq:half-space-png-antidiagonal-projective-identity}
\end{equation}
These consistent window laws define a law
$\mu_{c;\boldsymbol g}^{\delta}$ on $(\mathbb Z^{\mathbb N_0})^k$,
independent of $a$. This half-line law is invariant under the common
two-layer update $\Psi_\delta^{(c),[k]}$.
\end{lemma}

\begin{proof}
\emph{Step 1: consistency of the cut laws.}
At positive $\eta$, move the background parameter at $a+2k+1$
to $a+1$ through the swaps
$\tau_{a+2k},\tau_{a+2k-1},\ldots,\tau_{a+1}$,
where $\tau_\ell$ interchanges coordinates $\ell,\ell+1$.
At stage $j\in\{0,\ldots,2k\}$, the special labels occupy
\begin{equation*}
 p_h^{[j]}:=a+h+\mathbf1_{\{h>j\}},\qquad
 s_i^{[j]}:=(p_{\bar i}^{[j]},p_i^{[j]}),\qquad 1\leq h\leq2k.
\end{equation*}
The background parameter is at $a+j+1$. Let $X^{\eta,[j]}$
be this environment, and let $H_i^{\eta,[j]}(m)$ denote passage
from $s_i^{[j]}$ to $v_m^{M+1}$, pinned at $v_0^{M+1}$,
for $0\leq m\leq M$.

In the step $j\to j-1$, put $\ell=a+j$.
No source uses the background coordinate $\ell+1$, and
\begin{equation*}
 \min\{(v_m^{M+1})_1,(v_m^{M+1})_2\}
 =a+2k+M+1-m\geq a+2k+1>\ell.
\end{equation*}
Thus \eqref{eq:half-space-png-shifted-full-field-parameter-swap}
applies. If $\widetilde H_i^{\eta,[j]}$ is the post-swap field
from the old source $s_i^{[j]}$, it gives exact joint equality in law
between $H^{\eta,[j]}$ and $\widetilde H^{\eta,[j]}$.

Only $i_j=\min\{j,2k+1-j\}$ changes its source, since
\begin{equation*}
 s_i^{[j-1]}-s_i^{[j]}
 =(\mathbf1_{\{j=\bar i\}},\mathbf1_{\{j=i\}}).
\end{equation*}
Its new source has parameter $1-\eta$.
The accessibility argument of
Lemma~\ref{lem:half-space-png-accessible-cells-and-forced-limit}
extends to these intermediate orders, which insert one background
coordinate among the ordered special labels. From the affected old
source, an accessible mixed special cell still has labels
$\bar r,\ell'$ with $r\leq i_j\leq\ell'$.
Its limiting parameter $g_{\ell'}/g_r$ is strictly below one except
at the new paired source. Any additional cell involving the inserted
background coordinate has parameter bounded by $\delta g_1$ or
$\delta^2$; reflected and diagonal cells have the same proper bounds
as in that lemma. Hence all weights other than the new source have
parameters uniformly below one as $\eta\downarrow0$.

The source-shift observation above therefore yields
\begin{equation*}
 \mathbb P\bigl(\widetilde H^{\eta,[j]}
                    =H^{\eta,[j-1]}\bigr)\longrightarrow1,
\end{equation*}
where the arrays include all $i\leq k$ and $0\leq m\leq M$.
The same accessibility bounds give the pinned limits at each stage.
Combining these limits with the exact swap identities and telescoping
bounded-Lipschitz expectations over the $2k$ swaps identifies the
limiting laws at stages $2k$ and $0$.

Finally, $s_i^{[0]}=s_i^{[2k]}+(1,1)$.
Translation by $(-1,-1)$ restores the original parameter sequence
and sends $v_m^{M+1}$ to $v_m^M$, proving
\eqref{eq:half-space-png-antidiagonal-projective-identity}.
Kolmogorov consistency now defines the half-line law
$\mu_{c;\boldsymbol g}^{\delta}$; its independence of $a$ follows
from the diagonal translation noted above.

\emph{Step 2: invariance under the common update.}
Output coordinates $0,\ldots,M$ use only input coordinates
$0,\ldots,M+1$. Realize this input prefix on the cut
$N=N_{M+1}=a+2k+M+1$.
The fresh vertices are
\begin{equation*}
 o_m=(N+m+1,N-m),\qquad
 e_m=(N+m+1,N-m+1),\qquad 0\leq m\leq M.
\end{equation*}
Both coordinates exceed $a+2k$, so these are background sites:
every $o_m$ and each $e_m$ with $m\geq1$ have parameter $\delta^2$, while
$e_0=(N+1,N+1)$ has parameter $c\delta$.
The input passage values use vertices with $p+q\leq2N$, whereas
the fresh layers have $p+q=2N+1,2N+2$. Reflection preserves $p+q$,
so their reflection orbits are disjoint and the fresh weights are
independent of the entire input prefix.

The Bellman recurrence gives
\eqref{eq:half-space-png-discrete-png-half-step}--%
\eqref{eq:half-space-png-discrete-png-update} simultaneously for all
sources. The new height at the origin before recentering is
$d+K_h(1/2)$, whose subtraction pins each output.
The new cut is $N+1=N_{M+2}$, so the output
prefix has law
$\bigl(H_{i,M+2}^{\delta,a}(m):i\leq k,\ 0\leq m\leq M\bigr)$.
Step 1 identifies this with the original $M$-window law.

Each output coordinate depends on finitely many input coordinates
and fresh weights. Equality on all finite cylinders therefore proves
invariance of the full law; equivalently, for every bounded measurable $F$,
\begin{equation*}
 \int\mathbb E_\xi
 F\bigl(\Psi_{\delta,\xi}^{(c),[k]}(\boldsymbol h)\bigr)
 \,\mu_{c;\boldsymbol g}^{\delta}(\mathrm d\boldsymbol h)
 =\int F(\boldsymbol h)\,
       \mu_{c;\boldsymbol g}^{\delta}(\mathrm d\boldsymbol h).
\end{equation*}
\end{proof}

\paragraph*{\textbf{Deleting a spectral parameter}}
We next move the two labels representing one spectral parameter
outside the region accessible to the retained sources. Write
\begin{equation*}
 \boldsymbol g^{\widehat m}:=
 (g_1,\ldots,g_{m-1},g_{m+1},\ldots,g_k),\qquad
 \mathsf D_m(h_1,\ldots,h_k):=
 (h_1,\ldots,h_{m-1},h_{m+1},\ldots,h_k).
\end{equation*}

\begin{lemma}[Paired-parameter deletion]
\label{lem:half-space-png-geometric-paired-deletion}
Assume $k\geq2$ and
\eqref{eq:half-space-png-strict-interior-parameters}.
For every admissible $\delta>0$ and $1\leq m\leq k$,
\begin{equation}
 (\mathsf D_m)_*\mu_{c;\boldsymbol g}^{\delta}
 =\mu_{c;\boldsymbol g^{\widehat m}}^{\delta}
 \quad\text{on }(\mathbb Z^{\mathbb N_0})^{k-1}.
 \label{eq:half-space-png-geometric-sequence-deletion}
\end{equation}
\end{lemma}

\begin{proof}
By translation, take $a=0$. Put $N=\bar m$ and retain the
regularization $\gamma_i^{(\eta)}=g_i$,
$\gamma_{\bar i}^{(\eta)}=(1-\eta)/g_i$.
Move label $N$ to the first position and then label $m$ to the
second, through adjacent swaps:
\begin{equation}
 \begin{aligned}
 \sigma_{\mathrm{mid}}
 &:=(N,1,\ldots,N-1,N+1,\ldots,2k),\\
 (1,\ldots,2k)&\longrightarrow\sigma_{\mathrm{mid}}
 \longrightarrow
 (N,m,1,\ldots,m-1,m+1,\ldots,N-1,N+1,\ldots,2k).
 \end{aligned}
 \label{eq:half-space-png-deletion-parameter-train}
\end{equation}
For an intermediate order $\sigma$, assign parameter
$\gamma_{\sigma(p)}^{(\eta)}$ to position $1\leq p\leq2k$
and $\delta$ elsewhere, and denote the environment by
$X^{\eta,\sigma}$. The retained source for $i\neq m$ is
\begin{equation*}
 s_i^\sigma:=\bigl(\sigma^{-1}(\bar i),\sigma^{-1}(i)\bigr).
\end{equation*}
Use finite antidiagonal terminal sets with coordinates at least $2k$.
Before a swap at $d,d+1$, the moving deleted label is at $d+1$,
which no retained source uses. The terminals avoid $d\leq2k-1$.
Thus \eqref{eq:half-space-png-shifted-full-field-parameter-swap}
applies jointly to every retained source and terminal.

When $N$ crosses $m$, no retained source moves. In every other
swap, exactly one source moves to an adjacent paired source of
parameter $1-\eta$. The possibilities are as follows; the affected
component is $i$ in each row, and the last column gives the
post-swap parameter at its old source:
\begin{equation*}
 \begin{array}{c|c|c|c}
 \text{moving label}&\text{crossed label}&s^+-s^-&
       \text{old-source parameter}\\ \hline
 N&\bar i,\ i>m&(1,0)&(1-\eta)g_i/g_m\\
 N&i\leq k,\ i\neq m&(0,1)&(1-\eta)^2/(g_i g_m)\\
 m&i<m&(0,1)&(1-\eta)g_m/g_i
 \end{array}
\end{equation*}
All three limiting parameters are strictly below one.
The relative orders of the retained unreflected labels and of the
retained reflected labels stay unchanged. Hence an accessible mixed
cell with retained labels $\bar r,\ell$ satisfies $r\leq i\leq\ell$,
and its parameter $(1-\eta)g_\ell/g_r$ has a proper limit unless
$r=\ell=i$, the new source.
Products involving a moved deleted label are controlled by the
same inequalities $i>m$ or $i<m$ displayed in the table; the deleted
pair's own cell is inaccessible. Reflected, diagonal, and background
cells satisfy the proper bounds used above. Thus the sum of all
other accessible weights is tight on each finite terminal cone.

The source-shift observation now applies at every step: with
probability tending to one, all retained-terminal maximizers from
the affected old source visit the new source, and
\begin{equation*}
 L_{X_{\mathrm{post}}^{(\eta)}}(s^-;v)
 =X_{\mathrm{post}}^{(\eta)}(s^-)
  +L_{X_{\mathrm{post}}^{(\eta)}}(s^+;v).
\end{equation*}
Pinning removes the extra term. Combining these source shifts with
the exact joint swap identities and taking $\eta\downarrow0$
identifies the retained pinned fields before and after
\eqref{eq:half-space-png-deletion-parameter-train}.

In the final order, both coordinates of every retained source are
at least $3$, so the deleted positions $1,2$ are inaccessible.
Remove them and translate by $(-2,-2)$.
If $\operatorname{pos}(h)$ denotes the final position of label $h$,
then
\begin{equation*}
 (\operatorname{pos}(\bar i),\operatorname{pos}(i))
 =\begin{cases}
   (\bar i,i+2),&i<m,\\
   (\bar i+1,i+1),&i>m.
  \end{cases}
\end{equation*}
After translation, these are precisely the paired sources of the
reduced vector, indexed by $i$ for $i<m$ and $i-1$ for $i>m$.
The accessible parameter environment also becomes that of the
reduced vector. Finally, for
\begin{equation*}
 V_M^{(r)}:=\{(r+M+j,r+M-j):0\leq j\leq M\},
 \qquad V_M^{(2k)}-(2,2)=V_M^{(2k-2)},
\end{equation*}
the corresponding $M$-window laws agree. Since this holds for every
$M$, it proves \eqref{eq:half-space-png-geometric-sequence-deletion}
on the half-line sequence space.
\end{proof}

\subsubsection{A line representation of the geometric field}
\label{subsubsec:half-space-png-staggered-cut-readout}

We identify the geometric stationary law
$\mu_{c;\boldsymbol g}^{\delta}$ with a discrete version of the
Poisson-line construction. The representation uses independent
geometric bridges and signed geometric increments.
It is exact when the atoms in each mesh interval are placed in a
specified order. Moving them to the mesh points changes their mesh
readouts on a fixed window with probability at most $O(\delta)$.

Throughout, $c>0$, $g_1>\cdots>g_k>c\vee1$, and $\delta$ satisfies
the restrictions of the preceding subsection. Put $n=2k$,
$\bar j=n+1-j$, $\gamma_i=g_i$, and
$\gamma_{\bar i}=g_i^{-1}$. We first take $\eta\downarrow0$ at
fixed $\delta$ and set the cancelled source weights to zero, as in
Lemma~\ref{lem:half-space-png-accessible-cells-and-forced-limit}.
The limit $\delta\downarrow0$ is taken only in the next subsection.

\paragraph*{\textbf{The representation}}
For a nonnegative bridge array
$Z=(Z_{i,j}:1\leq i\leq k,\ i<j\leq\bar i)$, define
\begin{gather*}
 \mathrm dE_{Z,j}
 :=\sum_{\substack{1\leq i\leq k\\i<j\leq\bar i}}
 Z_{i,j}\,\delta_{-i},\qquad E_{Z,j}(0)=0,\\
 C^Z(i,q):=\mathcal L_{E_Z}\bigl((-i,i);(0,q)\bigr),
 \qquad 1\leq i\leq k,\quad 1\leq q\leq n.
\end{gather*}
Here $\mathcal L$ is the passage functional of
Definition~\ref{def:half-space-png-cadlag-lpp}, with marked bridge
coordinates $\{-k,\ldots,-1\}$. An impossible connection has value
$-\infty$.

For compatible matrices, write
\begin{equation*}
 (\mathsf U\odot\mathsf V)(p,q)
 :=\max_r\{\mathsf U(p,r)+\mathsf V(r,q)\}.
\end{equation*}
We also use juxtaposition for this max-plus product and its row-vector
version. For $1\leq p,q,j\leq n$, the transfer matrices for a
positive atom $u$ and a negative atom $-v$ on line $j$ are
\begin{align*}
 \mathsf P_j(u)(p,q)
 &:=\begin{cases}
 u\mathbf1_{\{p\leq j\leq q\}},&p\leq q,\\
 -\infty,&p>q,
 \end{cases}\\
 \mathsf N_j(v)(p,q)
 &:=\begin{cases}
 -v,&p=q=j,\\
 0,&p\leq q,\ (p,q)\neq(j,j),\\
 -\infty,&p>q.
 \end{cases}
\end{align*}
A positive atom can be collected whenever its line is crossed.
A negative atom can be avoided unless the entrance and exit lines
both equal $j$; this includes atoms at the right endpoint under the
one-sided convention in Definition~\ref{def:half-space-png-cadlag-lpp}.
For arrays $U,V$, set
\begin{align*}
 T_m:=\mathsf N_n(V_{n,m})\odot\cdots\odot\mathsf N_1(V_{1,m})
       \odot\mathsf P_1(U_{1,m})\odot\cdots\odot\mathsf P_n(U_{n,m}),\qquad
 T_{1:m}:=T_1\odot\cdots\odot T_m,\qquad T_{1:0}:=I,
\end{align*}
where $I$ is the max-plus identity, with diagonal entries zero and
all other entries $-\infty$.

\begin{theorem}[Finite line representation]
\label{thm:half-space-png-finite-staggered-straightening}
Fix $M\in\mathbb N$ and $a\in\mathbb Z$.
For every deterministic array \(X:\mathcal K_M\to\mathbb N_0\), where
\begin{equation*}
 \mathcal K_M:=\bigcup_{\substack{1\leq i\leq k\\0\leq m\leq M}}
 \{z\in\mathbb Z^2:s_i^a\leq z\leq v_m^M,\ z_1\geq z_2\},
\end{equation*}
and with \(X(s_i^a)=0\) for every \(i\), the finite rearrangement below produces arrays $Z^{\rm out},U,V$
and constants \(\kappa_i\), independent of \(m\), satisfying, for all \(1\leq i\leq k\) and \(0\leq m\leq M\),
\begin{gather}
\kappa_i:=L_X(s_i^a;v_0^M)-C^{Z^{\rm out}}(i,2k)
 \label{eq:half-space-png-straightening-row-gauge}\\
L_X(s_i^a;v_m^M)
 =\kappa_i+
 \bigl[C^{Z^{\rm out}}\odot T_{1:m}\bigr](i,2k).
 \label{eq:half-space-png-straightening-boundary-measurement}
\end{gather}

\smallskip
For the product-geometric array with source weights set to zero, all entries
of \(Z^{\rm out}\), \(U\), and \(V\) are independent, with
\begin{gather}
Z^{\rm out}_{i,j}\sim
 \begin{cases}
  \operatorname{Geo}(\gamma_j/g_i),&i<j<\bar i,\\
  \operatorname{Geo}(c/g_i),&j=\bar i,
 \end{cases}
 \quad 1\leq i\leq k,\quad i<j\leq\bar i,
 \label{eq:half-space-png-straightening-final-core-law}\\
V_{j,m}\sim\operatorname{Geo}(\delta/\gamma_j),
 \qquad U_{j,m}\sim\operatorname{Geo}(\delta\gamma_j),
 \quad 1\leq j\leq2k,\quad 1\leq m\leq M.
 \label{eq:half-space-png-straightening-final-packet-laws}
\end{gather}
Consequently, the geometric cut field has the joint representation
\begin{equation}
 \bigl(H_{i,M}^{\delta,a}(m):i\leq k,\ 0\leq m\leq M\bigr)
 \stackrel{\mathrm d}{=}
 \bigl([C^{Z^{\rm out}}\odot T_{1:m}](i,2k)
       -C^{Z^{\rm out}}(i,2k):i\leq k,\ 0\leq m\leq M\bigr).
 \label{eq:half-space-png-unconditional-exact-staggered-readout}
\end{equation}
\end{theorem}

\paragraph*{\textbf{The finite rearrangement}}
We now give the deterministic construction used in the theorem.
The bridge array is
\begin{equation*}
 Z_{i,j}:=X(a+\bar i,a+j),\qquad
 1\leq i\leq k,\quad i<j\leq\bar i.
\end{equation*}
The indices \((i,j)\) correspond to the marked line atom \(Z(-i,j)\)
in the Poisson-line construction.

For the fixed window size $M$, put
\begin{equation*}
 p_h:=a+2k+M-h+1,\qquad r_m:=a+2k+M+m,
 \qquad 1\leq h,m\leq M.
\end{equation*}
Thus \(v_0^M=(p_1,p_1)\) and \(v_m^M=(r_m,p_m-1)\) for \(1\leq m\leq M\).
For \(1\leq j\leq2k\) and \(1\leq h,m\leq M\), initialize
\begin{equation*}
 A_j^{(h)}:=X(p_h,a+j),\qquad
 D_h:=X(p_h,p_h),\qquad U_{j,m}:=X(r_m,a+j).
\end{equation*}  Under the geometric input law, these variables, the entries of $Z$,
and the remaining weights on $\mathcal K_M$ are independent,
with
\begin{equation*}
 \begin{aligned}
 Z_{i,j}&\sim\operatorname{Geo}(\gamma_j/g_i) &&(i<j<\bar i),&
 Z_{i,\bar i}&\sim\operatorname{Geo}(c/g_i),\\
 A_j^{(h)},U_{j,m}&\sim\operatorname{Geo}(\delta\gamma_j),&&&
 D_h&\sim\operatorname{Geo}(c\delta),\\
 X(p_\ell,p_h),X(r_m,p_s)&\sim\operatorname{Geo}(\delta^2)
 &&(\ell<h,\ m<s).&&
 \end{aligned}
\end{equation*}
For \(u,v,w\in\mathbb N_0\), set
\begin{equation}
 \mathscr R(u,v,w)
 :=\bigl(w+(u-v)_+,\,w+(v-u)_+,\,u\wedge v\bigr).
 \label{eq:half-space-png-finite-corner-map}
\end{equation}
This involution preserves $u+w$ and $v+w$. Each assignment below
replaces its three displayed variables by the components of
$\mathscr R$. The working variables $C$ and $w_{\mathrm{int}}$
carry the active coordinates through the updates; they are not
additional independent copies.

The construction has three stages: fold each background coordinate
$p_h$ through the bridge block, move its resulting negative factor
past the remaining $p$-coordinates, and interleave the positive
$r$-factors with these negative factors.

For \(h=M,M-1,\ldots,1\), apply the following updates.  For \(i=1,\ldots,k\), set \(C=A_i^{(h)}\), update in increasing \(j\), and then apply the last update:
\begin{gather}
(Z_{i,j},C,A_j^{(h)})\longleftarrow
 \mathscr R(Z_{i,j},C,A_j^{(h)}),
 \quad j=i+1,\ldots,\bar i-1,
 \label{eq:half-space-png-straightening-forward-fold-word}\\
(Z_{i,\bar i},C,D_h)\longleftarrow
 \mathscr R(Z_{i,\bar i},C,D_h),
 \qquad A_i^{(h)}\longleftarrow C.
 \label{eq:half-space-png-straightening-terminal-fold-letter}
\end{gather}
For \(i=k,k-1,\ldots,1\), start with \(C=A_{\bar i}^{(h)}\), apply
\begin{equation*}
 (Z_{i,j},A_{\bar j}^{(h)},C)\longleftarrow
 \mathscr R(Z_{i,j},A_{\bar j}^{(h)},C),
 \quad j=\bar i-1,\ldots,i+1,
\end{equation*}
in decreasing \(j\), and put \(A_{\bar i}^{(h)}\leftarrow C\).  Set
\begin{equation}
 V_{j,h}:=A_{\bar j}^{(h)},\qquad 1\leq j\leq2k.
 \label{eq:half-space-png-straightening-fold-output}
\end{equation}
Next, for \(\ell=h-1,h-2,\ldots,1\), start with \(w_{\mathrm{int}}=X(p_\ell,p_h)\) and apply
\begin{equation}
 (V_{j,h},A_j^{(\ell)},w_{\mathrm{int}})\longleftarrow
 \mathscr R(V_{j,h},A_j^{(\ell)},w_{\mathrm{int}}),
 \qquad j=2k,\ldots,1.
 \label{eq:half-space-png-straightening-low-low-word}
\end{equation}
After processing every \(p_h\), for \(m=1,\ldots,M-1\) and \(s=M,M-1,\ldots,m+1\), start with \(w_{\mathrm{int}}=X(r_m,p_s)\) and apply
\begin{equation}
 (U_{j,m},V_{j,s},w_{\mathrm{int}})\longleftarrow
 \mathscr R(U_{j,m},V_{j,s},w_{\mathrm{int}}),
 \qquad j=2k,\ldots,1.
 \label{eq:half-space-png-straightening-low-high-word}
\end{equation}
\paragraph*{\textbf{Deterministic passage identities}}
We first record two matrix identities used to verify the rearrangement. Put \(n=2k\) and write
\begin{equation*}
 \mathsf P(a)=\mathsf P_1(a_1)\odot\cdots\odot \mathsf P_n(a_n),\qquad
 \mathsf N(v)=\mathsf N_n(v_n)\odot\cdots\odot \mathsf N_1(v_1).
\end{equation*}
For nonnegative vectors, their entries are
\begin{equation*}
 \mathsf P(a)(p,q)=\sum_{j=p}^q a_j,
 \qquad
 \mathsf N(v)(p,q)=\begin{cases}-v_p,&p=q,\\0,&p<q.\end{cases}
\end{equation*}
As usual entries below the diagonal are $-\infty$.

First start with any $c_0\geq0$ and, in increasing order of $j$, set \((b_j,c_j,a'_j)=\mathscr R(z_j,c_{j-1},a_j)\). Then the following matrix identity holds:
\begin{equation}
 \label{eq:half-space-png-positive-row-interchange}
 \mathsf P(z)\odot\mathsf P(a)
 =\mathsf P(a')\odot\mathsf P(b).
\end{equation}
Fix $p\leq q$, put $c=c_{p-1}$, and set
\begin{equation*}
 Z_j=\sum_{\ell=p}^jz_\ell,\quad
 A_j=\sum_{\ell=p}^ja_\ell,\quad
 h_j=Z_j-A_{j-1},\quad
 H_j=\max_{p\leq r\leq j}h_r,\quad
 D_j=(H_j-c)_+,
\end{equation*}
with $A_{p-1}=D_{p-1}=0$. The recurrence for $c_j$ gives
\begin{equation*}
 \sum_{\ell=p}^ja'_\ell=Z_j-D_j,
 \qquad \sum_{\ell=p}^jb_\ell=A_j+D_j.
\end{equation*}
The left entry in \eqref{eq:half-space-png-positive-row-interchange} is
$A_q+H_q$ and the right entry is
\begin{equation*}
 A_q+D_q+\max_{p\leq r\leq q}
      \{h_r-D_r-D_{r-1}\}.
\end{equation*}
The last maximum is $\min\{H_q,c\}$: before the first crossing of
$c$ all $D$'s vanish, and at the first crossing the displayed expression
is exactly $c$; all subsequent expressions are at most $c$.
This also covers the case with no crossing.

Second, start with any $t_{n+1}\geq0$ and, in decreasing order of $j$,
set \((v'_j,a'_j,t_j)=\mathscr R(v_j,a_j,t_{j+1})\). Then
\begin{equation}
 \label{eq:half-space-png-negative-positive-interchange}
 \mathsf N(v)\odot\mathsf P(a)
 =\mathsf P(a')\odot\mathsf N(v').
\end{equation}
The diagonal entries agree because $a'_j-v'_j=a_j-v_j$.
For $p<q$ the two entries are respectively
\begin{equation*}
 \sum_{j=p+1}^q a_j+(a_p-v_p)_+,
 \qquad
 \sum_{j=p}^{q-1}a'_j+(a'_q-v'_q)_+.
\end{equation*}
They agree by
$a'_j=\min\{v_{j+1},a_{j+1}\}+(a_j-v_j)_+$ for $j<n$,
and $\min\{v_j,a_j\}+(a_j-v_j)_+=a_j$.
The initial value $t_{n+1}\geq0$ is arbitrary.
To track a terminal at the intersection being eliminated, we also need
an identity for the intermediate endpoint: writing \(w_{\mathrm{int}}=t_{n+1}\), for any incoming row \(\ell\),
\begin{equation}
 [\ell\mathsf P(a')]_n
 =w_{\mathrm{int}}+\max\{\ell_n,[\ell\mathsf N(v)\mathsf P(a)]_n\}.
 \label{eq:half-space-png-intermediate-interchange-endpoint}
\end{equation}
Indeed, for \(p<n\), \(\mathsf P(a')(p,n)
=w_{\mathrm{int}}+[\mathsf N(v)\mathsf P(a)](p,n)\), while its diagonal entry is \(w_{\mathrm{int}}+(a_n-v_n)_+\).

\begin{lemma}[Passage identities for the updates]
\label{lem:half-space-png-all-port-fold}
Let $n=2k$ and apply \eqref{eq:half-space-png-straightening-forward-fold-word}--%
\eqref{eq:half-space-png-straightening-fold-output} to one input
$(Z,A,D)$, with outputs $(Z',V)$; $A$ in the identities below
denotes its initial value.
Record $\alpha_i$ as the value of $A_i$ at the start of forward row $i$,
and $\beta_i$ as the value of $C$ at the end of backward row $i$
(before the subsequent updates). Put
$\chi_i=\alpha_i+\beta_i$.
Then, simultaneously for every source $i$ and every $q\geq i$,
\begin{align}
 [C^Z\odot\mathsf P(A)](i,q)
 &=\chi_i+[C^{Z'}\odot\mathsf N(V)](i,q),
 \label{eq:half-space-png-all-port-fold}\\
 D+[C^Z\odot\mathsf P(A)](i,n)
 &=\chi_i+C^{Z'}(i,n).
 \label{eq:half-space-png-diagonal-fold}
\end{align}
In particular, $\chi_i$ does not depend on $q$.
\end{lemma}

\begin{proof}
For $k=1$, write $Z_{1,2}=z$, $A=(a_1,a_2)$, and $D=d$.
The updates give
\begin{equation*}
 z'=d+(z-a_1)_+,\qquad V_1=a_2,\qquad
 V_2=d+(a_1-z)_+,\qquad \chi_1=a_1+a_2.
\end{equation*}
Since $C^Z=(0,z)$ and $C^{Z'}=(0,z')$, the all-port identity at
$q=1$ is $a_1=\chi_1-V_1$, and at $q=2$ it is
$a_2+\max\{a_1,z\}=\chi_1+(z-a_1)_+$.
Adding $d$ gives the diagonal identity.

Assume $k\geq2$ and the result for $k-1$. Write the first row of the bridge array as $z=(0,Z_{1,2},\ldots,Z_{1,n})$. For a source $i>1$, let $\ell_i^Z$
be its passage vector in the bridge array restricted to lines $2,\ldots,n-1$, extended by
$\ell_i^Z(1)=-\infty$ and $\ell_i^Z(n)=\ell_i^Z(n-1)$. Then
$C^Z(i,\cdot)=\ell_i^Z\odot\mathsf P(z)$.
The first source has $C^Z(1,\cdot)=e\odot\mathsf P(z)$, where
$e=(0,\ldots,0)$.

Apply the first forward updates through $j=2,\ldots,n-1$, with initial
value $C=A_1$. Write its first and third output components as $b_j$
and $a_j$, respectively, and let $x$ be the final value of $C$. Complete these
vectors by
\begin{equation*}
 a_1=0,\quad b_1=A_1,\quad
 d=a_n=\min\{Z_{1,n},x\},\quad
 b_n=A_n+(Z_{1,n}-x)_+.
\end{equation*}
For the positive-row identity, initialize its running variable at
zero. Its first update is
$\mathscr R(0,0,A_1)=(A_1,A_1,0)$, so the remaining updates start
at $C=A_1$ and give $a_1=0$, $b_1=A_1$.
Thus \eqref{eq:half-space-png-positive-row-interchange} gives
$\mathsf P(z)\mathsf P(A)=\mathsf P(a)\mathsf P(b)$.
The actual last update gives
\begin{equation*}
 Z'_{1,n}=D+(Z_{1,n}-x)_+,\qquad
 V_n=D+(x-Z_{1,n})_+,
\end{equation*}
and leaves $d$ as the diagonal input for the inner array.

Apply the induction hypothesis to the inner bridge array, the vector
$(a_2,\ldots,a_{n-1})$, and diagonal $d$. Denote its updated passage vector
by $\overline C_i$ and the weights of its negative atoms by $w_2,\ldots,w_{n-1}$.
Extend \(w_1=A_1,\qquad w_n=A_n+(x-Z_{1,n})_+\), and, for each $i>1$, set \(\overline C_i(1)=-\infty,\qquad \overline C_i(n)=\overline C_i(n-1)+A_n\). The two induction identities give the vector equality
\begin{equation*}
 \ell_i^Z\odot\mathsf P(a)
 =\chi_i+\overline C_i\odot\mathsf N(w).
\end{equation*}
Indeed, $\overline C_i$ is nondecreasing on its feasible coordinates up to
$n-1$, and $w_n\geq A_n$, so
\begin{equation*}
 [\overline C_i\mathsf N(w)](n)
 =\max\{\overline C_i(n-1),\overline C_i(n)-w_n\}
 =\overline C_i(n-1).
\end{equation*}
The last left coordinate is $d+[\ell_i^Z\mathsf P(a)](n-1)$,
which agrees by the inner diagonal identity. The other coordinates are
the interface induction hypothesis for $q<n$.

Now commute $\mathsf N(w)\mathsf P(b)$ using
\eqref{eq:half-space-png-negative-positive-interchange}, with initial value
$t_{n+1}=D$. The $n$-th update produces exactly $(V_n,Z'_{1,n})$
and leaves $A_n$ as the next initial value. The updates for $j=n-1,\ldots,2$ are exactly
the first backward row. At $j=1$, since $w_1=b_1=A_1$, both output
components equal $V_1$. Thus the resulting vector of positive atom weights is \(b'=(V_1,Z'_{1,2},\ldots,Z'_{1,n})\), and the negative vector is $V$.
For $i>1$, $\overline C_i\mathsf P(b')=C^{Z'}(i,\cdot)$.
Before the last coordinate this is immediate; at the last,
the entrance value $\overline C_i(n-1)+A_n$ is dominated
by the route through line $n-1$, because the first backward update
has $Z'_{1,n-1}\geq A_n$. This proves
\eqref{eq:half-space-png-all-port-fold} for all inner sources.

For the first source, put $\zeta_j=\sum_{\ell=2}^jZ_{1,\ell}$ and
$H_j=[C^Z\mathsf P(A)](1,j)$. The scalar path recurrence is
$H_j=A_j+\max\{H_{j-1},\zeta_j\}$, and its auxiliary variable is
$c_j=H_j-\zeta_j$. Thus $H_1=A_1$ and
$H_j-H_{j-1}=A_j+(Z_{1,j}-c_{j-1})_+=b_j$, giving
\begin{equation*}
 [C^Z\mathsf P(A)](1,q)=A_1+\sum_{j=2}^q b_j
 =A_1+[e\mathsf N(w)\mathsf P(b)](q).
\end{equation*}
Here $e\mathsf N(w)=(-A_1,0,\ldots,0)$.
After the same interchange,
$e\mathsf P(b')=V_1+C^{Z'}(1,\cdot)$, proving the first-row identity
with $\chi_1=A_1+V_1$. The constants for the inner rows are unchanged on adding lines $1$ and $n$,
so they are precisely the recorded $\alpha_i+\beta_i$.

Finally every row of the bridge array satisfies
$C^{Z'}(i,n)=C^{Z'}(i,n-1)+Z'_{1,n}$, since the only atom on the
last line of the bridge environment is $Z'_{1,n}$. Therefore
\begin{equation*}
 [C^{Z'}\mathsf N(V)](i,n)
 =C^{Z'}(i,n)-\min\{V_n,Z'_{1,n}\}
 =C^{Z'}(i,n)-D,
\end{equation*}
which proves \eqref{eq:half-space-png-diagonal-fold}.
\end{proof}

\begin{proof}[Proof of Theorem~\ref{thm:half-space-png-finite-staggered-straightening}]
\emph{Step 1: passage preservation.}
We maintain the following invariant: at every stage, all
source-to-interface passage values and all retained terminal values
are preserved up to one additive constant for each source.
The folding identities introduce source constants; the matrix
interchanges preserve the passage vectors, and the intermediate-endpoint
identity tracks a terminal at an eliminated intersection.
Max-plus multiplication carries an additive source constant unchanged.
In the original grid the passage matrix for the special
block is $C^Z$: its row at label $\bar i$ has exactly
the atoms $Z_{i,j}$, and its source cell is zero.  An ordinary row
with special-block weights $a$ contributes the matrix $\mathsf P(a)$.
The scalar recurrence is $H_j=a_j+\max\{H_{j-1},\ell_j^{\mathrm{in}}\}$, where
$\ell_j^{\mathrm{in}}$ is the incoming passage value at line $j$; hence these identities
hold simultaneously for every source.

The initial background order is
$p_M,\ldots,p_1,r_1,\ldots,r_M$.
Eliminate the rows indexed by $p_h$ in the order $h=M,\ldots,1$.
Lemma~\ref{lem:half-space-png-all-port-fold} replaces the
block consisting of the special coordinates and $p_h$ by the updated bridge array followed by
$\mathsf N(V_{\cdot,h})$, with one constant added to each source row.
Its second identity gives the same constants at the diagonal endpoint;
thus the constants are the same for every retained terminal.
Move this negative factor past the remaining positive rows
$p_{h-1},\ldots,p_1$ using
\eqref{eq:half-space-png-negative-positive-interchange}, initialized
at their intersection $w_{\mathrm{int}}=X(p_\ell,p_h)$.  The updates in decreasing $j$ are
exactly \eqref{eq:half-space-png-straightening-low-low-word}.
In particular each interchange preserves every exit coordinate of every
source row, including the coordinates needed at later interchanges.
Descending induction in $h$ gives the negative-factor order
$p_1,\ldots,p_M$. The constants from successive folds accumulate
into one constant $\kappa_i$ for each source.

Next move the positive factor for $r_m$ left past
$p_M,\ldots,p_{m+1}$, for $m=1,\ldots,M-1$.
The same identity applies with $a=U_{\cdot,m}$,
$v=V_{\cdot,s}$, and initial value $X(r_m,p_s)$.
Swapping the first two inputs of $\mathscr R$ swaps its first two
outputs, so the order of variables is precisely
\eqref{eq:half-space-png-straightening-low-high-word}.
The final factor order is
\begin{equation}
 p_1,r_1,p_2,r_2,\ldots,p_M,r_M.
 \label{eq:half-space-png-straightening-final-wire-word}
\end{equation}
All original grid weights are included in this composition.  The bridge
array contains the weights from the special block.  The updates for
$p_h$ use its intersections with that block and its diagonal weight;
the first interchanges use the $(p_\ell,p_h)$ cells, and $U_{j,m}$
gives the intersections of $r_m$ with the special block.  A mixed cell is in the cone of
terminal $\ell$ exactly when \((r_m,p_s)\leq v_\ell^M \quad\Longleftrightarrow\quad m\leq\ell<s\). These are exactly the second interchanges.  No cell $(r_m,r_s)$ is accessible.  Thus this composition is also a direct elimination of
the grid Bellman recurrences: the forward and backward identities eliminate one
block consisting of the special coordinates and $p_h$, and the row-interchange identity eliminates
one remaining intersection with its displayed initial value.
Each background intersection is eliminated once and every retained source remains.
At an intersection of weight \(w_{\mathrm{int}}\),
\eqref{eq:half-space-png-intermediate-interchange-endpoint} is the
Bellman recurrence for the endpoint between the exchanged factors.
It therefore gives the passage values at retained terminals there, as well
as those after both factors.

\emph{Step 2: terminal identification and pinning.}
The terminal $v_0^M$ is the end of the updated bridge array, by the diagonal
identity \eqref{eq:half-space-png-diagonal-fold} at each step.
Terminal $v_m^M$ is the end of the first $m$ pairs in
\eqref{eq:half-space-png-straightening-final-wire-word}, since its
background labels are $(r_m,p_m-1)$.  The factors indexed by $p_m$ and $r_m$ are
$\mathsf N(V_{\cdot,m})$ and $\mathsf P(U_{\cdot,m})$,
respectively; their product is $T_m$.
All the preceding replacements preserve every passage-matrix entry,
so they give, with the same row constants for every \(m\),
\begin{equation*}
 L_X(s_i^a;v_m^M)=\kappa_i+[C^{Z^{\rm out}}\odot T_{1:m}](i,2k).
\end{equation*}
The case $m=0$ determines
\eqref{eq:half-space-png-straightening-row-gauge}.

\emph{Step 3: the geometric product law.}
Write \((u',v',w')=\mathscr R(u,v,w)\).  For independent inputs with respective laws $\operatorname{Geo}(r)$,
$\operatorname{Geo}(s)$, and $\operatorname{Geo}(rs)$, where
$0\leq r,s<1$, the outputs have the same product law.
For $r,s>0$, use bijectivity and
\begin{equation*}
 r^u s^v(rs)^w=r^{u+w}s^{v+w}
 =r^{u'+w'}s^{v'+w'}=r^{u'}s^{v'}(rs)^{w'};
\end{equation*}
the cases \(r=0\) or \(s=0\) follow directly from
\eqref{eq:half-space-png-finite-corner-map}.

The forward, final, and backward updates use parameter triples \((r,s,rs)\):
\begin{equation*}
 \left(\frac{\gamma_j}{g_i},\delta g_i,\delta\gamma_j\right),
 \qquad
 \left(\frac c{g_i},\delta g_i,c\delta\right),
 \qquad
 \left(\frac{\gamma_j}{g_i},\frac\delta{\gamma_j},
       \frac\delta{g_i}\right),
\end{equation*}
while the two sequences of background updates use \((\delta/\gamma_j,\delta\gamma_j,\delta^2)\) or the same first two
parameters in reverse order. At each step, the three active variables
have the displayed product law and are independent of all untouched
active variables. The involution therefore preserves the product law
of the entire active array. Induction through the finite sequence of
updates proves mutual independence of all retained outputs and
\eqref{eq:half-space-png-straightening-final-core-law}--%
\eqref{eq:half-space-png-straightening-final-packet-laws}.
Combining the deterministic pinned identity with the source-cancellation
limit in Lemma~\ref{lem:half-space-png-accessible-cells-and-forced-limit}
gives \eqref{eq:half-space-png-unconditional-exact-staggered-readout}.
\end{proof}

\paragraph*{\textbf{Staggered atoms and mesh approximation}}
We first realize the ordered product $T_m$ exactly by signed line atoms. Take the bridge array and families
$(U_{j,m},V_{j,m})_{j=1}^{2k}$,
$m\geq1$, independently with the laws
\eqref{eq:half-space-png-straightening-final-core-law}--%
\eqref{eq:half-space-png-straightening-final-packet-laws}.  For \(m\geq1\), place the atoms with index \(m\) at
\begin{equation*}
 x_{m,j}^-:=(m-1)\delta+\frac{2k+1-j}{4k+1}\delta,
 \qquad
 x_{m,j}^+:=(m-1)\delta+\frac{2k+j}{4k+1}\delta.
\end{equation*}
Their left-to-right order is \(-V_{2k,m},\ldots,-V_{1,m},U_{1,m},\ldots,U_{2k,m}\).  Set
\begin{equation*}
 \mathrm d\widehat S_j^\delta
 :=-\sum_{m\geq1}V_{j,m}\delta_{x_{m,j}^-}
   +\sum_{m\geq1}U_{j,m}\delta_{x_{m,j}^+},
\end{equation*}
normalize \(\widehat S_j^\delta(0)=0\), retain the bridge atoms at
\(-k,\ldots,-1\), and call the resulting line environment
\(\widehat E^\delta\).  Define
\begin{equation*}
 \widehat G_i^\delta(x)
 :=\mathcal L_{\widehat E^\delta}(s_i;(x,2k)),
 \qquad
 \widehat h_i^\delta(x)
 :=\widehat G_i^\delta(x)-\widehat G_i^\delta(0).
\end{equation*}
Theorem~\ref{thm:half-space-png-finite-staggered-straightening} gives
\begin{equation}
 \bigl(H_{i,M}^{\delta,a}(m):i\leq k,\ 0\leq m\leq M\bigr)
 \stackrel{\mathrm d}{=}
 \bigl(\widehat h_i^\delta(\delta m):
       i\leq k,\ 0\leq m\leq M\bigr).
 \label{eq:half-space-png-exact-staggered-readout}
\end{equation}

Next move all atoms in $((m-1)\delta,m\delta)$ to $m\delta$.
Keep $U$ and $V$ as separate marked inputs; the resulting c\`adl\`ag
line has their net increment $U_{j,m}-V_{j,m}$ at each mesh point.  With the same bridge array and the resulting
environment \(\widetilde E^\delta\), set
\begin{equation*}
\widetilde S_j^\delta(x)
 :=\sum_{1\leq m\leq\lfloor x/\delta\rfloor}(U_{j,m}-V_{j,m}), \qquad
\widetilde G_i^\delta(x)
 :=\mathcal L_{\widetilde E^\delta}(s_i;(x,2k)),
 \qquad
 \widetilde h_i^\delta(x)
 :=\widetilde G_i^\delta(x)-\widetilde G_i^\delta(0).
 \end{equation*}

This collapse need not preserve passage values when several atoms
occur in one interval. For example, with two lines,
$U_{2,m}=V_{2,m}=1$ and all other increments zero, a path from line
$1$ to line $2$ can avoid the staggered negative atom and collect the
positive atom, obtaining weight $1$; after collapse the net jump is
zero. The lemma below controls all such exceptional intervals.

\begin{figure}[t]
\centering
\begin{tikzpicture}[
 x=1cm,y=0.64cm,
 every node/.style={font=\footnotesize},
 >={Latex[length=1.4mm,width=1.1mm]},
 carrier/.style={gray!50,line width=0.45pt},
 plus/.style={green!40!black,line width=0.9pt,->},
 minus/.style={orange!85!black,line width=0.9pt,->}
]
\node at (2.25,5.15) {(a) Staggered packet};
\fill[black!2] (0,0.68) rectangle (4.5,4.40);
\foreach \j in {1,2,3,4}{
 \draw[carrier] (0,\j)--(4.5,\j);
 \node[anchor=east] at (-0.10,\j) {\(\j\)};
 \pgfmathsetmacro{\xm}{(5-\j)/9*4.5}
 \pgfmathsetmacro{\xp}{(4+\j)/9*4.5}
 \draw[minus] (\xm,\j+0.17)--(\xm,\j-0.17);
 \draw[plus] (\xp,\j-0.17)--(\xp,\j+0.17);
}
\draw[gray!70,dashed] (0,0.65)--(0,4.42);
\draw[gray!70,dashed] (4.5,0.65)--(4.5,4.42);
\node[anchor=north] at (0,0.51) {\((m-1)\delta\)};
\node[anchor=north] at (4.5,0.51) {\(m\delta\)};
\draw[->,black!65] (4.80,2.45)--(5.45,2.45);
\begin{scope}[xshift=6cm]
 \node at (1.15,5.15) {(b) Mesh collapse};
 \fill[black!2] (0,0.68) rectangle (2.3,4.40);
 \foreach \j in {1,2,3,4}{
  \draw[carrier] (0,\j)--(2.3,\j);
  \fill[green!40!black] (1.85,\j+0.13) circle[radius=1.5pt];
  \fill[orange!85!black] (1.85,\j-0.13) circle[radius=1.5pt];
 }
 \draw[gray!70,dashed] (1.85,0.65)--(1.85,4.42);
 \node[anchor=north] at (1.85,0.51) {\(m\delta\)};
 \node at (1.15,-0.55) {retain line and sign};
\end{scope}
\draw[->,black!65] (8.65,2.45)--(9.35,2.45);
\node[anchor=south] at (9,2.65) {\(\delta\downarrow0\)};
\begin{scope}[xshift=9.9cm]
 \node at (2.10,5.15) {(c) Signed Poisson lines};
 \fill[black!2] (0,0.68) rectangle (4.2,4.40);
 \foreach \j in {1,2,3,4}
  \draw[carrier] (0,\j)--(4.2,\j);
 \foreach \xx/\yy in {0.35/1,2.25/1,1.05/2,3.70/2,1.65/3,3.20/3,0.70/4,2.75/4}
  \draw[plus] (\xx,\yy-0.16)--(\xx,\yy+0.16);
 \foreach \xx/\yy in {1.45/1,3.55/1,2.65/2,0.55/3,2.10/4,3.85/4}
  \draw[minus] (\xx,\yy+0.16)--(\xx,\yy-0.16);
 \draw[->,black!65] (0,0.5)--(4.35,0.5);
 \node[anchor=north] at (0,0.46) {\(0\)};
 \node[anchor=west] at (4.35,0.5) {\(x\)};
 \node at (2.1,-0.55)
  {\(\text{line }j:\quad +:\gamma_j,\quad -:\gamma_j^{-1}\)};
\end{scope}
\end{tikzpicture}
\caption{From geometric packets to Poisson lines, illustrated for \(k=2\).
The exact packet order is
\(-V_{4,m},\ldots,-V_{1,m},U_{1,m},\ldots,U_{4,m}\).
Moving atoms to \(m\delta\) retains their line and sign; equality of mesh
passage values holds outside the exceptional cells of
Lemma~\ref{lem:half-space-png-compact-window-collapse}.
The limiting positive and negative atoms are independent Poisson
processes.  The geometric bridge array is retained throughout.}
\label{fig:half-space-png-packet-poissonization}
\end{figure}
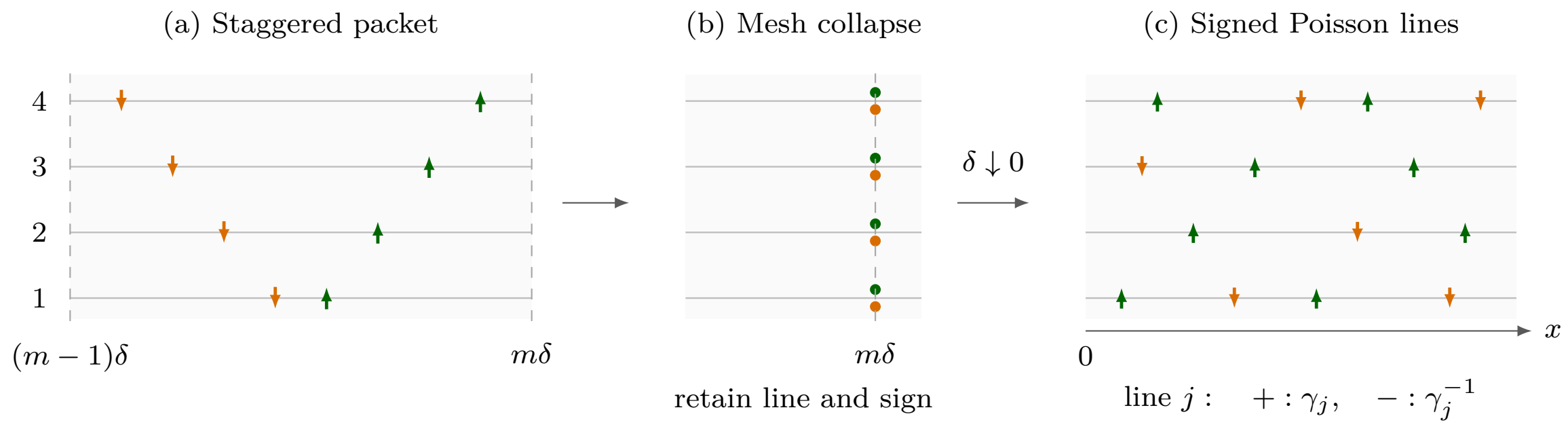

Figure~\ref{fig:half-space-png-packet-poissonization} shows the
mesh collapse and the subsequent Poisson limit.

\begin{lemma}[Approximation on the mesh]
\label{lem:half-space-png-compact-window-collapse}
For every $0<L<\infty$, with \(M_\delta(L):=\lceil L/\delta\rceil+1\),
\begin{equation*}
 \mathbb P\left(
  \widehat h_i^\delta(\delta m)\neq
  \widetilde h_i^\delta(\delta m)
  \text{ for some }i\leq k,\ 0\leq m\leq M_\delta(L)
 \right)
 \leq C_{\boldsymbol g,L}\delta.
\end{equation*}
\end{lemma}

\begin{proof}
Call the $m$-th interval exceptional if
$B_m:=\{\sum_{j=1}^{2k}(U_{j,m}+V_{j,m})\geq2\}$ occurs.  Otherwise its
passage value from line \(p\) to line \(q\), \(p\leq q\), is
\begin{equation*}
 \begin{array}{c|c}
  \text{atoms in the interval}&\text{passage value}(p,q)\\ \hline
  \text{empty}&0\\
  \text{one \(+1\) atom on line \(j\)}
    &\mathbf1_{\{p\leq j\leq q\}}\\
  \text{one \(-1\) atom on line \(j\)}
    &-\mathbf1_{\{p=q=j\}}.
 \end{array}
\end{equation*}
Moving these atoms to the right endpoint preserves this value: for a negative atom
and \(p<q\), a switch avoids it, including at the endpoint under the
one-sided convention.  Composition preserves all source-to-terminal passage
values through successive nonexceptional intervals, before and after pinning.

Let \(\Lambda_{\boldsymbol g}:=
\sum_{j=1}^{2k}(\gamma_j+\gamma_j^{-1})\).  Since \(\mathbb P(\operatorname{Geo}(q)\geq1)=q\) and \(\mathbb P(\operatorname{Geo}(q)\geq2)=q^2\), independence and a union bound give
\begin{equation*}
 \mathbb P(\text{a fixed interval is exceptional})
 \leq\delta^2\left[
  \sum_{j=1}^{2k}(\gamma_j^2+\gamma_j^{-2})
  +\frac12\Lambda_{\boldsymbol g}^{\,2}\right]
 =:C_{\boldsymbol g}\delta^2.
\end{equation*}
A final union bound over the first \(M_\delta(L)\) intervals yields \(C_{\boldsymbol g}(L+2)\delta\) for \(0<\delta\leq1\), proving the
claim.
\end{proof}

For a sequence \(f:\mathbb N_0\to\mathbb R\), let \(\mathcal I_\delta f\)
be the upper-semicontinuous regularization of the function equal to
\(f(m)\) on \([\delta m,\delta(m+1))\).  For a function on
\(\delta\mathbb N_0\), use the same notation after identifying its
value at \(\delta m\) with the sequence entry at \(m\); apply
\(\mathcal I_\delta\) componentwise to vectors of profiles.
The collapsed fields are constant on these intervals, and therefore
$ \mathcal I_\delta\bigl(\widetilde h_i^\delta|_{\delta\mathbb N_0}\bigr)
 =\operatorname{usc}\widetilde h_i^\delta$.
Using the same infinite family $(Z,U,V)$, construct the coupling by setting
$H_i^\delta(m):=\widehat h_i^\delta(\delta m)$ for $m\in\mathbb N_0$.
The exact readout identity for every $M$ identifies the entire sequence
law as $H^\delta\sim\mu_{c;\boldsymbol g}^{\delta}$.
On the event that the first $M_\delta(L)$ intervals are nonexceptional,
the mesh values agree with those of the collapsed profiles. The extra
mesh interval also covers the adjacent plateau values used by
upper-semicontinuous regularization at $L$. Hence
\begin{equation}
 \mathbb P\left(
  \mathcal I_\delta H_i^\delta
  =\operatorname{usc}\widetilde h_i^\delta
  \text{ on }[0,L]\text{ for every }i\leq k
 \right)
 \geq1-C_{\boldsymbol g,L}\delta.
 \label{eq:half-space-png-unconditional-line-approximation}
\end{equation}

\subsubsection{Passage to the common-noise PNG evolution}
\label{subsubsec:half-space-png-poissonization}

We transfer joint stationarity from the geometric model to PNG by
passing the initial profiles and their common-noise evolution to the
limit together. Throughout, $c>0$ and
$g_1>\cdots>g_k>c\vee1$. Fix a deterministic $t\geq0$ and choose
$n_\delta\in\mathbb N_0$ with $t_\delta:=n_\delta\delta\to t$.
Let $H^\delta\sim\mu_{c;\boldsymbol g}^\delta$ and write
$R=(R_i^{(c;\boldsymbol g)})_{i=1}^k$ for the Poisson-line candidate.

For each time step, sample a fresh noise array $\xi_r$, independently
of the input and the other time steps, and use that same array for
all $k$ coordinates. We use the shorthand
\begin{equation*}
 (\Psi_\delta^{(c),[k]})^n
 :=\Psi_{\delta,\xi_n}^{(c),[k]}\circ\cdots\circ
   \Psi_{\delta,\xi_1}^{(c),[k]},
\end{equation*}
with the analogous convention for one coordinate and the identity
map when $n=0$. Our main claim is
\begin{equation}
 \left(\mathcal I_\delta H^\delta,
       \mathcal I_\delta(\Psi_\delta^{(c),[k]})^{n_\delta}H^\delta\right)
 \Longrightarrow
 \left(R,\Theta_t^{(c),[k]}(R;\Pi)\right),
 \label{eq:half-space-png-joint-input-future-convergence}
\end{equation}
in the product local-hypograph topology. Here $\Pi$ is independent
of $R$ and is shared by every coordinate. We first prove convergence
of the initial profiles, then of the independent future environment,
and finally of the initial and evolved profiles jointly, obtaining
\eqref{eq:half-space-png-joint-input-future-convergence}. We then
transfer stationarity and parameter-deletion consistency to PNG.

\paragraph*{\textbf{Convergence of the initial profiles}}
Recall the independent line increments from
\eqref{eq:half-space-png-straightening-final-packet-laws}, and set
\begin{equation*}
 \mathcal N_{j,+}^\delta:=\sum_{m\geq1}U_{j,m}\,\delta_{m\delta},
 \qquad
 \mathcal N_{j,-}^\delta:=\sum_{m\geq1}V_{j,m}\,\delta_{m\delta},
 \qquad 1\leq j\leq2k.
\end{equation*}
For $Y\sim\operatorname{Geo}(q)$, as $q\downarrow0$,
\begin{equation}
 \log\mathbb E e^{-uY}
 =\log\frac{1-q}{1-qe^{-u}}
 =-q(1-e^{-u})+O(q^2),
 \label{eq:half-space-png-geometric-laplace-expansion}
\end{equation}
uniformly for bounded $u\geq0$. Apply this expansion with
$q=\delta\gamma_j$ and $q=\delta/\gamma_j$, insert a nonnegative
continuous compactly supported test function in place of $u$, and
sum over the mesh sites. The resulting Riemann sums give joint vague
convergence
\begin{equation*}
 (\mathcal N_{j,+}^\delta,\mathcal N_{j,-}^\delta)_{j=1}^{2k}
 \Longrightarrow
 (\mathcal N_{j,+},\mathcal N_{j,-})_{j=1}^{2k},
 \qquad
 \begin{matrix}
  \operatorname{rate}(\mathcal N_{j,+})=\gamma_j,\\
  \operatorname{rate}(\mathcal N_{j,-})=\gamma_j^{-1}.
 \end{matrix}
\end{equation*}
The limiting point measures are independent. The finite bridge array
is independent of these measures and retains the laws
\eqref{eq:half-space-png-bridge-laws-spectral}.
We keep the two signs separate: convergence of their signed difference
alone could conceal cancellations at coincident prelimit sites.

\begin{lemma}[Continuity of the line profiles]
\label{lem:half-space-png-compact-line-network-continuity}
Fix $0<L<K<\infty$ and finitely many sources. Suppose that the
finite bridge arrays converge coordinatewise and that the positive
and negative counting measures on every line converge separately
and vaguely on $[0,K]$. If the limiting measures consist of finitely
many unit atoms, with distinct locations across all lines and signs
and no atom at $0,L,K$, then the associated pinned last-passage
profiles converge jointly on $[0,L]$ in compact hypograph topology.
\end{lemma}

\begin{proof}
Integer-valued marked convergence gives an eventual matching of all
ordinary atoms, preserving their line, sign, and spatial order.
Choose increasing piecewise affine homeomorphisms $\phi_\ell$,
equal to the identity outside $(0,K)$ and fixing $L$, which map each
limiting atom to its matched prelimit atom. Then
$\|\phi_\ell-\mathrm{id}\|_{[0,K]}\to0$.
These maps preserve switching patterns and atom ownership.
If $b_\ell$ is the sum of the absolute bridge-weight errors,
corresponding path weights differ by at most $b_\ell$. Maximization
and subtraction of the value at zero therefore give
\begin{equation*}
 \max_i\sup_{0\leq x\leq L}
 \left|\widetilde R_{i,\ell}(\phi_\ell(x))
             -\widetilde R_i(x)\right|
 \leq2b_\ell\longrightarrow0,
\end{equation*}
where the tildes denote the pinned c\`adl\`ag profiles.
Upper-semicontinuous regularization commutes with $\phi_\ell$,
which proves the claim.
\end{proof}

Poisson configurations satisfy these hypotheses almost surely.
The continuous mapping theorem and the coupling estimate
\eqref{eq:half-space-png-unconditional-line-approximation} give
\begin{align}
 (\operatorname{usc}\widetilde h_i^\delta)_{i=1}^k
 &\Longrightarrow R,
 \label{eq:half-space-png-collapsed-line-law-convergence}\\
 \mathcal I_\delta H^\delta
 &\Longrightarrow R,
 \qquad H^\delta\sim\mu_{c;\boldsymbol g}^\delta,
 \label{eq:half-space-png-antidiagonal-line-law-convergence}
\end{align}
locally in product hypograph topology.

\paragraph*{\textbf{Convergence of the fresh future environment}}
Sample all fresh update arrays independently of the entire initial
line coupling. The time-zero cut belongs to the past; the next
$2n_\delta$ lattice layers form the half-open future slab
$(0,t_\delta]$. Under the embedding
\eqref{eq:half-space-png-null-coordinate-embedding}, push one
representative of each off-diagonal reflection orbit, with its
multiplicity, to its physical point. Do the same for the diagonal
weights, obtaining $\Pi_{\rm bulk}^\delta$ and $\Pi_\partial^\delta$.

The bulk lattice has covolume $\delta^2/2$ and geometric parameter
$\delta^2$; the boundary lattice has time spacing $\delta$ and
geometric parameter $c\delta$.
Thus \eqref{eq:half-space-png-geometric-laplace-expansion}, with
$q=\delta^2$ and $q=c\delta$, gives
\begin{align*}
 \Pi_{\rm bulk}^\delta&\Longrightarrow\Pi_{\rm bulk},
 &\Pi_{\rm bulk}&\sim
 \operatorname{PPP}(2\,\mathrm ds\,\mathrm dz)
       |_{(0,t]\times(0,\infty)},\\
 \Pi_\partial^\delta&\Longrightarrow\Pi_\partial,
 &\Pi_\partial&\sim\operatorname{PPP}(c\,\mathrm ds)|_{(0,t]}.
\end{align*}
Here $\operatorname{PPP}(\mu)$ denotes a Poisson point measure of
intensity $\mu$. The limits are independent of one another and of
the initial line network. On a compact slab, the probability of any
weight exceeding one is $O(\delta^2)$ in the bulk and $O(\delta)$
at the boundary. Embed the moving slabs in a fixed compact envelope,
including the time faces, and extend their measures by zero.
The limiting measures have no atoms at times $0,t$, so the changing
terminal time causes no restriction discontinuity.

\paragraph*{\textbf{Convergence of the evolution}}
We first express the discrete dynamics through last passage in the
fresh slab. Let $\mathsf L_{\delta,n}^{(c)}(\ell,m)$ be its passage
weight from
\begin{equation*}
 (N+\ell,N-\ell)\quad\text{to}\quad
 (N+n+m,N+n-m)
\end{equation*}
in $p\geq q$, where the initial cut is at level $N$.
Count only weights strictly after the initial cut, including the
terminal weight. Set the value to $-\infty$ when $|\ell-m|>n$;
for $n=0$, set it to zero if $\ell=m$ and $-\infty$ otherwise.
For pinned $f:\mathbb N_0\to\mathbb Z$, define
\begin{equation*}
 (\mathcal T_{\delta,n}^{(c)}f)(m)
 :=\max_{\substack{\ell\geq0\\|\ell-m|\leq n}}
       \{f(\ell)+\mathsf L_{\delta,n}^{(c)}(\ell,m)\}
   -\max_{0\leq\ell\leq n}
       \{f(\ell)+\mathsf L_{\delta,n}^{(c)}(\ell,0)\}.
\end{equation*}
The Bellman recurrences
\eqref{eq:half-space-png-discrete-png-half-step}--%
\eqref{eq:half-space-png-discrete-png-update} and passage composition
give the pathwise identity
\begin{equation}
 \mathcal T_{\delta,n}^{(c)}f=(\Psi_\delta^{(c)})^n f.
 \label{eq:half-space-png-prelimit-bellman-composition}
\end{equation}
The initial heights are supplied by $f$, so initial-cut weights are
not counted again. Intermediate recentering constants cancel because max-plus updates
commute with addition of constants.

The coupling from the preceding subsection also identifies the
\emph{raw} mesh heights. For $K<\infty$, let
\begin{equation*}
 \mathcal E_{\delta,K}:=
 \left\{H_i^\delta(m)=\widetilde h_i^\delta(m\delta)
  \text{ for all }i\leq k,\quad
  0\leq m\leq\lceil K/\delta\rceil+1\right\}.
\end{equation*}
The mesh-approximation lemma gives
\begin{equation}
 \mathbb P(\mathcal E_{\delta,K})
 \geq1-C_{\boldsymbol g,K}\delta.
 \label{eq:half-space-png-raw-grid-coupling-probability}
\end{equation}
This raw equality is needed because the lattice update acts on mesh
heights, whereas upper-semicontinuous regularization can change their
values at downward jumps. If $K>L+t+1$, finite propagation and
\eqref{eq:half-space-png-prelimit-bellman-composition} give identical
interpolated outputs on $[0,L]$ for these two inputs, on
$\mathcal E_{\delta,K}$ and for all sufficiently small $\delta$,
when they use the same future slab.

We next state the continuity property used for these line inputs.
For a future point $p=(s,z)$, put $u(p)=s+z$ and $v(p)=s-z$.
Let $A\subset(0,K)$ be the set of ordinary incoming atom locations
on every line and of either sign, and let $P$ be the finite future
configuration in the compact envelope. Call $(A,P)$ \emph{generic}
if the configurations are simple, $P$ has no atom at times $0$ or $t$,
distinct future atoms have distinct $u$- and distinct $v$-coordinates,
and the random entries in the critical sets
\begin{equation*}
 \begin{split}
  \mathscr C_0&:=A\cup\{u(p),-v(p):p\in P\}\cup\{0,K\},\\
  \mathscr C_t&:=\{a-t,a+t:a\in A\}
    \cup\{u(p)-t,t-v(p):p\in P\}\cup\{0,L\}
 \end{split}
\end{equation*}
do not coincide with one another or with a deterministic endpoint,
apart from forced repetitions arising from a single atom or $t=0$.
Only entries in the corresponding spatial window are relevant.
These conditions stabilize causal order, separate the feasible initial
interval endpoints from incoming jumps, and prevent an output
breakpoint at zero or $L$. Independent Poisson inputs are almost
surely generic for each fixed deterministic $t,L,K$.

\begin{lemma}[Continuity of the evolution on compact intervals]
\label{lem:half-space-png-compact-cone-evolution-continuity}
Fix $L>0$, $t\geq0$, and $K>L+t$. Suppose that the bridge arrays
and separately marked incoming line measures converge on $[0,K]$
as in Lemma~\ref{lem:half-space-png-compact-line-network-continuity}.
Suppose also that the future point measures in $(0,t_\ell]$,
$t_\ell\to t$, including their multiplicities and bulk or boundary
marks, converge vaguely on a fixed compact envelope containing the
relevant backward cones and its time faces. If the limiting
configuration is generic, the pinned initial and evolved profiles
converge jointly on $[0,L]$ in compact hypograph topology.

In the lattice case, use raw mesh inputs
$f_{i,\ell}(m)=\widetilde R_{i,\ell}(m\delta_\ell)$ and outputs
$\mathcal I_{\delta_\ell}\mathcal T_{\delta_\ell,n_\ell}^{(c)}f_{i,\ell}$,
where $\delta_\ell\downarrow0$ and $n_\ell\delta_\ell=t_\ell$.
In the continuum case, use the half-space PNG max-plus evolution.
\end{lemma}

\begin{proof}
Initial-profile convergence follows from the preceding continuity
lemma. First assume $t>0$. There are finitely many causal chains in
the limiting future. For each such chain
$C=(p_1\prec\cdots\prec p_r)$, with $p_j=(s_j,z_j)$, define
\begin{equation*}
 b_{i,C}:=r+\sup_{(z_1-s_1)_+\leq y\leq z_1+s_1}R_i(y),
 \qquad
 I_C:=\bigl[(z_r-(t-s_r))_+,\ z_r+t-s_r\bigr].
\end{equation*}
The unpinned output is the maximum of these constant-height
contributions $b_{i,C}$ on their intervals $I_C$, together with the
no-future-atom contribution $\sup_{y \geq 0, \, |y-x| \leq t} R_i(y)$.

Marked convergence matches all future atoms. Their strict
null-coordinate comparisons are stable, so the collection of feasible
chains eventually agrees with that of the limit. Genericity of
$\mathscr C_0$ separates each feasible starting-interval endpoint
from the incoming jumps. Thus the corresponding initial suprema
converge; in the lattice case, the mesh captures every relevant
plateau because its limiting length is positive. Consequently,
the chain heights converge to $b_{i,C}$ and their support endpoints
converge to those of $I_C$.
The no-atom contribution is likewise a finite maximum over incoming
plateaus, whose heights and expanded interval endpoints converge.
Taking the finite maximum gives hypograph convergence of the
unpinned outputs, jointly for all sources. At interval endpoints,
upper-semicontinuous regularization selects the larger adjacent value.

Genericity of $\mathscr C_t$ makes the limiting output constant in
a neighborhood of zero. The matched finite description then also
gives convergence of its value at zero. Subtracting this value
therefore preserves convergence and proves the pinned assertion.

If $t=0$, the limiting future measure is empty. Since convergence
is on a compact envelope including time zero, integer-valued atom
matching eliminates all prelimit future atoms for large $\ell$.
The remaining evolution is a dilation of radius $t_\ell\to0$,
with a mesh restriction in the lattice case. The matched plateau
description proves hypograph convergence, and the initial atom-free
neighborhood of zero gives convergence after pinning.
\end{proof}

Joint marked convergence of the incoming lines, bridges, and the
independent future, followed by this lemma and
\eqref{eq:half-space-png-raw-grid-coupling-probability}, proves
\eqref{eq:half-space-png-joint-input-future-convergence} on $[0,L]$.
For local convergence, take a deterministic increasing exhaustion
with $K>L+t+1$ at each stage. Its endpoints are simultaneously
noncritical almost surely, so the claim holds in product
local-hypograph topology.

\paragraph*{\textbf{Stationarity and parameter consistency}}
Lemma~\ref{lem:half-space-png-equal-time-parameter-train}, applied
successively with fresh independent slabs, gives
\begin{equation*}
 (\Psi_\delta^{(c),[k]})^nH^\delta
 \overset{\mathrm d}{=}H^\delta,
 \qquad n\in\mathbb N_0.
\end{equation*}
Hence the two marginals on the left of
\eqref{eq:half-space-png-joint-input-future-convergence} have the
same law. Their limits must agree, proving
\begin{equation}
 \Theta_t^{(c),[k]}(R;\Pi)\overset{\mathrm d}{=}R,
 \qquad t\geq0,
 \label{eq:half-space-png-simultaneous-one-slab-identity}
\end{equation}
under \eqref{eq:half-space-png-strict-interior-parameters}.
The argument applies to each fixed deterministic $t$; convergence
of complete time trajectories is not needed for this invariance.

Finally, for $k\geq2$ and $1\leq m\leq k$, coordinate deletion is
continuous in the product topology
and commutes with $\mathcal I_\delta$. Applying
\eqref{eq:half-space-png-antidiagonal-line-law-convergence} to the
two sides of \eqref{eq:half-space-png-geometric-sequence-deletion}
gives
\begin{equation}
 \left(R_i^{(c;\boldsymbol g)}:i\neq m\right)
 \overset{\mathrm d}{=}
 \left(R_j^{(c;\boldsymbol g^{\widehat m})}:1\leq j\leq k-1\right),
 \label{eq:half-space-png-poisson-deletion-consistency}
\end{equation}
with the retained coordinates listed in their original order.
Moreover, for every input and every realization of the common future,
\begin{equation*}
 \mathsf D_m\Theta_t^{(c),[k]}(\boldsymbol h;\Pi)
 =\Theta_t^{(c),[k-1]}(\mathsf D_m\boldsymbol h;\Pi).
\end{equation*}
Thus deletion also preserves the retained common-noise evolutions.

\subsubsection{Endpoints and completion of the proof}
\label{subsubsec:half-space-png-deletion-and-endpoints}

We extend joint stationarity and deletion consistency from
$c>0$, $g_1>\cdots>g_k>c\vee1$ to all admissible parameters.
The argument couples the initial profiles and their common future
so that they agree locally with their limits with probability tending
to one. Proper bridge limits are treated by direct coupling; at
$g_k=c\geq1$, the diverging terminal bridge is removed by pinning
and the source shift of
Lemma~\ref{lem:half-space-png-forced-terminal-bridge}.

\begin{lemma}[Parameter limits with a common future]
\label{lem:half-space-png-parameter-limit-common-future}
Let $c_n>0$, $c_n\to c$, and
$\boldsymbol g^{(n)}\to\boldsymbol g$, where
$g_1^{(n)}>\cdots>g_k^{(n)}>c_n\vee1$,
$g_1>\cdots>g_k\geq r_c$.
Assume either
\begin{enumerate}[label=\textup{(\roman*)},leftmargin=*]
\item $g_k>c$; or
\item $g_k=c\geq1$, $c_n=c$, $g_i^{(n)}=g_i$ for $i<k$, and
      $g_k^{(n)}\downarrow c$.
\end{enumerate}
Set $\boldsymbol R^{(n)}=(R_i^{(c_n;\boldsymbol g^{(n)})})_{i=1}^k$
and $\boldsymbol R=(R_i^{(c;\boldsymbol g)})_{i=1}^k$.
In case \textup{(ii)}, define the last limiting profile using the
shifted source in \eqref{eq:half-space-png-singular-shifted-source}.
There is a coupling of these initial profiles and future environments
$\Pi_n,\Pi$, with boundary rates $c_n,c$, such that each future is
independent of its initial vector and shared by all its coordinates,
and, for each fixed deterministic $t\geq0$,
\begin{equation}
 \left(\boldsymbol R^{(n)},
       \Theta_t^{(c_n),[k]}(\boldsymbol R^{(n)};\Pi_n)\right)
 \xrightarrow[n\to\infty]{\mathbb P}
 \left(\boldsymbol R,
       \Theta_t^{(c),[k]}(\boldsymbol R;\Pi)\right)
 \label{eq:half-space-png-parameter-limit-with-common-future}
\end{equation}
in the product local-hypograph topology.
\end{lemma}

\begin{proof}
\emph{Regular limits.}
In case \textup{(i)}, every limiting bridge parameter is below one.
For $q_n\to q<1$, the masses $(1-q_n)q_n^m$ converge to
$(1-q)q^m$, giving
\begin{equation*}
 d_{\mathrm{TV}}\bigl(\operatorname{Geo}(q_n),
                      \operatorname{Geo}(q)\bigr)\longrightarrow0.
\end{equation*}
Independent maximal couplings therefore make the finite bridge
arrays agree with probability tending to one.
Couple each incoming Poisson clock by a common rate-cutoff measure,
using independent coupling randomness for distinct clocks and bridges:
for rates $a_n\to a$, its mismatch probability on $[0,K]$ is
$1-e^{-K|a_n-a|}\leq K|a_n-a|$.
Since all line rates converge, their configurations on $[0,K]$
also agree with probability tending to one.

\emph{The singular limit.}
In case \textup{(ii)}, apply the same coupling to the retained
bridges and line clocks. The remaining bridge
$Y_n=Z_n(-k,k+1)$ diverges in probability.
Lemma~\ref{lem:half-space-png-forced-terminal-bridge} shows that,
with probability tending to one, all maximizing paths from the last
source collect $Y_n$ throughout a fixed compact window.
Pinning cancels this common contribution, leaving the shifted-source
profile. The other sources cannot access $Y_n$.
Thus in both cases the entire pinned input vector agrees with its
limit on $[0,K]$ with probability tending to one.

\emph{The common future.}
Use identical bulk futures and couple the boundary clocks by common
rate cutoffs. Their mismatch probability up to time $t$ is at most
$t|c_n-c|$. Choose the randomness for this future coupling independently
of the entire incoming coupling. For $K>L+t+1$, finite propagation
then gives identical pinned outputs on $[0,L]$ whenever the inputs
agree on $[0,K]$ and the future configurations agree in the relevant
backward cone. Both events have probability tending to one.
This proves the joint convergence on every compact interval and
hence in the product local-hypograph topology.
\end{proof}

\paragraph*{\textbf{Approximating the endpoints}}
For a distinct ordered admissible vector, all cases outside the
strict interior are covered by the following choices, with
$\varepsilon_n\downarrow0$ and $0<\varepsilon_n<1$:
\begin{enumerate}[label=\textup{(\roman*)},leftmargin=*]
\item If $c=0$, take $c_n=\varepsilon_n$ and
      $g_i^{(n)}=g_i+\varepsilon_n$ for every $i$.
      This includes $g_k=1$ and gives proper limiting bridges.
\item If $0<c<1$ and $g_k=1$, fix $c_n=c$ and the first $k-1$
      parameters, and take $g_k^{(n)}=1+\varepsilon_n$.
      The last terminal bridge converges to $\operatorname{Geo}(c)$.
\item If $c\geq1$ and $g_k=c$, fix $c_n=c$ and the first $k-1$
      parameters, and take $g_k^{(n)}=c+\varepsilon_n$.
      The limit is the shifted-source construction; in particular,
      $c=1$, $g_k=1$ belongs to this singular case.
\end{enumerate}
For all sufficiently large $n$, these vectors remain strictly ordered
and satisfy \eqref{eq:half-space-png-strict-interior-parameters}.

\begin{proof}[Proof of Theorem~\ref{thm:intro-half-space-png-joint-stationary-family}]
\emph{Distinct ordered parameters.}
In the strict interior, joint stationarity and deletion consistency
are \eqref{eq:half-space-png-simultaneous-one-slab-identity} and
\eqref{eq:half-space-png-poisson-deletion-consistency}.
For an endpoint vector, choose an approximation above.
For each fixed deterministic $t$, the two components of the
approximating input--output pair have the same law.
Passing to their limits in
\eqref{eq:half-space-png-parameter-limit-with-common-future}
gives $\Theta_t^{(c),[k]}(\boldsymbol R;\Pi)
\overset{\mathrm d}{=}\boldsymbol R$. 
For $k\geq2$ and $1\leq m\leq k$, obtain the reduced approximating
vector by deleting the same coordinate from $\boldsymbol g^{(n)}$.
The parameter-limit lemma applies to both the full and reduced
vectors. At a singular endpoint, deleting $m<k$ retains the same
singular approximation in the smaller dimension, while deleting
$m=k$ leaves a fixed strictly interior vector.
Coordinate deletion is continuous, so the prelimit deletion identity
passes to
\begin{equation}
 \left(R_i^{(c;\boldsymbol g)}:i\neq m\right)
 \overset{\mathrm d}{=}
 \left(R_j^{(c;\boldsymbol g^{\widehat m})}:1\leq j\leq k-1\right),
 \label{eq:half-space-png-endpoint-deletion-consistency}
\end{equation}
with the retained coordinates in their original order.

\emph{Arbitrary parameter lists.}
Construct the joint law at the distinct values in decreasing order,
then permute and duplicate coordinates according to the requested
list. These operations commute with the common-noise evolution;
in particular, identical copies remain identical under the same
future. Joint stationarity follows.
For deletion from a list of length at least two, if the parameter
still occurs elsewhere, deletion removes an identical copy.
If its occurrence is unique, apply
\eqref{eq:half-space-png-endpoint-deletion-consistency} to the
distinct-parameter family. This proves deletion consistency for all
lists. Permutation and repetition consistency hold by construction.

\emph{One-profile marginals.}
Iterated deletion reduces each coordinate to its one-parameter law.
For $g>c$, Lemma~\ref{lem:half-space-png-one-parameter-check}
identifies this law as $\nu_{c,g}$.
At $g=c\geq1$, \eqref{eq:half-space-png-singular-k-one} gives the
same identification through the shifted source $(0,2)$.
This proves all the assertions of the theorem.
\end{proof}

\subsection{Half-space PNG horizon}
\label{subsec:half-space-png-horizon}

We now realize the jointly stationary laws
$\nu_{c;g_1,\ldots,g_k}$ as one measurable field indexed by
$g\in\mathcal G_c$. The field is determined by its profile at
$r_c=c\vee1$ and a positive measure recording its parameter increments.
Planar crossing gives positivity of the measure, while right continuity
in the parameter selects a canonical realization. We then prove
uniqueness in law of this regular representation and its stationarity,
as asserted in Theorem~\ref{thm:intro-half-space-png-horizon}.

\subsubsection{Right-continuous representatives and planar crossing}
\label{subsubsec:planar-crossing}

For $h\in\mathcal X_{\mathrm{step},0}$ without coincident bulk pairs or a boundary
mark at zero, define its right-continuous representative by
\begin{equation}
 h^{\to}(x):=\lim_{y\downarrow x,\,y>x}h(y),\qquad
 h=\operatorname{usc}h^{\to},\qquad
 (\operatorname{usc}f)(x)=
 \begin{cases}f(0),&x=0,\\\max\{f(x-),f(x)\},&x>0.\end{cases}
\label{eq:half-space-png-horizon-measurable-usc}
\end{equation}
Every finite-parameter law is supported on such profiles.  Step
coordinates (Definition~\ref{def:half-space-png-step-coordinate-space})
make $\mathfrak c(h):=h^{\to}$ a Borel map, extended by the zero
profile on the excluded set. For a field $F$, set
\begin{equation*}
 \Delta_\square F(g_1,g_2;x_1,x_2)
 :=F(g_1,x_1)+F(g_2,x_2)-F(g_1,x_2)-F(g_2,x_1).
\end{equation*}

\begin{lemma}[Planar crossing]
\label{lem:half-space-png-horizon-monge}
For $r_c\leq a<b$ and $(R_a,R_b)\sim\nu_{c;a,b}$, almost surely
\begin{equation*}
 R_a^{\to}(x_1)+R_b^{\to}(x_2)
 -R_a^{\to}(x_2)-R_b^{\to}(x_1)\geq0,
 \qquad 0\leq x_1<x_2.
\end{equation*}
Thus $R_b^{\to}-R_a^{\to}$ is nonnegative and nondecreasing, and
$R_a\leq R_b$.
\end{lemma}

\begin{proof}
First take $a>c$.  The descending construction $(b,a)$ has sources
$s_1=(-1,1)$, $s_2=(-2,2)$, terminal line $4$, and unrecentered last-passage profiles
$G_b,G_a$.  Maximizers $P_b:s_1\to(x_1,4)$ and
$P_a:s_2\to(x_2,4)$ must have intersecting closed horizontal intervals
on a common line, including switching endpoints: otherwise $P_a$ lies strictly left of
$P_b$ on line $2$, and this order propagates through their common
lines, forcing $x_2<x_1$. Switch their tails on the first line where they
intersect. The overlap makes both crossed intervals admissible.
On the meeting line, with entrances $\alpha_i$ and exits
$\beta_i$,
\begin{equation*}
 \mathrm dE((\alpha_1,\beta_1])+
 \mathrm dE((\alpha_2,\beta_2])
 =\mathrm dE((\alpha_1,\beta_2])+
 \mathrm dE((\alpha_2,\beta_1]).
\end{equation*}
Entrance singleton contributions are unchanged and subsequent segments
are exchanged, preserving signed weight. Thus
$G_b(x_1)+G_a(x_2)\leq G_b(x_2)+G_a(x_1)$, and pinning leaves this
cross increment unchanged.

For $a=c\geq1$, take $a_n\downarrow c$ with $c<a_n<b$.
Lemma~\ref{lem:half-space-png-forced-terminal-bridge}, with the common
clock-and-bridge coupling on a larger compact interval, gives agreement
of both pinned right-continuous representatives for $(b, a_n)$ and $(b, c)$ with probability
tending to one.  Pass the inequality at rational spatial points and
then everywhere by right-continuity.  The difference is nondecreasing and vanishes at zero; applying
$\operatorname{usc}$ preserves the order.
\end{proof}

\begin{lemma}[One-parameter right continuity]
\label{lem:half-space-png-horizon-parameter-continuity}
If $g_n\downarrow g$ in $\mathcal G_c$, there is a coupling such that,
for every $L<\infty$,
\begin{equation*}
 \mathbb P\left(
 H^{(c,g_n),\to}|_{[0,L]}=H^{(c,g),\to}|_{[0,L]}
 \right)\longrightarrow1.
\end{equation*}
Consequently $\nu_{c,g_n}\Rightarrow\nu_{c,g}$ locally in hypograph
topology.
\end{lemma}

\begin{proof}
An eventually constant tail is coupled identically; otherwise $g_n>g$.
Let $M^{(g)}$ be the running maximum in
\eqref{eq:intro-half-space-png-stationary-profile-raw}.
For $g>c$, common Poisson random measures with rate cutoffs couple the
four counters with mismatch probability at most
$2L(|g_n-g|+|g_n^{-1}-g^{-1}|)$.
Maximal coupling makes the geometric variables agree with probability
tending to one, since $c/g_n\to c/g<1$.
At $g=c\geq1$, $Y^{(c,g_n)}\to\infty$ in probability, whereas
$M^{(g_n)}(L)$ is tight, bounded by four counters with bounded rates.
The regulator therefore vanishes on $[0,L]$ with probability tending
to one, leaving the coupled $S_2^{(c)}=H^{(c,c),\to}$.
Applying $\operatorname{usc}$ gives the hypograph convergence.
\end{proof}

\subsubsection{Canonical regularization and stationarity}
\label{subsubsec:canonical-regularization}

Put $r=r_c$.  Let $\mathcal D_{+,0}$ be the standard Borel space of
pinned, locally finite, integer-valued c\`adl\`ag step functions, with
its jump-coordinate Borel structure.  Let $\mathcal E_c$ consist of
$(\mathfrak b,\mu)$, where the base profile $\mathfrak b\in\mathcal D_{+,0}$ and $\mu$ is a positive
locally finite Borel measure on $\mathcal G_c\times\mathbb R_+$,
zero on the lower faces $\{r\}\times\mathbb R_+$ and
$\mathcal G_c\times\{0\}$.  The vague measure structure and Borel lower-face conditions make
$\mathcal E_c$ standard Borel. Set
\begin{equation*}
 E_{\mathfrak b,\mu}(g,x):=\mathfrak b(x)+\mu((r,g]\times(0,x]),\qquad
 \mathfrak R_c(\mathfrak b,\mu)(g,\cdot):=\operatorname{usc}E_{\mathfrak b,\mu}(g,\cdot).
\end{equation*}
Local finiteness makes each $E_{\mathfrak b,\mu}(g,\cdot)$ c\`adl\`ag.
Rectangle evaluations, spatial left limits, and
\eqref{eq:half-space-png-horizon-measurable-usc} make the reconstruction
jointly Borel in $(\mathfrak b,\mu,g,x)$.

\begin{definition}[Regular half-space PNG horizon]
\label{def:half-space-png-horizon}
A regular half-space PNG horizon of boundary rate $c$ is a random
$(\mathfrak b_c^{\mathrm{PNG}},\mu_c^{\mathrm{PNG}})\in\mathcal E_c$ whose
reconstructed field $H_c^{\mathrm{PNG}}:=
\mathfrak R_c(\mathfrak b_c^{\mathrm{PNG}},\mu_c^{\mathrm{PNG}})$ belongs to
$\mathcal X_{\mathrm{step},0}$ at every $g\in\mathcal G_c$ and satisfies
\begin{equation}
 \operatorname{Law}(H_c^{\mathrm{PNG}}(g_1,\cdot),\ldots,
                    H_c^{\mathrm{PNG}}(g_k,\cdot))
       =\nu_{c;g_1,\ldots,g_k}
\label{eq:half-space-png-horizon-finite-restrictions}
\end{equation}
for every deterministic finite list.  Write
$H_c^{\mathrm{PNG},\to}:=E_{\mathfrak b_c^{\mathrm{PNG}},\mu_c^{\mathrm{PNG}}}$.
\end{definition}

The representing measure satisfies
\begin{equation}
 \mu_c^{\mathrm{PNG}}((g_1,g_2]\times(x_1,x_2])
 =\Delta_\square H_c^{\mathrm{PNG},\to}(g_1,g_2;x_1,x_2).
\label{eq:half-space-png-horizon-cross-increment-measure}
\end{equation}
For $g_n\downarrow g$, continuity from above gives
\begin{equation*}
 \sup_{0\leq x\leq L}|H_c^{\mathrm{PNG},\to}(g_n,x)
                        -H_c^{\mathrm{PNG},\to}(g,x)|
 =\mu_c^{\mathrm{PNG}}((g,g_n]\times(0,L])\longrightarrow0.
\end{equation*}

\begin{proof}[Proof of Theorem~\ref{thm:intro-half-space-png-horizon}]
\emph{Construction.}
Let $\mathbb Q_c=\{r\}\cup(\mathbb Q\cap(r,\infty))$.
By projective consistency, the standard Borel step space, and
Kolmogorov's extension theorem, there is a family
$(H_c^0(q))_{q\in\mathbb Q_c}$ with the prescribed finite-dimensional laws.  Put $f_q=\mathfrak c(H_c^0(q))$, $\mathfrak b=f_r$, and
$F_q^{\mathrm{inc}}=f_q-\mathfrak b$.  On one full-probability event these are the right-continuous
representatives and $F_{q_2}^{\mathrm{inc}}-F_{q_1}^{\mathrm{inc}}$ is nonnegative and nondecreasing.
For $q\in\mathbb Q_c$ and rational $x\geq0$, choose
$p_n\in\mathbb Q_c$ decreasing to $q$ from above.  Then
$f_q(x)\leq f_{p_n}(x)\downarrow Z$.
Lemma~\ref{lem:half-space-png-horizon-parameter-continuity} gives
$Z\stackrel{\mathrm d}{=}f_q(x)$, so order forces equality.
A countable intersection yields
$F_q^{\mathrm{inc}}(x)=\inf_{p\in\mathbb Q_c,\,p>q}F_p^{\mathrm{inc}}(x)$ for
$q\in\mathbb Q_c$ and $x\in\mathbb Q_{\geq0}$. Define
\begin{equation*}
 F(g,x):=\inf_{\substack{q\in\mathbb Q_c\\q>g}}
         \inf_{\substack{y\in\mathbb Q_{\geq0}\\y>x}}F_q^{\mathrm{inc}}(y).
\end{equation*}
Commuting the infima and using spatial right-continuity gives
$F(g,x)=\inf_{q\in\mathbb Q_c,\,q>g}F_q^{\mathrm{inc}}(x)$ and
$F(q,x)=F_q^{\mathrm{inc}}(x)$ for $q\in\mathbb Q_c$.
For the second identity, also use the preceding right-continuity identity at
rational $y\downarrow x$.  The function $F$ is finite-valued, integer-valued, nondecreasing and
right-continuous in each variable;
agreement on $\mathbb Q_c$ gives $F(r,x)=F(g,0)=0$.
Taking synchronized $q_{i,n}\downarrow g_i$ with $q_{1,n}<q_{2,n}$
in the crossing inequality gives $\Delta_\square F\geq0$.

For fixed $G>g$ and sufficiently small $u,v>0$, the rectangle inequality gives
\begin{equation*}
 0\leq F(g+u,x+v)-F(g,x)
 \leq F(G,x+v)-F(G,x)+F(g+u,x)-F(g,x).
\end{equation*}
Separate right-continuity therefore gives joint continuity from above.
The nonnegative rectangle increments define an additive content with
its Carath\'eodory extension, a unique positive
locally finite measure $\mu$, zero on the lower faces, with
\begin{equation*}
 \mu((g_1,g_2]\times(x_1,x_2])
 =\Delta_\square F(g_1,g_2;x_1,x_2),\qquad
 F(g,x)=\mu((r,g]\times(0,x]).
\end{equation*}
For finite $G>r$ and $L>0$, local finiteness follows from $F(G,L)<\infty$.
Masses of rectangles with endpoints in the chosen countable dense sets
are measurable and generate the vague Borel structure. The preceding
construction defines the regularization map
\begin{equation*}
 \mathfrak K_c:(\mathcal X_{\mathrm{step},0})^{\mathbb Q_c}\longrightarrow\mathcal E_c,
 \qquad (H(q))_{q\in\mathbb Q_c}\longmapsto(\mathfrak b,\mu).
\end{equation*}
The good event is specified by regularity of each skeleton profile,
the crossing inequalities and parameter right-continuity identities,
both tested at rational spatial points. These are countably many Borel conditions.
Thus $\mathfrak K_c$ is Borel, taking the value $(0,0)$ off this event.

Set $(\mathfrak b_c^{\mathrm{PNG}},\mu_c^{\mathrm{PNG}})=(\mathfrak b,\mu)$ and reconstruct
by $\mathfrak R_c$.  For every $g,L$ and $q_n\downarrow g$ in
$\mathbb Q_c$, one has
$\mu((g,q_n]\times(0,L])=F_{q_n}^{\mathrm{inc}}(L)-F(g,L)\downarrow0$.
These masses are nonnegative integers and eventually vanish. On this
single skeleton event the argument applies to every $g\in\mathcal G_c$ and finite $L\geq0$.
Thus every slice agrees on each compact interval with a skeleton profile,
proving admissibility at all parameters and eventual compact equality
from the right.

\emph{Identification and uniqueness.}
For a deterministic list $g_1,\ldots,g_k$, extend the skeleton to
$\mathbb Q_c\cup\{g_1,\ldots,g_k\}$.
For $q_n\downarrow g_i$, planar crossing gives
$\mathfrak c(H_c^0(g_i))(x)\leq \mathfrak b(x)+F(g_i,x)=\lim_nf_{q_n}(x)$.
Marginal continuity identifies their laws, so order forces equality at
rational $x$, then everywhere by right-continuity:
\begin{equation}
 \mathfrak c(H_c^0(g_i))(x)=\mathfrak b(x)+F(g_i,x)
       =H_c^{\mathrm{PNG},\to}(g_i,x),\qquad x\geq0.
\label{eq:half-space-png-horizon-regularization-identification}
\end{equation}
For fixed $g_i$, this identification event is Borel: equality of the
right-continuous profiles is determined by their rational spatial values.
Intersect over the finite list and apply $\operatorname{usc}$ to obtain
\eqref{eq:half-space-png-horizon-finite-restrictions}; rectangle
subtraction gives
\eqref{eq:half-space-png-horizon-cross-increment-measure}.
Each $E_{\mathfrak b,\mu}(g,\cdot)$ is c\`adl\`ag, so
$(\operatorname{usc}E_{\mathfrak b,\mu}(g,\cdot))^{\to}=E_{\mathfrak b,\mu}(g,\cdot)$.
Thus the skeleton profiles recover $\mathfrak b=f_r$ and all rectangle
masses of $\mu$ on the chosen dense sets, proving uniqueness of the law of $(\mathfrak b,\mu)$.

\emph{Stationarity.}
Fix $t\geq0$ and let $H_{c,t}^{\mathrm{ev}}$ be the evolved field in
Theorem~\ref{thm:intro-half-space-png-horizon}.
Graphical evolution followed by $\operatorname{Rec}$ is Borel on the
step space, and joint stationarity gives $(H_{c,t}^{\mathrm{ev}}(q))_q$ the original
skeleton's countable-product law. Put
$(\widehat{\mathfrak b}_t,\widehat\mu_t)=\mathfrak K_c((H_{c,t}^{\mathrm{ev}}(q))_q)$ and
$\widehat H_{c,t}^{\mathrm{PNG}}=\mathfrak R_c(\widehat{\mathfrak b}_t,\widehat\mu_t)$.
Borel reconstruction preserves the law of $(\mathfrak b,\mu)$. For each fixed deterministic $g$, joint stationarity gives equality in
law of the augmented countable families
$((H_{c,t}^{\mathrm{ev}}(q))_q,H_{c,t}^{\mathrm{ev}}(g))$ and
$((H_c^{\mathrm{PNG}}(q))_q,H_c^{\mathrm{PNG}}(g))$.
Thus the measurable identification event
\eqref{eq:half-space-png-horizon-regularization-identification} transfers
to the evolved field.  Hence $\widehat H_{c,t}^{\mathrm{PNG}}(g)=H_{c,t}^{\mathrm{ev}}(g)$
almost surely, proving
\eqref{eq:intro-half-space-png-horizon-stationarity} for each fixed time.
\end{proof}

\section{KPZ universality}\label{sec:kpz-universality}

Using the scalings of Section~\ref{subsec:intro-kpz-universality},
we first prove convergence of the PNG horizon to the half-space
stationary horizon in Section~\ref{subsec:stationary-horizon-convergence}.
Section~\ref{subsec:kpz-fixed-point-convergence} establishes convergence
of the height evolution to the half-space KPZ fixed point through
transition-probability asymptotics, stationary comparison, and tightness.
Finally, Section~\ref{subsec:directed-landscape-convergence} proves
convergence of the passage field to the half-space directed landscape,
together with joint convergence of evolutions from several initial
profiles driven by the same environment.

\subsection{Convergence to the half-space stationary horizon}
\label{subsec:stationary-horizon-convergence}

We prove convergence of the stationary PNG horizon in two senses:
joint convergence at every deterministic finite parameter list, and
convergence of its base profile and representing measure in
$\mathscr D(\lambda_\rho)$. For the first, we pass from the shared
Poisson-line construction to the Brownian-line construction of the
half-space stationary horizon. For the second, positivity of rectangle
increments gives tightness of the measures, and thin-strip estimates
exclude mass on the two lower faces. Both conclusions are retained,
since evaluation at a fixed parameter need not be continuous in the
measure topology.

Fix $\rho\in\mathbb R\cup\{-\infty\}$ and use
$c_\varepsilon=c_\varepsilon(\rho)$, $\lambda_\rho$, and
$g_\varepsilon$ from
\eqref{eq:intro-png-critical-boundary-scaling}--%
\eqref{eq:intro-png-slope-scaling}. Write
\begin{equation*}
 r_\varepsilon:=c_\varepsilon\vee1=g_\varepsilon(\lambda_\rho),
 \qquad \ell_\varepsilon(g):=2\varepsilon^{-1/2}\log g,
 \qquad \ell_\varepsilon^{-1}=g_\varepsilon.
\end{equation*}

\subsubsection{The rescaled PNG horizon}
\label{subsubsec:rescaled-png-horizon}

We first identify the representing pair of
$\mathcal H_\varepsilon^{\mathrm{PNG},\rho}$ and verify its
measurability. Spatial interpolation spreads each microscopic spatial
atom uniformly across its mesh interval. The atom's rescaled mass is
$\varepsilon^{1/2}$ times its original mass, so the corresponding
density on an interval of length $\varepsilon$ has factor
$\varepsilon^{-1/2}$.

For $y>0$, let $j(y)=n$ on $(n,n+1]$, $n\in\mathbb N_0$, and put
$j(0)=0$. For a Borel set
$A\subset[\lambda_\rho,\infty)\times\mathbb R_+$, define
\begin{equation*}
 \mu_\varepsilon(A):=\varepsilon^{-1/2}
 \int\!\int_{\varepsilon j(y)}^{\varepsilon(j(y)+1)}
 \mathbf1_A(\ell_\varepsilon(g),x)\,\mathrm dx\,
 \mu_{c_\varepsilon}^{\mathrm{PNG}}(\mathrm dg,\mathrm dy).
\end{equation*}
For $\lambda_\rho\leq\lambda_1<\lambda_2$, put
$g_a=g_\varepsilon(\lambda_a)$ and
\begin{equation*}
 F_n^{\mathrm{inc}}:=H_{c_\varepsilon}^{\mathrm{PNG},\to}(g_2,n)
                   -H_{c_\varepsilon}^{\mathrm{PNG},\to}(g_1,n).
\end{equation*}
At $x=\varepsilon(n+\theta)$, $0\leq\theta<1$, the microscopic
rectangle identity gives
\begin{align*}
 \mu_\varepsilon((\lambda_1,\lambda_2]\times(0,x])
 =\varepsilon^{1/2}
   [(1-\theta)F_n^{\mathrm{inc}}+\theta F_{n+1}^{\mathrm{inc}}]
 =\mathcal H_\varepsilon^{\mathrm{PNG},\rho}(\lambda_2,x)
   -\mathcal H_\varepsilon^{\mathrm{PNG},\rho}(\lambda_1,x).
\end{align*}
Thus, with
$\mathfrak b_\varepsilon(x):=
\mathcal H_\varepsilon^{\mathrm{PNG},\rho}(\lambda_\rho,x)$,
\begin{equation*}
 \mathcal H_\varepsilon^{\mathrm{PNG},\rho}(\lambda,x)
 =\mathfrak b_\varepsilon(x)
  +\mu_\varepsilon((\lambda_\rho,\lambda]\times(0,x]).
\end{equation*}
The measure is positive and locally finite. Its parameter lower face
has zero mass by the microscopic normalization, and integration with
respect to $\mathrm dx$ gives zero mass to the spatial lower face.
Interpolation makes every spatial profile continuous and pinned.
Mesh evaluations make $\mathfrak b_\varepsilon$ measurable, and Borel
integration makes $\int f\,\mathrm d\mu_\varepsilon$ measurable for
every $f\in C_c([\lambda_\rho,\infty)\times\mathbb R_+)$.
A countable determining family therefore makes
$(\mathfrak b_\varepsilon,\mu_\varepsilon)$ a Borel
$\mathscr D(\lambda_\rho)$-valued random element.

\subsubsection{Finite-parameter convergence}
\label{subsubsec:stationary-finite-parameters}

Fix $\lambda_1>\cdots>\lambda_k\geq\lambda_\rho$ and put
$\bar j=2k+1-j$. The limiting shared environment consists of independent
Brownian lines and exponential bridges. Let $B_1,\ldots,B_{2k}$ be
independent Brownian motions started at zero, with drift $d_j$ and
variance specified by
\begin{equation*}
 d_j:=\begin{cases}
        \lambda_j,&j\leq k,\\
        -\lambda_{\bar j},&j>k,
       \end{cases}
 \qquad
 \operatorname{Var}(B_j(y)-B_j(x))=2(y-x),\quad 0\leq x<y.
\end{equation*}
For $1\leq i\leq k$ and $i<j\leq\bar i$, define
\begin{equation*}
 a_{ij}^\rho:=
 \begin{cases}
  (\lambda_i-\lambda_j)/2,&i<j\leq k,\\
  (\lambda_i+\lambda_{\bar j})/2,&k<j<\bar i,\\
  \lambda_i/2-\rho,&j=\bar i,\ \rho\in\mathbb R,\\
  +\infty,&j=\bar i,\ \rho=-\infty.
 \end{cases}
\end{equation*}
For each positive rate take
$X^\rho(-i,j)\sim\operatorname{Exp}(a_{ij}^\rho)$, independently of
one another and of the Brownian lines; here
$\operatorname{Exp}(+\infty)$ is the point mass at zero.
Set $s_i=(-i,i)$. The only zero rate occurs when
$\rho\geq0$, $\lambda_k=2\rho$,
$i=k$, and $j=k+1$. In this case omit that atom and use the shifted
source $s_k=(-k+1,k+1)$.

As in the Poisson-line construction, concatenate the cumulative bridge
functions on $(-\infty,0]$, normalized at zero, with the Brownian lines
to obtain $E^\rho$. Using
Definition~\ref{def:half-space-png-cadlag-lpp} with marked columns
$\{-k,\ldots,-1\}$, set
\begin{equation*}
 R_i^\rho(x):=\mathcal L_{E^\rho}(s_i;(x,2k))
              -\mathcal L_{E^\rho}(s_i;(0,2k)).
\end{equation*}
The bridge-entry convention agrees with
\cite[Section~2.5]{dauvergne2026directed}: the positive-axis lines
are continuous and all atoms lie at the marked bridge columns.
The terminal line is included, consistently with the construction and
shifted-source formula in \cite[Proposition~2.18]{dauvergne2026directed}.
By \cite[Proposition~2.18 and Theorem~2.19]{dauvergne2026directed},
\begin{equation}
 (R_i^\rho)_{i=1}^k\overset{\mathrm d}{=}
 (\mathcal H^\rho(\lambda_i,\cdot))_{i=1}^k|_{\mathbb R_+}.
\label{eq:png-horizon-dz-finite-line-identification}
\end{equation}

\begin{proposition}[Finite-parameter invariance principle]
\label{prop:png-horizon-finite-slope-limit}
For every deterministic list
$\lambda_1>\cdots>\lambda_k\geq\lambda_\rho$,
\begin{equation*}
 (\mathcal H_\varepsilon^{\mathrm{PNG},\rho}
       (\lambda_i,\cdot))_{i=1}^k
 \Longrightarrow(\mathcal H^\rho(\lambda_i,\cdot))_{i=1}^k
 \quad\text{in }C_{\mathrm{loc}}(\mathbb R_+)^k.
\end{equation*}
\end{proposition}

\begin{proof}
Set $g_{i,\varepsilon}=g_\varepsilon(\lambda_i)$,
$\boldsymbol g_\varepsilon=(g_{1,\varepsilon},\ldots,g_{k,\varepsilon})$,
and $\gamma_{j,\varepsilon}=e^{d_j\varepsilon^{1/2}/2}$.

\emph{Step 1: convergence of the Poisson lines.}
Let $S_{j,\varepsilon}:=S_j^{\boldsymbol g_\varepsilon}$ and put
\begin{equation*}
 \widehat S_{j,\varepsilon}(x)
 :=\varepsilon^{1/2}S_{j,\varepsilon}(\varepsilon^{-1}x),\qquad
 \widehat M_{j,\varepsilon}(x)
 :=\widehat S_{j,\varepsilon}(x)
   -\varepsilon^{-1/2}(\gamma_{j,\varepsilon}
                      -\gamma_{j,\varepsilon}^{-1})x.
\end{equation*}
The martingale jumps have size $\varepsilon^{1/2}$, and locally uniformly,
\begin{equation*}
 \varepsilon^{-1/2}(\gamma_{j,\varepsilon}
                  -\gamma_{j,\varepsilon}^{-1})x\longrightarrow d_jx,
 \qquad
 \langle\widehat M_{j,\varepsilon}\rangle(x)
 = (\gamma_{j,\varepsilon}+\gamma_{j,\varepsilon}^{-1})x
 \longrightarrow2x.
\end{equation*}
The vanishing jumps give the Lindeberg condition, and independence
gives zero cross variations. The martingale functional central limit
theorem therefore yields joint convergence of the lines to
$(B_j)_{j=1}^{2k}$ in $D_{\mathrm{loc}}(\mathbb R_+)^{2k}$ with the
local Skorokhod $J_1$ topology.

\emph{Step 2: convergence of the bridges and treatment of endpoints.}
For each finite positive $a_{ij}^\rho$, the geometric parameter in
\eqref{eq:half-space-png-bridge-laws-spectral} is
$e^{-a_{ij}^\rho\varepsilon^{1/2}}$. If
$Z_\varepsilon\sim\operatorname{Geo}(e^{-a\varepsilon^{1/2}})$ with
$a>0$, then, for $x\geq0$,
\begin{equation*}
 \mathbb P(\varepsilon^{1/2}Z_\varepsilon>x)
 =e^{-a\varepsilon^{1/2}
       (\lfloor x\varepsilon^{-1/2}\rfloor+1)}
 \longrightarrow e^{-ax}.
\end{equation*}
For $\rho=-\infty$, terminal bridge parameters are at most
$c_\star<1$, so
$\mathbb E[\varepsilon^{1/2}Z]\leq
\varepsilon^{1/2}c_\star/(1-c_\star)\to0$.
If $a_{k,k+1}^\rho=0$, then $g_{k,\varepsilon}=c_\varepsilon\geq1$
for every $\varepsilon$: both constructions omit the terminal bridge
and use the same shifted source. This includes $\rho=0$, $\lambda_k=0$.
Independence makes convergence of all the lines and retained bridges
joint, with the specified independent limiting inputs.

\emph{Step 3: convergence of the passage profiles.}
Let $E_\varepsilon=(E_{j,\varepsilon})_{j=1}^{2k}$ be the corresponding
shared PNG environment, using the shifted-source convention when
needed. The negative coordinates label a finite bridge array, so we
keep them fixed and rescale only the positive axis. Set
\begin{equation*}
 \phi_\varepsilon(u):=
 \begin{cases}u,&u\leq0,\\\varepsilon u,&u\geq0,\end{cases}
 \qquad
 \widehat E_{j,\varepsilon}(x)
 :=\varepsilon^{1/2}E_{j,\varepsilon}(\phi_\varepsilon^{-1}(x)).
\end{equation*}
This increasing coordinate change preserves switching paths and the
bridge-entry convention. Hence
\begin{equation*}
 \mathcal L_{\widehat E_\varepsilon}(s_i;(x,2k))
 -\mathcal L_{\widehat E_\varepsilon}(s_i;(0,2k))
 =\varepsilon^{1/2}
  \widetilde R_i^{(c_\varepsilon;\boldsymbol g_\varepsilon)}
                    (\varepsilon^{-1}x).
\end{equation*}
Along any sequence $\varepsilon\downarrow0$, use a joint Skorokhod
representation of the lines and bridges. Continuity of the Brownian
limits gives uniform convergence of the positive lines on compacts.
Together with convergence of the bridge masses, this gives, almost surely,
\begin{equation*}
 \eta_{\varepsilon,L}:=\max_{1\leq j\leq2k}
 \sup_{-k-1\leq x\leq L}
 |\widehat E_{j,\varepsilon}(x)-E_j^\rho(x)|\longrightarrow0,
\end{equation*}
also controlling left limits at the fixed bridge columns.
Every path uses at most $2k$ lines, and each line contributes at most
$4\eta_{\varepsilon,L}$ to the weight error. Taking suprema and then
pinning therefore gives
\begin{equation*}
 \sup_{0\leq x\leq L}
 \left|\varepsilon^{1/2}
 \widetilde R_i^{(c_\varepsilon;\boldsymbol g_\varepsilon)}
           (\varepsilon^{-1}x)-R_i^\rho(x)\right|
 \leq16k\eta_{\varepsilon,L}.
\end{equation*}

\emph{Step 4: interpolation and identification of the limiting law.}
By \eqref{eq:half-space-png-horizon-finite-restrictions} and
\eqref{eq:half-space-png-horizon-regularization-identification}, the
canonical right-continuous horizon representatives have these pinned
profiles' finite joint laws. The Brownian passage profiles $R_i^\rho$
are continuous, by uniform continuity of the finitely many Brownian
lines on compacts.
Write $V_{i,\varepsilon}$ for the scaled c\`adl\`ag profile and
$\overline V_{i,\varepsilon}$ for its mesh interpolation. For a
continuous $f$, write
\begin{equation*}
 \omega_f(a;L):=
 \sup_{\substack{x,y\in[0,L]\\|x-y|\leq a}}|f(x)-f(y)|.
\end{equation*}
For $\varepsilon\leq1$,
\begin{equation*}
 \|\overline V_{i,\varepsilon}-R_i^\rho\|_{[0,L]}
 \leq\|V_{i,\varepsilon}-R_i^\rho\|_{[0,L+\varepsilon]}
     +\omega_{R_i^\rho}(\varepsilon;L+1)
 \longrightarrow0.
\end{equation*}
The joint convergence and
\eqref{eq:png-horizon-dz-finite-line-identification} prove the proposition.
\end{proof}

\subsubsection{Convergence of the representing pair}
\label{subsubsec:stationary-process-convergence}

The preceding proposition identifies every prescribed finite collection
of profiles. We now use positivity to obtain tightness and identify
the law of the complete representing pair.

\begin{lemma}[Convergence in $\mathscr D(\lambda_0)$]
\label{lem:png-horizon-positive-measure-upgrade}
Let $G_n,G$ be Borel $\mathscr D(\lambda_0)$-valued random elements,
where $\mathscr D(\lambda_0)$ is the space defined in
Section~\ref{subsubsec:intro-stationary-horizon-convergence} with
lower parameter endpoint $\lambda_0$.
Assume that, for every deterministic $\lambda\geq\lambda_0$, their
spatial slices are $C_{\mathrm{loc}}(\mathbb R_+)$-valued random
elements. If, for every deterministic finite parameter list,
\begin{equation*}
 (G_n(\lambda_i,\cdot))_{i=1}^m
 \Longrightarrow(G(\lambda_i,\cdot))_{i=1}^m
 \quad\text{in }C_{\mathrm{loc}}(\mathbb R_+)^m,
\end{equation*}
then $G_n\Longrightarrow G$ in $\mathscr D(\lambda_0)$.
\end{lemma}

\begin{proof}
Write $(\mathfrak b_n,\mu_n)$ and $(\mathfrak b,\mu)$ for the
representing pairs. Work first in the ambient Polish space
\begin{equation*}
 C_{\mathrm{loc},0}(\mathbb R_+)
 \times\mathcal M_{\mathrm{loc},+}
       ([\lambda_0,\infty)\times\mathbb R_+),
\end{equation*}
where the first factor consists of continuous profiles pinned at zero
with the local uniform topology, and the second of positive locally
finite measures with the vague topology.

\emph{Step 1: tightness.}
The hypothesis at $\lambda_0$ makes $(\mathfrak b_n)$ tight. For
finite $\Lambda>\lambda_0$ and $L>0$,
\begin{equation*}
 \mu_n([\lambda_0,\Lambda]\times[0,L])
 =G_n(\Lambda,L)-G_n(\lambda_0,L)
\end{equation*}
is tight by the two-parameter hypothesis. Positivity and a compact
exhaustion give vague tightness of $(\mu_n)$, and hence tightness of
the pairs in the ambient space.

\emph{Step 2: exclusion of lower-face mass.}
Zero lower-face masses need not be preserved by vague limits, so we
control the mass of thin strips. Put
\begin{equation*}
 X_{n,\delta}^L:=G_n(\lambda_0+\delta,L)-G_n(\lambda_0,L),
 \qquad
 Y_{n,\delta}^{\Lambda}:=G_n(\Lambda,\delta)-G_n(\lambda_0,\delta).
\end{equation*}
These are the masses of
$(\lambda_0,\lambda_0+\delta]\times(0,L]$ and
$(\lambda_0,\Lambda]\times(0,\delta]$, respectively.
For fixed $\delta$ they converge in law to their target counterparts
$X_\delta^L,Y_\delta^{\Lambda}$, which decrease to zero as
$\delta\downarrow0$ by continuity from above. Portmanteau therefore
gives, for every $\eta>0$,
\begin{equation*}
 \lim_{\delta\downarrow0}\limsup_n
 \mathbb P(X_{n,\delta}^L>\eta)
 =\lim_{\delta\downarrow0}\limsup_n
 \mathbb P(Y_{n,\delta}^{\Lambda}>\eta)=0.
\end{equation*}

Take an ambient subsequential limit $(\mathfrak b_*,\mu_*)$ and a
Skorokhod coupling. A continuous compactly supported cutoff in
$[0,1]$, equal to one on $\{\lambda_0\}\times[0,L]$ and supported in
$[\lambda_0,\lambda_0+\delta)\times[0,L+1)$, has $\mu_n$ integral
at most $X_{n,\delta}^{L+1}$. Similarly, a cutoff equal to one on
$[\lambda_0,\Lambda]\times\{0\}$ and supported in
$[\lambda_0,\Lambda+1)\times[0,\delta)$ has integral at most
$Y_{n,\delta}^{\Lambda+1}$. Almost sure vague convergence yields
\begin{align*}
 \mathbb P\bigl(\mu_*(\{\lambda_0\}\times[0,L])>\eta\bigr)
 &\leq\limsup_n
       \mathbb P(X_{n,\delta}^{L+1}>\eta/2),\\
 \mathbb P\bigl(\mu_*([\lambda_0,\Lambda]\times\{0\})>\eta\bigr)
 &\leq\limsup_n
       \mathbb P(Y_{n,\delta}^{\Lambda+1}>\eta/2).
\end{align*}
Let $\delta\downarrow0$ and exhaust the two faces. Both lower faces
are $\mu_*$-null almost surely.

\emph{Step 3: identification of the subsequential limit.}
Choose deterministic countable dense grids
$Q\subset(\lambda_0,\infty)$ and $T\subset(0,\infty)$ whose vertical
and horizontal lines are $\mu_*$-null almost surely.
Such grids exist: the coordinate projections of the finite measures
\begin{equation*}
 \overline\mu_K(A):=
 \mathbb E\left[\frac{\mu_*(A\cap K)}{1+\mu_*(K)}\right]
\end{equation*}
have only countably many atoms, and one avoids these atoms over a
compact exhaustion $K$.
The resulting rectangles are continuity sets, so on the coupling,
\begin{equation*}
 G_n(q,t)
 =\mathfrak b_n(t)+\mu_n((\lambda_0,q]\times(0,t])
 \longrightarrow
 \mathfrak b_*(t)+\mu_*((\lambda_0,q]\times(0,t])=:G_*(q,t).
\end{equation*}
The finite-parameter hypothesis identifies every finite grid law,
jointly with the base-profile values, with the corresponding law of
$G$. Dense base values determine the continuous base profile, and
rectangle increments on the grids determine the locally finite measure.
Together with the zero lower-face masses, this proves
$(\mathfrak b_*,\mu_*)\overset{\mathrm d}{=}(\mathfrak b,\mu)$.
Every ambient subsequential limit therefore lies in
$\mathscr D(\lambda_0)$ and has the target pair law. This proves the
claimed convergence in its induced topology.
\end{proof}

\begin{proof}[Proof of Theorem~\ref{thm:intro-half-space-png-stationary-horizon-limit}]
Restrict the horizon of \cite[Theorem~2.19]{dauvergne2026directed}
to $x\geq0$. Restrict its representing measure to
$(\lambda_\rho,\infty)\times(0,\infty)$ and extend by zero on the
lower faces; this leaves all reconstruction rectangles unchanged.
Positivity of cross increments from \cite[Definition~2.13]{dauvergne2026directed} gives
\begin{equation*}
 \mu^\rho([\lambda_\rho,\Lambda]\times[0,L])
 =\mathcal H^\rho(\Lambda,L)-\mathcal H^\rho(\lambda_\rho,L)<\infty,
\end{equation*}
so the measure is locally finite. \cite[Proposition~2.18]{dauvergne2026directed} gives
continuous spatial profiles at each deterministic parameter.
The continuous base profile and rectangle increments on countable
dense grids give a Borel $\mathscr D(\lambda_\rho)$-valued
representing pair. The rescaled PNG pair is Borel, Proposition~\ref{prop:png-horizon-finite-slope-limit} and Lemma~\ref{lem:png-horizon-positive-measure-upgrade} therefore give
convergence in $\mathscr D(\lambda_\rho)$.
The additional convergence at every prescribed finite parameter list
is the separate conclusion of
Proposition~\ref{prop:png-horizon-finite-slope-limit}.
\end{proof}

\subsection{Convergence to the half-space KPZ fixed point}
\label{subsec:kpz-fixed-point-convergence}

The transition limit requires separate critical and subcritical
analyses of the boundary operators. In the critical window we prove
uniform estimates for their Airy and first-hitting factors and pass
to the Fredholm Pfaffian. At a fixed subcritical rate the boundary
normalization changes, and stationary comparison supplies terminal
wall observations. Localization and spatial tightness give profile
limits; Markov iteration and short-time bounds then prove
Theorem~\ref{thm:intro-half-space-png-kpz-fixed-point-limit}, following
the transition-based approach of \cite[Section~3]{matetski2021kpz}.

Use
\begin{equation*}
 \delta=\varepsilon^{1/2},\qquad
 T=t\delta^{-3},\qquad X_i=x_i\delta^{-2},\qquad
 c_\delta=e^{\rho\delta}.
\end{equation*}
We index the operator estimates by $\delta$ and the process estimates
by $\varepsilon=\delta^2$, with profiles reindexed by
$f_\varepsilon=f_{\sqrt\varepsilon}$.
We first treat finite $\rho$ and then a fixed subcritical boundary
rate. We use the Markov property and space--time continuity of the
half-space KPZ fixed point from
\cite[Theorem~4.11]{dauvergne2026directed}.

\subsubsection{Critical-window transition probabilities}
\label{subsubsec:critical-transition-limit}

The boundary poles approach the cubic saddle under critical scaling.
We control the corresponding residue terms and the first-hitting
factors in weighted operator norms before passing to the Pfaffian.

\begin{theorem}[Fixed-time transition limit for general initial data]
\label{thm:png-general-data-kpz-limit}
Fix $\rho\in\mathbb R$ and $t>0$. Let
$h_{0,\delta}\in\mathcal X_{\mathrm{adm}}$ be deterministic initial
profiles, and put
$f_\delta(s)=\delta h_{0,\delta}(s\delta^{-2})$ for $s\geq0$.
Suppose $f_\delta\to f\in\mathrm{UC}_+^{\mathbb R}$ locally in
hypograph and, for one finite constant $C\geq0$,
\begin{equation}
 f_\delta(s)\leq C(1+s),\qquad s\geq0,
 \label{eq:png-general-data-common-ceiling}
\end{equation}
for all sufficiently small $\delta$.  Then, for every
$0\leq x_1<\cdots<x_m$,
\begin{equation}
 \left(\delta\left[h^{(c_\delta)}
      (T,x_i\delta^{-2};h_{0,\delta})-2T\right]\right)_{i=1}^m
 \Longrightarrow
 \left(\mathfrak h^\rho(t,x_i;f)\right)_{i=1}^m.
 \label{eq:png-general-data-kpz-limit}
\end{equation}
Here $\mathfrak h^\rho$ is the half-space KPZ fixed point of
\cite[Theorem~4.11]{dauvergne2026directed}; its transition Pfaffian is
Zhang's formula \cite{zhang2024tasep}, displayed below in the convention of
Section~\ref{subsec:intro-kpz-universality}.
\end{theorem}

\paragraph{\textbf{The transition formula and Airy factors}}
\label{subsubsec:png-companion-kernel}

We first recall the symmetric kernel of
\cite[Theorem~1.1 and Section~1.3]{chenZhangPNG}.
The free, hitting, and boundary factors encode propagation, the initial
profile, and the wall, respectively.
For this recall, fix \(c\geq0\), \(h\in\mathcal X_{\mathrm{adm}}\),
\(T>0\), \(0\leq X_1<\cdots<X_m\), and \(r_j\in\mathbb Z\).

On functions on \(\mathbb Z\), set
\begin{equation*}
 \begin{gathered}
 (\tau_\pm f)(u)=f(u\pm1),\qquad
 \nabla=\tfrac12(\tau_+-\tau_-),\qquad
 \Delta=\tau_++\tau_--2I,\\
 \mathsf S_{T,X}=e^{2T\nabla+X\Delta},\qquad
 1_r(u,v)=\mathbf1_{\{u>r\}}\mathbf1_{\{u=v\}}.
 \end{gathered}
\end{equation*}
Discrete kernel products use counting measure.  The boundary kernel is
\begin{equation*}
 \mathsf H_c(u,v)=
 \begin{cases}
  \tfrac12\operatorname{sgn}(v-u)c^{|v-u|-1},&u\ne v,\\
  0,&u=v,
 \end{cases}
 \qquad 0^0=1.
\end{equation*}
For \(0\leq c<1\), it represents
\((I-c\tau_+)^{-1}\nabla(I-c\tau_-)^{-1}\); the displayed kernel
continues it to every \(c\geq0\).
The lower and upper one-sided inverse kernels are
\begin{equation*}
 \mathsf H_{c,-}^{-1}
 =2\sum_{k\geq0}\tau_-^{2k+1}(I-c\tau_+)(I-c\tau_-),\qquad
 \mathsf H_{c,+}^{-1}
 =-2\sum_{k\geq0}\tau_+^{2k+1}(I-c\tau_+)(I-c\tau_-).
\end{equation*}
Their supports are \(u\geq v\) and \(u\leq v\), respectively, and
\(\mathsf H_c^*=-\mathsf H_c\),
\((\mathsf H_{c,-}^{-1})^*=-\mathsf H_{c,+}^{-1}\).
These one-sided kernels are used on the corresponding exponentially
weighted spaces; grouped products for \(c\geq1\) are interpreted by
analytic continuation.

Let \(\mathsf X\) be the continuous-time symmetric nearest-neighbor walk
on \(\mathbb Z\), with jumps \(+1\) and \(-1\) each occurring at rate
one. Its generator is \(\Delta\). We use its
lower-semicontinuous version at jump times.  Write \(\mathbb P_{a,u}\)
for its law started from \(\mathsf X(a)=u\). For \(0\leq a\leq b\), put
\begin{equation*}
 \begin{aligned}
 \mathsf P_{a,b}^{h}(u,v)
 &:=\mathbb P_{a,u}\bigl(
   \mathsf X(y)\leq h(y)\text{ for some }y\in[a,b],\
   \mathsf X(b)=v\bigr),\\
 \mathsf Q_{a,b}^{h}&:=e^{(b-a)\Delta}-\mathsf P_{a,b}^{h},\\
 \mathsf T_{a,b}&:=e^{a\Delta}\mathsf P_{a,b}^{h}
                         e^{-b\Delta}\mathsf H_{c,-}^{-1},
 \qquad \mathsf T_b:=\mathsf T_{0,b}.
 \end{aligned}
\end{equation*}
Thus hitting is tested on the closed interval, and
\(\mathsf P_{a,a}^{h}\) is multiplication by
\(\mathbf1_{\{u\leq h(a)\}}\).
For an empty interval \(a>b\geq0\), set
\(\mathsf P_{a,b}^{h}=\mathsf T_{a,b}=0\).
For nonzero real \(p,q\), define
\begin{equation*}
 \mathsf W_{p,q}:=
 \begin{cases}
  \mathsf T_q-\mathsf T_p^*
     +\mathsf T_p^*\mathsf H_c\mathsf T_q,&p,q>0,\\
  \mathsf T_{-p,q},&p<0<q,\\
  -\mathsf T_{-q,p}^*,&q<0<p,\\
  0,&p,q<0.
 \end{cases}
\end{equation*}
When \(T\ne X_j\) for every \(j\), the symmetric kernel uses
\begin{equation}
 \mathcal S_j=
 \begin{pmatrix}\mathsf S_{T,-X_j}&0\\0&\mathsf H_c\mathsf S_{T,X_j}\end{pmatrix},
 \qquad
 \mathsf V_{ij}=
 \begin{pmatrix}
  \mathsf W_{T-X_i,T-X_j}&\mathsf W_{T-X_i,T+X_j}\\
  \mathsf W_{T+X_i,T-X_j}&\mathsf W_{T+X_i,T+X_j}
 \end{pmatrix}.
 \label{eq:png-companion-grouped-factors}
\end{equation}
Set \(J_0=\left(\begin{smallmatrix}0&1\\-1&0\end{smallmatrix}\right)\).
The projected kernel is \(K_{ij}=\mathsf R_{ij}+\widetilde K_{ij}\), where
\begin{equation*}
 \begin{aligned}
 \mathsf R_{ij}
 &=1_{r_i}
 \begin{pmatrix}
  \mathbf1_{\{j<i\}}e^{(X_i-X_j)\Delta}
    &e^{X_i\Delta}\mathsf H_ce^{X_j\Delta}\\
  0&\mathbf1_{\{i<j\}}e^{(X_j-X_i)\Delta}
 \end{pmatrix}1_{r_j},\\
 \widetilde K_{ij}
 &=1_{r_i}J_0\mathcal S_i^*\mathsf V_{ij}\mathcal S_j1_{r_j}.
 \end{aligned}
\end{equation*}
Here the adjoint includes height variables and internal components.
Let \(J_m\) act by \(J_0\) at each observation-height pair.
The identity \(\mathsf W_{p,q}^*=-\mathsf W_{q,p}\) gives
\((J_mK)^*=-J_mK\).
Then
\begin{equation*}
 \mathbb P_h^{(c)}
 \bigl(h^{(c)}(T,X_j)\leq r_j,\ 1\leq j\leq m\bigr)
 =\operatorname{Pf}(J_m+J_mK)
\end{equation*}
on \(\bigoplus_{j=1}^m\ell^2(\mathbb Z)\otimes\mathbb C^2\), with
\(\operatorname{Pf}(J_m)=1\).
The exponentially weighted formulation is justified by
\cite[Proposition~1.2]{chenZhangPNG}.
For sufficiently small \(\delta\), our scaling has \(T>X_m\).
The two hitting compositions that we will estimate directly are
\begin{equation*}
 \mathsf T_i^-:=\mathsf T_{T-X_i}\mathsf S_{T,-X_i},\qquad
 \mathsf T_i^+:=\mathsf T_{T+X_i}\mathsf H_c\mathsf S_{T,X_i}.
\end{equation*}

\paragraph{\textbf{The limiting half-space KPZ kernel}}
We state Zhang's limiting transition kernel \cite{zhang2024tasep}\footnote{We use the kernel formulation in a newer revision of this
work kindly provided by X.~Zhang and not yet available online.}
in the conventions used here. Fix \(t>0\), \(\rho\in\mathbb R\), and
an upper-semicontinuous \(h_0\not\equiv-\infty\) with a linear upper bound.
Here \(0\leq y_1<\cdots<y_m\) and
\(s_1,\ldots,s_m\in\mathbb R\). The parameter and variable
correspondence with our PNG scaling is
\begin{equation*}
 \varrho=-\rho,\qquad h_0=f,\qquad y_j=x_j,\qquad s_j=r_j.
\end{equation*}
The corresponding microscopic height thresholds are
\(\lfloor2T+\delta^{-1}s_j\rfloor\).
The affine boundary parameter of
\cite[Section~1.3]{chenZhangPNG} is
\(\delta^{-1}(1-c_\delta)\to\varrho\).
The sign change in the boundary parameter agrees with the convention
of \cite[Section~4.2]{dauvergne2026directed}.

Let \(D=\partial_u\), and let \(\mathcal C_a^{\pi/3}\) be the two Airy
rays \(a+re^{\pm i\pi/3}\), oriented from the lower ray to the upper
ray.  For \(t>0\), define
\begin{equation*}
 \mathbf S^{t,x}(u,v)
 =\frac{1}{2\pi i}\int_{\mathcal C_a^{\pi/3}}
    e^{tw^3/3+xw^2+(u-v)w}\,\mathrm dw,\qquad a>0.
\end{equation*}
Write \(\mathbf H=\mathbf H_\varrho\) for the boundary kernel and
use the lower-triangular inverse convention
\begin{equation*}
 \begin{gathered}
 \mathbf H(u,v)
 =\tfrac12\operatorname{sgn}(v-u)e^{-\varrho|u-v|},\\
 \mathbf H^{-1}=D^{-1}(\varrho-D)(\varrho+D),
 \qquad D^{-1}(u,v)=\mathbf1_{\{u>v\}}.
 \end{gathered}
\end{equation*}
In particular, the inverse Airy factor has the kernel
\begin{equation*}
 (\mathbf H^{-1}\mathbf S)^{t,x}(u,v)
 =\frac{1}{2\pi i}\int_{a+i\mathbb R}
   \frac{\varrho^2-w^2}{w}
   e^{tw^3/3+xw^2+(u-v)w}\,\mathrm dw,
 \qquad a>\max(0,-x/t).
\end{equation*}
For the boundary-modified Airy factor, let \(\mathcal C_\varrho\)
be a clockwise loop enclosing only \(w=\varrho\), and set
\begin{equation*}
 \mathbf S_{\mathbf H}^{t,x}(u,v)
 =\frac{1}{2\pi i}
  \int_{\mathcal C_a^{\pi/3}
    +(1-\frac12\mathbf1_{\{\varrho=0\}})\mathcal C_\varrho}
  \frac{w\,e^{tw^3/3+xw^2+(u-v)w}}
       {(\varrho+w)(\varrho-w)}\,\mathrm dw,
 \qquad a>|\varrho|.
\end{equation*}
The loop has half weight when \(\varrho=0\).  Equivalently,
\begin{equation}
 \mathbf S_{\mathbf H}^{t,x}(u,v)
 =\frac{1}{2\pi i}\int_{a+i\mathbb R}
   \frac{w}{\varrho^2-w^2}
   e^{tw^3/3+xw^2+(u-v)w}\,\mathrm dw 
 +\frac12
   e^{t\varrho^3/3+x\varrho^2+\varrho(u-v)},
 \qquad a>|\varrho|+|x|/t .
 \label{eq:png-limiting-boundary-airy}
\end{equation}
At \(\varrho=0\) the last term is \(1/2\).
The notation \(\mathbf S_{\mathbf H}\) always uses this contour
prescription; for nonpositive \(\varrho\), it does not assert an
absolutely convergent convolution on unweighted \(L^2\).

The initial profile enters through Brownian first hitting.
Let \(B\) be Brownian motion of variance \(2\), started at \(u\), and put
\(\tau=\inf\{s\geq0:B(s)\leq h_0(s)\}\).  Define
\begin{equation*}
 \begin{aligned}
 \mathbf S^{\mathrm{hypo}(h_0),t,x}(u,v)
 &:=\mathbb E_u\!
   \left[\mathbf S^{t,x-\tau}(B(\tau),v);
         \tau<\infty\right],\\
 (\mathbf H^{-1}\mathbf S)^{\mathrm{hypo}(h_0),t,x}(u,v)
 &:=\mathbb E_u\!
   \left[(\mathbf H^{-1}\mathbf S)^{t,x-\tau}(B(\tau),v);
         \tau<\infty\right].
 \end{aligned}
\end{equation*}
Let \(U_s=e^{sD}\), so \((U_s\phi)(u)=\phi(u+s)\).  Set
\begin{equation*}
 \begin{aligned}
 \mathbf F_{-j}&=\mathbf S^{t,-y_j}U_{-s_j},
 &\mathbf F_j&=-\mathbf S_{\mathbf H}^{t,y_j}U_{-s_j},\\
 \mathbf G_{-j}&=-(\mathbf H^{-1}\mathbf S)
              ^{\mathrm{hypo}(h_0),t,-y_j}U_{-s_j},
 &\mathbf G_j&=\mathbf S^{\mathrm{hypo}(h_0),t,y_j}U_{-s_j},
 \end{aligned}
\end{equation*}
and collect the four factors in
\begin{equation*}
 \boldsymbol{\mathcal S}_j
 :=\begin{pmatrix}\mathbf F_{-j}&\mathbf F_j\\\mathbf G_{-j}&\mathbf G_j\end{pmatrix}.
\end{equation*}
The translations move the height thresholds \(s_j\) to zero.
The complete limiting kernel is
\begin{equation*}
 \mathbf K_{ij}^{\mathrm{fp}}
 =\mathbf R_{ij}^{\mathrm{fp}}+\widetilde{\mathbf K}_{ij}^{\mathrm{fp}},
 \qquad
 \widetilde{\mathbf K}_{ij}^{\mathrm{fp}}
 =J_0\boldsymbol{\mathcal S}_i^*J_0
   \begin{pmatrix}I&\mathbf H\\0&I\end{pmatrix}
   \boldsymbol{\mathcal S}_j,
\end{equation*}
where
\begin{equation*}
 \mathbf R_{ij}^{\mathrm{fp}}
 =\begin{pmatrix}
 \mathbf1_{\{j<i\}}U_{s_i}e^{(y_i-y_j)D^2}U_{-s_j}
 &
 -U_{s_i}e^{y_iD^2}\mathbf H e^{y_jD^2}U_{-s_j}\\
 0&
 \mathbf1_{\{i<j\}}U_{s_i}e^{(y_j-y_i)D^2}U_{-s_j}
 \end{pmatrix}.
\end{equation*}
The adjoint includes the height variable and internal matrix indices;
the final height variables lie in \([0,\infty)\).
With \(J_m\) acting by \(J_0\) at each observation-height pair, the
transition formula is
\begin{equation*}
 \mathbb P\!\left(
   \mathfrak h^\rho(t,y_j;h_0)\leq s_j,\ 1\leq j\leq m\right)
 =\operatorname{Pf}(J_m+J_m\mathbf K^{\mathrm{fp}})
       _{L^2(\{1,\ldots,m\}\times[0,\infty);\mathbb C^2)}.
\end{equation*}
We use Zhang's Pfaffian transition laws \cite{zhang2024tasep} together with
their realization as the continuous Markov process $\mathfrak h^\rho$ in
\cite[Theorem~4.11]{dauvergne2026directed}.

\begin{lemma}[Uniform Airy-factor limits]
\label{lem:png-finite-wedge-airy-estimates}
Put $n_\delta=\lfloor2T\rfloor$ and, for $u,v\in\mathbb R$, put
$u_\delta=\lfloor u/\delta\rfloor$ and
$v_\delta=n_\delta+\lfloor v/\delta\rfloor$.
For $x$ in a fixed compact set, the following limits are locally uniform:
\begin{align}
 \delta^{-1}\mathsf S_{T,x\delta^{-2}}(u_\delta,v_\delta)
 &\longrightarrow\mathbf S^{t,x}(u,v),
 \label{eq:png-free-airy-limit}\\
 (\mathsf H_{c_\delta}\mathsf S_{T,x\delta^{-2}})(u_\delta,v_\delta)
 &\longrightarrow\mathbf S_{\mathbf H}^{t,x}(u,v),
 \notag\\
 \delta^{-2}(\mathsf H_{c_\delta,-}^{-1}\mathsf S_{T,x\delta^{-2}})(u_\delta,v_\delta)
 &\longrightarrow(\mathbf H^{-1}\mathbf S)^{t,x}(u,v).
 \notag
\end{align}
If the terminal variable is restricted to $v\geq r$, then each free or
inverse factor, multiplied on the right by $e^{\eta v}$, has row
$L^2$ norm at most $C_a e^{au}$, for every sufficiently large fixed
$a>\eta$.  The factor involving the boundary kernel multiplied on the right by
$e^{-\theta v}$ has row norm at most
$C_a(e^{au}+e^{\varrho u})$ when $\theta>|\varrho|$.
These bounds are uniform for sufficiently small $\delta$, for the
limits, and for $r$ in a compact set.  The corresponding weighted rows
converge in $L^2$ for each fixed $u$.
\end{lemma}

\begin{proof}
The Laurent symbol of $\mathsf S_{T,X}$ is
$\exp\{T(z-z^{-1})+X(z+z^{-1}-2)\}$.  On the circle
$z=e^{\delta(a+i\xi)}$, $|\xi|\leq\pi/\delta$, its centered exponent is
\begin{equation*}
 E_\delta(w)=2T\sinh(\delta w)
 +2x\delta^{-2}(\cosh(\delta w)-1)
 +(u_\delta-v_\delta)\delta w.
\end{equation*}
It converges locally to $tw^3/3+xw^2+(u-v)w$.  Moreover,
\begin{equation*}
 \operatorname{Re}\bigl(E_\delta(a+i\xi)-E_\delta(a)\bigr)
 =-2\bigl[T\sinh(\delta a)+x\delta^{-2}\cosh(\delta a)\bigr]
       (1-\cos(\delta\xi)) 
 \leq-\kappa\xi^2
\end{equation*}
for $a$ sufficiently large, with $\kappa>0$ independent of small
$\delta$.  Also
$\operatorname{Re}E_\delta(a)\leq C_a+a(u-v)$.
The factor $\delta^{-1}$ in \eqref{eq:png-free-airy-limit} cancels
the Jacobian $dz/z=\delta\,dw$.  Dominated convergence proves the
first limit and the bound $C_a e^{a(u-v)}$.

The symbols of the boundary operator and its lower-triangular inverse are, respectively,
\begin{equation*}
 h_c(z)=\frac{z-z^{-1}}{2(1-cz)(1-cz^{-1})},\qquad
 u_c(z)=\frac{2z^{-1}(1-cz)(1-cz^{-1})}{1-z^{-2}}.
\end{equation*}
On the same contour,
\begin{equation*}
 \delta h_{c_\delta}(e^{\delta w})\longrightarrow
 \frac{w}{\varrho^2-w^2},\qquad
 \delta^{-1}u_{c_\delta}(e^{\delta w})\longrightarrow
 \frac{\varrho^2-w^2}{w}.
\end{equation*}
For $a>|\varrho|$, the first multiplier is bounded uniformly on the
contour.  The second is bounded by $C(1+|\xi|^2)$: on
$|\delta\xi|\leq\pi/2$ this follows by expanding the factors at $1$;
on the complementary arc the denominator has modulus at least
$1-e^{-2\delta a}$ and the bound is absorbed by $1+|\xi|^2$.
The Gaussian bound therefore applies to both integrals.

For $c\ne0$, put
\begin{equation*}
 F_\delta(c):=
 e^{T(c^{-1}-c)+X(c+c^{-1}-2)}c^{v_\delta-u_\delta},
 \qquad X=x\delta^{-2}.
\end{equation*}
For $0<c<1$, the Laurent integrand at the selected pole satisfies
\begin{equation*}
 \operatorname*{Res}_{z=c^{-1}}
 \left[h_c(z)e^{T(z-z^{-1})+X(z+z^{-1}-2)}
                         z^{u_\delta-v_\delta}\frac1z\right]
 =-\frac{F_\delta(c)}{2c}.
\end{equation*}
Moving the positively oriented circle past this pole therefore adds
$F_\delta(c)/(2c)$ to the centered boundary factor.
For fixed $T,X$, expansion of
$\mathsf S_{T,X}=e^{-2X}e^{(T+X)\tau_+}e^{(X-T)\tau_-}$ gives
$\sum_n|\mathsf S_{T,X}(0,n)|R^{|n|}<\infty$ for every $R\geq1$.
Thus the explicit convolution $\mathsf H_c\mathsf S_{T,X}$ is entrywise
entire in $c$, with normal convergence on compact $c$-sets.
For the fixed contour radius $R=e^{a\delta}>1$, the contour integral
plus the displayed residue term is holomorphic on $R^{-1}<|c|<R$.
Equality for $R^{-1}<c<1$ extends to $c=c_\delta$ when $a>|\varrho|$,
including $c=1$. This is a fixed-$\delta$ continuation; the preceding
contour estimates supply the bounds uniform in small $\delta$.
At $c=c_\delta$, $F_\delta(c_\delta)=e^{E_\delta(\varrho)}$.
The added term therefore converges to the last term of
\eqref{eq:png-limiting-boundary-airy}; its absolute value is bounded
by $C e^{\varrho(u-v)}$.  This proves the second limit, including the residue of the boundary pole.  Integration against the indicated exponential
weights gives the row bounds and, by dominated convergence, weighted
row convergence.
\end{proof}

\paragraph{\textbf{First-hitting measures and uniform estimates}}

The initial profile enters through the joint time and height of first
contact. We prove convergence of these hitting measures under local
hypograph convergence and check the walk's state at a contact occurring
at a jump. Weighted row estimates at arbitrary hitting positions
remove compact cutoffs using the common linear upper bound.
Write $M_\beta$ for multiplication by $e^{\beta u}$
and $b^\pm=\max(\pm b,0)$ for the positive and negative parts.

\begin{lemma}[Uniform kernel estimates at arbitrary hitting positions]
\label{lem:png-general-hit-row-bound}
Fix $t>0$, a compact set of $x$ values, a lower bound $r$ for the
terminal height, and a finite $C\geq0$.  Put $n_\delta=\lfloor2T\rfloor$,
$b_\delta=\lfloor b/\delta\rfloor$, and
$v_\delta=n_\delta+\lfloor v/\delta\rfloor$.  For
$0\leq s\leq t\delta^{-1}+x$, define
\begin{align*}
 g_{\delta,+}^{x}(s,b;v)
 =\delta^{-1}\mathsf S_{T,(x-s)\delta^{-2}}(b_\delta,v_\delta),\qquad
 g_{\delta,-}^{x}(s,b;v)
 =\delta^{-2}(\mathsf H_{c_\delta,-}^{-1}
                 \mathsf S_{T,(x-s)\delta^{-2}})(b_\delta,v_\delta).
\end{align*}
Extend these rows by zero when $s>t\delta^{-1}+x$.
For every $\mu>0$ and every fixed terminal exponential weight,
there are $c,K,\delta_0>0$ such that both prescribed weighted row
$L^2([r,\infty))$ norms satisfy
\begin{equation}
 \|g_{\delta,\pm}^{x}(s,b;\cdot)M_{\pm\theta}\|_2
 \leq K e^{-cs^3-\mu b^-},\qquad
 b\leq C(1+s),\quad0<\delta<\delta_0.
 \label{eq:png-general-hit-row-bound}
\end{equation}
Either sign of a fixed weight can be used for either row.  The limiting rows
$\mathbf S^{t,x-s}(b,\cdot)$ and
$(\mathbf H^{-1}\mathbf S)^{t,x-s}(b,\cdot)$ satisfy the same
bound.  On compact sets of $(s,b)$ the weighted microscopic rows
converge uniformly in $L^2$ to these limiting rows.
\end{lemma}

\begin{proof}
Choose $\kappa$ larger than $\mu$, all absolute terminal weight
parameters, $|\rho|$, and $1$.  The centered Laurent phase on
$z=e^{\eta+i\delta\xi}$, where $a=\eta/\delta$, is
\begin{equation*}
 E_\delta(a+i\xi)
 =2t\delta^{-3}\{\sinh(\delta(a+i\xi))
                         -\delta(a+i\xi)\}
  +2(x-s)\delta^{-2}\{\cosh(\delta(a+i\xi))-1\}
  +(b-v+\gamma_\delta\delta)(a+i\xi),\qquad
       |\gamma_\delta|\leq3.
\end{equation*}
The bounded correction $\gamma_\delta$ records the height rounding;
its contribution at the positive real point is $O(\eta)$.
Put $S=s-x$.  If $S\geq0$ and $r_\delta=\delta S/t<1$, choose
\begin{equation}
 \eta=\max\{\kappa\delta,2\operatorname{arctanh}r_\delta\}.
 \label{eq:png-general-hit-adaptive-contour}
\end{equation}
If $S<0$, choose $\eta=\kappa\delta$.  On the adaptive part of
the contour one has the exact identities
\begin{align}
 2t\delta^{-3}(\sinh\eta-\eta)
       -2S\delta^{-2}(\cosh\eta-1)
 &=4t\delta^{-3}(r_\delta-\operatorname{arctanh}r_\delta)
 \leq-\frac{4S^3}{3t^2},
 \label{eq:png-general-hit-cubic-phase}\\
 t\delta^{-3}\sinh\eta-S\delta^{-2}\cosh\eta
 &=S\delta^{-2}.
 \label{eq:png-general-hit-gaussian-coefficient}
\end{align}
On the part where $\eta=\kappa\delta$, the choice in
\eqref{eq:png-general-hit-adaptive-contour} implies
$S\leq t\delta^{-1}\tanh(\kappa\delta/2)$ when $S\geq0$.
Consequently the coefficient on the left of
\eqref{eq:png-general-hit-gaussian-coefficient} is at least
$t\delta^{-3}\tanh(\kappa\delta/2)$.
It has a uniform positive multiple of $\delta^{-2}$ as a lower
bound in both cases.  Thus, for $|\xi|\leq\pi/\delta$,
\begin{equation}
 \operatorname{Re}(E_\delta(a+i\xi)-E_\delta(a))
 =-2\{t\delta^{-3}\sinh\eta
                   -S\delta^{-2}\cosh\eta\}(1-\cos\delta\xi)
 \leq-k\xi^2
 \label{eq:png-general-hit-gaussian-bound}
\end{equation}
with a fixed $k>0$.

The linear upper bound is absorbed by the cubic decay, including near
the causal endpoint.  Since $a\geq\kappa$,
$b\leq C(1+s)$ and $v\geq r$ give
\begin{equation*}
 a(b-v)\leq a\{C(1+s)+|r^-|\}-\mu b^- -\kappa v.
\end{equation*}
For $r_\delta\leq1/2$ one has $a\leq K(1+s)$, so the additional
positive term is quadratic in $s$ and is absorbed by the negative
cubic in \eqref{eq:png-general-hit-cubic-phase}; the region where
$a=\kappa$ is bounded in $s$.  For $r_\delta>1/2$, there is an
absolute $k_1>0$ such that
$\operatorname{arctanh}r_\delta-r_\delta \geq k_1\operatorname{arctanh}r_\delta$.
Here $s\leq t\delta^{-1}+x$ and
$a=2\delta^{-1}\operatorname{arctanh}r_\delta$.  The contribution from the linear upper bound and the rounding term
are bounded by
$K\delta^{-2}\operatorname{arctanh}r_\delta$.
For small $\delta$ they absorb at most half of the negative
$4k_1t\delta^{-3}\operatorname{arctanh}r_\delta$.
Since $s=O(\delta^{-1})$ in this region, the remainder bounds a
negative fixed multiple of $s^3$.  These observations show that
\begin{equation}
 \operatorname{Re}E_\delta(a)
 \leq K-cs^3-\mu b^- -\kappa v.
 \label{eq:png-general-hit-real-bound}
\end{equation}
In the region $r_\delta>1/2$, split the remaining negative phase
into two equal parts before obtaining this estimate.  One part
bounds $-cs^3$; the other retains the stronger additional term
$-k_2\delta^{-3}\operatorname{arctanh}r_\delta$.

The contour integral for the rescaled $S$ kernel has multiplier $1$.  For the scaled lower
inverse, the multiplier is $\delta^{-1}u_{c_\delta}(z)$, with
$u_c$ as in the Airy-factor lemma.  Uniformly on these contours,
$|\delta^{-1}u_{c_\delta}(e^{\eta+i\delta\xi})| \leq K(a+\xi^2)$.
Indeed, for $\eta\leq1$, the numerator is bounded by
$K(\eta^2+\min(\delta|\xi|,1)^2)$ and
$|1-z^{-2}|\geq k\eta$; division by $\delta$ gives
$K(a+\xi^2/a)$.  For $\eta\geq1$ the unscaled multiplier is
bounded, and $\delta^{-1}\leq a$.
Gaussian integration of \eqref{eq:png-general-hit-gaussian-bound}
therefore gives a prefactor at most $K(1+a)$ in front of
\eqref{eq:png-general-hit-real-bound}.  When $r_\delta\leq1/2$
this is polynomial in $s$; when $r_\delta>1/2$ its logarithm is
absorbed by the unused negative
$k\delta^{-3}\operatorname{arctanh}r_\delta$ term.  Decreasing
$c$ once gives the uniform pointwise bound
\begin{equation*}
 |g_{\delta,\pm}^{x}(s,b;v)|
 \leq K e^{-cs^3-\mu b^- -\kappa v}.
\end{equation*}
Integration against the squared terminal weights proves
\eqref{eq:png-general-hit-row-bound}.

At $s=t\delta^{-1}+x$, the free operator is
$\mathsf S_{T,-T}=e^{2T(I-\tau_-)}$, which is lower triangular.
The lower inverse is also lower triangular.  The common linear upper bound gives $b_\delta\leq K\delta^{-2}$ for a
constant $K$ independent of $\delta$, whereas the terminal height
is $2t\delta^{-3}+O(\delta^{-1})$.  Both rows therefore vanish
for sufficiently small $\delta$.  This proves the bound also at
the causal endpoint without an infinite contour radius.

For the continuum rows, choose
$a=\max\{\kappa,2(s-x)^+/t\}$.  The real cubic phase has the
same negative cubic bound after the linear upper bound is included,
and its off-real part is $-[ta-(s-x)]\xi^2$, with a fixed positive
coefficient.  This proves the limiting bound.  On compact
$(s,b)$ sets, the Airy contour argument applies uniformly,
and the terminal exponential bound permits uniform removal of
the $v$ tail.  This proves the stated weighted row convergence.
\end{proof}

\begin{lemma}[Convergence of first-hitting measures]
\label{lem:png-general-hit-measures}
Let $\mathsf X_\delta(s)=\delta \mathsf X(s\delta^{-2})$ be the rescaled walk,
started at $\delta\lfloor u/\delta\rfloor$, and let
\begin{equation*}
 \tau_\delta=\inf\{s\geq0:\mathsf X_\delta(s)\leq f_\delta(s)\},
 \qquad
 \nu_{\delta,u}(ds,db)
 =\mathbb P_u(\tau_\delta\in ds,\mathsf X_\delta(\tau_\delta)\in db,
                                  \tau_\delta<\infty).
\end{equation*}
Define $\tau$ and $\nu_u$ analogously for variance-$2$ Brownian
motion and the limiting barrier $f$.  Under the hypotheses of
Theorem~\ref{thm:png-general-data-kpz-limit},
$\nu_{\delta,u}\to\nu_u$ vaguely for Lebesgue-almost every $u$.
For every $c>0$, $\mu>0$, and $0<\lambda<\mu$,
\begin{equation}
 \int e^{-cs^3-\mu b^-}\nu_{\delta,u}(ds,db)
 \leq K e^{-\lambda|u|}
 \label{eq:png-general-hit-measure-envelope}
\end{equation}
uniformly for sufficiently small $\delta$.  The same bound holds
for $\nu_u$.
\end{lemma}

\begin{proof}
The lower-semicontinuous path and upper-semicontinuous barrier
give a closed set of hitting times, so every finite first hit is attained. The
first hitting times are stopping times for the completed natural
filtrations. Along any sequence $\delta\downarrow0$, the invariance principle
couples the centered walk increments $\mathsf X_\delta^\circ:=\mathsf X_\delta-\mathsf X_\delta(0)$
and a centered variance-$2$ Brownian motion $B^\circ$ so that
$\mathsf X_\delta^\circ\to B^\circ$ uniformly on compact intervals
almost surely. Set $B(s):=u+B^\circ(s)$ for the Brownian motion
started at $u$. Replacing the walk by a continuous interpolation changes its value by at most $\delta$.
Choose a deterministic countable dense set $\mathcal T$ of positive
continuity points of the extended-real upper-semicontinuous function
$f$.  Such points are dense because $f$ is a Baire-one function
(one can first truncate it below).  Local hypograph convergence then
also holds after restriction to $[0,q]$, for every $q\in\mathcal T$.
For each such $q$, the common compact upper bound and hypograph
convergence imply
\begin{equation*}
 A_\delta(q):=\sup_{0\leq s\leq q}
                     \{f_\delta(s)-\mathsf X_\delta^\circ(s)\}
 \longrightarrow
 A(q):=\sup_{0\leq s\leq q}\{f(s)-B^\circ(s)\}.
\end{equation*}
For the upper bound, take a convergent subsequence of approximate
maximizers and use the hypograph limsup inequality; if $A(q)=-\infty$,
this rules out every finite lower bound. Hypograph recovery
sequences give the reverse inequality.

The event $\{\tau_\delta\leq q\}$ is the height comparison
$\{\delta\lfloor u/\delta\rfloor\leq A_\delta(q)\}$.
For each Brownian sample, exclude the countably many real values
$A(q)$, $q\in\mathcal T$.  For every remaining $u$ these height
comparisons converge at all $q\in\mathcal T$.  Fubini's theorem
therefore gives this conclusion almost surely for almost every
deterministic $u$.  The exceptional deterministic heights can be
chosen to be the atoms of the countably many random variables
$A(q)$, so they depend only on the limiting Brownian law and not on
the chosen sequence.  Monotonicity of the hitting events and density
of $\mathcal T$ imply $\tau_\delta\to\tau$, including the value
$\infty$.  At a finite limiting time, compact uniform convergence
of the paths gives
$\mathsf X_\delta(\tau_\delta)\to B(\tau)$.  Thus bounded continuous
functions with compact support of the hitting pair converge in
expectation, proving vague convergence.  The lower-semicontinuous
jump convention changes interpolated values by at most $\delta$ and
does not change this argument.

For the estimate, the exponential martingales of the walk have
compensator
$\psi_\delta(\ell)=2\delta^{-2}(\cosh(\delta\ell)-1)$.
For any fixed $\ell$, $\psi_\delta(\ell)$ is uniformly bounded
and converges to $\ell^2$.  Optional stopping, first at bounded
times and then using Fatou's lemma, gives
\begin{equation*}
 \mathbb E_u[e^{\pm\ell \mathsf X_\delta(\tau_\delta)
            -\psi_\delta(\ell)\tau_\delta};\tau_\delta<\infty]
 \leq K e^{\pm\ell u}.
\end{equation*}
A first hit at an upward jump from height $a\in\delta\mathbb Z$
would require its time to be a left-isolated point of the deterministic
closed set $\{f_\delta\geq a\}$. These points, over all such $a$,
are countable and are almost surely avoided by the jump times.
Thus at every finite first hit the lower-semicontinuous path agrees
with the right-continuous Markov state; at a downward jump both use
the post-jump value. Initial height rounding contributes at most
$e^{\ell\delta}$ to the martingale bound.
At a hit, $b\leq C(1+s)$ by
\eqref{eq:png-general-data-common-ceiling}.  If $u\leq0$, use
$e^{-\mu b^-}\leq e^{\lambda b}$ and
$\sup_{s\geq0}e^{-cs^3+\psi_\delta(\lambda)s}<K$.
The positive exponential-martingale bound gives
$Ke^{\lambda u}$.  If $u\geq0$, use
\begin{equation*}
 e^{-cs^3-\mu b^-}
 \leq e^{-\lambda b-\psi_\delta(\lambda)s}
       e^{-cs^3+\lambda C(1+s)+\psi_\delta(\lambda)s},
\end{equation*}
whose second factor is uniformly bounded.  The negative
exponential-martingale bound gives $Ke^{-\lambda u}$.
Brownian motion has compensator $\ell^2$, so the same proof
applies verbatim to its hitting measure.
\end{proof}

\paragraph{\textbf{Convergence of the Fredholm Pfaffians}}
\label{subsubsec:png-kernel-limit-preliminaries}

The $11$ and $22$ blocks converge in trace norm, while the $12$ and
$21$ blocks require Hilbert--Schmidt convergence. We prove these
mixed estimates and construct the normalized Pfaffian as an analytic
square root, fixing its branch and justifying the conjugations below.
Apply the transition formula \cite[Theorem~1.1]{chenZhangPNG} with thresholds
$\lfloor2T+\delta^{-1}r_i\rfloor$.  Embed the height lattice into
$L^2(\mathbb R)$ by constant functions on cells of length $\delta$,
center heights by $\lfloor2T\rfloor$, and use the internal conjugation
$C_\delta=\operatorname{diag}(i\delta^{1/2},-i\delta^{-1/2})$.
At the kernel level the rescaling is
\begin{equation*}
 \widehat K_\delta(u,v)
 :=\delta^{-1}C_\delta
 K(\lfloor2T\rfloor+\lfloor u/\delta\rfloor,
   \lfloor2T\rfloor+\lfloor v/\delta\rfloor)C_\delta^{-1}.
\end{equation*}
Thus the four entries of the rescaled kernel carry the factors
$\delta^{-1},-1,-\delta^{-2},\delta^{-1}$, respectively.
Indeed, the isometric embedding sends $(v_n)$ to
$\delta^{-1/2}v_n$ on $[\delta n,\delta(n+1))$, giving the density
factor $\delta^{-1}$. On finite height cutoffs using the same projection in both components,
$C_\delta^{\mathsf T}J_0C_\delta=J_0$ and $\det C_\delta=1$;
the congruence on the skew kernel is
\(J_mK\mapsto C_\delta^{-\mathsf T}(J_mK)C_\delta^{-1}\),
which preserves the normalized Fredholm Pfaffian.

Choose $|\rho|<\theta_1<\cdots<\theta_m$, and conjugate the two
components at observation $i$ by $e^{-\theta_i v}$ and
$e^{\theta_i v}$, respectively.  A superscript $\delta$ on a kernel
below denotes this centered, embedded and conjugated version.
Write $\mathsf R=(\mathsf R_{ij})$ for the free block
matrix. Its diagonal blocks are heat kernels between distinct observation
positions. Write $p_\ell^\delta(u,v)$ for the embedded cell kernel at
microscopic time $\ell\delta^{-2}$ and
\begin{equation*}
 p_\ell(z)=\frac{1}{\sqrt{4\pi\ell}}e^{-z^2/(4\ell)},\qquad \ell>0,
\end{equation*}
for its Gaussian limit.  The local central limit theorem and the same shifted-circle
estimate as above give, for every fixed $b>0$,
$|p_\ell^\delta(u,v)|\leq C_b e^{-b|u-v|}$, uniformly for small
$\delta$. This is $\delta^{-1}$ times the lattice heat kernel at
$(\lfloor u/\delta\rfloor,\lfloor v/\delta\rfloor)$; its semigroup
identity holds on the embedded cell subspace. The Gaussian limit and
this bound give convergence in every weighted difference-variable
$L^2$ norm, with cell-rounding errors absorbed in the same bound.
To obtain trace convergence for the $11$ block with
$i>j$, put $\ell=x_i-x_j>0$ and choose
$\theta_j<\eta<\theta_i$.  Factor the conjugated heat block as
\begin{equation*}
 \bigl(\mathbf1_{>r_i}M_{-\theta_i}p_{\ell/2}^\delta M_\eta\bigr)
 \bigl(M_{-\eta}p_{\ell/2}^\delta M_{\theta_j}\mathbf1_{>r_j}\bigr).
\end{equation*}
The two factors converge in Hilbert--Schmidt norm: integrating out
the difference variable leaves, respectively, the integrable weights
$e^{2(\eta-\theta_i)u}$ on $u>r_i$ and
$e^{2(\theta_j-\eta)v}$ on $v>r_j$.
Their products therefore converge in trace norm.  The $22$ blocks
follow by transposition.  The $O(\delta)$ changes of the scaled thresholds caused by
lattice rounding vanish under the same bounds.  The remaining block of $\mathsf R$ is the boundary kernel convolved with heat
kernels. The random-walk exponential moments and the explicit
formula for $\mathsf H_c$ give, after the stated conjugation,
\begin{equation*}
 |\mathsf R_{12}^{\delta}(i,u;j,v)|
 \leq C e^{-\theta_i u-\theta_j v+|\rho||u-v|}
\end{equation*}
on the projected half-lines. This bound is square-integrable there
because $\theta_i,\theta_j>|\rho|$.
These kernels converge almost everywhere by the Gaussian limit
for the heat factors and the pointwise limit of $\mathsf H_c$.  This
bound proves Hilbert--Schmidt convergence, including $x_i=x_j=0$.
For operators, write $\|\cdot\|_1$ for the trace norm and
$\|\cdot\|_2$ for the Hilbert--Schmidt norm.  Thus the blocks of $\mathsf R$
converge in the corresponding components of the mixed norm
\begin{equation*}
 \|A\|_{1|2}=\|A_{11}\|_1+\|A_{22}\|_1
              +\|A_{12}\|_2+\|A_{21}\|_2.
\end{equation*}

Write $J=J_m$ on the continuum observation-height space
$\bigoplus_{i=1}^mL^2(\mathbb R)\otimes\mathbb C^2$.
For kernels in this mixed Schatten class such that $(JK)^{\mathsf T}=-JK$,
the generalized normalized Pfaffian is constructed as follows. Choose real
finite-rank orthogonal projections $P_N\uparrow I$ on the scalar observation-height
space and put $\Pi_N=P_N\otimes I_2$, $K_N=\Pi_NK\Pi_N$.
The normalized finite-dimensional Pfaffian $p_N$ satisfies
\begin{equation*}
 p_N(0)=1,\qquad
 p_N(z)^2=\det{}_2(I+zK_N)e^{z\tau_N},\qquad
 \tau_N=\operatorname{tr}(K_N)_{11}+\operatorname{tr}(K_N)_{22}.
\end{equation*}
Here $\det_2(I+A)=\det((I+A)e^{-A})$ is the regularized
Fredholm determinant for Hilbert--Schmidt $A$.
Mixed-norm convergence of $K_N$ gives locally uniform convergence of
these squares. The finite-dimensional inequality
$|(1+w)e^{-w}|\leq e^{|w|^2/2}$ and Schur triangularization give
\begin{equation*}
 |p_N(z)|\leq\exp\left\{
 \tfrac14|z|^2\|K\|_2^2+
 \tfrac12|z|(\|K_{11}\|_1+\|K_{22}\|_1)\right\}.
\end{equation*}
Montel's theorem gives an entire limiting square root, unique near zero
by its value $1$ there and hence unique everywhere. Thus the full sequence
converges locally uniformly to a projection-independent $p_K$; define
$\Pf(J+JK):=p_K(1)$. The same argument proves continuity under
$\|\cdot\|_{1|2}$ convergence. It agrees with the ordinary normalized
Pfaffian for trace-class kernels.

For a bounded invertible component-diagonal $D$ with
$D^{\mathsf T}JD=J$, similarity preserves $\det_2$ and each diagonal
trace, and $JDKD^{-1}=D^{-\mathsf T}(JK)D^{-1}$ remains skew.
Uniqueness of the normalized square root gives $p_{DKD^{-1}}=p_K$.
Apply this on bounded height intervals and then remove the height cutoffs
in mixed norm to justify the exponential conjugations whenever the
resulting kernels lie in the stated mixed class.

\begin{proposition}[Identification of the limiting transition kernel]
\label{prop:png-kpz-kernel-identification}
For an upper-semicontinuous $f\not\equiv-\infty$ with a linear upper
bound, the limiting projected kernel has Fredholm Pfaffian
\begin{equation*}
 \mathbb P\bigl(\mathfrak h^\rho(t,x_i;f)\leq r_i,
                          \ 1\leq i\leq m\bigr).
\end{equation*}
\end{proposition}
\begin{proof}
For sufficiently small \(\delta\), every \(T\pm X_j\) is positive.
Substitute
\(\mathsf W_{p,q}=\mathsf T_q-\mathsf T_p^*
 +\mathsf T_p^*\mathsf H_c\mathsf T_q\)
in the symmetric kernel
\eqref{eq:png-companion-grouped-factors}.
Its columns involve only \(\mathsf S_{T,-X_j}\),
\(\mathsf H_c\mathsf S_{T,X_j}\), and the two direct hitting compositions
\(\mathsf T_j^-,\mathsf T_j^+\).
Replace these columns by the Airy and Brownian first-hitting factors
displayed above. They give the four factors
\(\mathbf F_{-j},\mathbf F_j,\mathbf G_{-j},\mathbf G_j\) in the
limiting kernel, after the
internal conjugation and the threshold translations \(U_{s_j}\).
Their signs follow from \(\mathsf H_c^*=-\mathsf H_c\) and the
internal conjugation by \(C_\delta\).
Expanding the two-by-two product therefore gives
\begin{equation*}
 \widetilde{\mathbf K}_{ij}^{\mathrm{fp}}
 =J_0\boldsymbol{\mathcal S}_i^*J_0
   \begin{pmatrix}I&\mathbf H\\0&I\end{pmatrix}
   \boldsymbol{\mathcal S}_j,
\end{equation*}
and the free block gives \(\mathbf R_{ij}^{\mathrm{fp}}\).
Thus the complete limit is the displayed
\(\mathbf K^{\mathrm{fp}}\), with
\(y_j=x_j\), \(s_j=r_j\), \(h_0=f\), and \(\varrho=-\rho\).
The contour prescription, including the half residue at zero, agrees
with the boundary Airy factor
\eqref{eq:png-limiting-boundary-airy}.
\end{proof}

\begin{proof}[Proof of Theorem~\ref{thm:png-general-data-kpz-limit}]
The transition formula applies to each $h_{0,\delta}$. Since
$f\not\equiv-\infty$, a bounded scaled interval contains a finite
initial point for small $\delta$, making every observed height finite.

The direct hitting compositions have rescaled horizons
\begin{equation*}
 \delta^2(T\pm X_i)=t/\delta\pm x_i\longrightarrow\infty.
\end{equation*}
The strong Markov decomposition at the first hitting time gives, with $q=T\pm X_i$,
\begin{equation*}
 \mathsf P_{0,q}^{h_{0,\delta}}e^{-q\Delta}\mathsf S_{T,\pm X_i}(u,v)
 =\mathbb E_u[\mathsf S_{T,\pm X_i-\tau}(\mathsf X(\tau),v);
                                                    \tau\leq q].
\end{equation*}
The same identity holds with $\mathsf H_{c_\delta,-}^{-1}$ inserted,
since this convolution commutes with the free operators.  For a
fixed $\delta$, these equalities can first be tested on finitely
supported vectors in the height space with squared norm
$\sum_n e^{2\eta n}|v_n|^2$, $\eta<0$.  The lower-triangular branch of the inverse
is bounded there, its geometric shift series converging in operator
norm; the free operators are bounded there as well.  The hitting
operator is dominated entrywise by the heat kernel and hence is
bounded on the same space.  The strong Markov identity and these
weighted bounds permit passage to the completion, justifying the
interchanges with the inverse. The argument for the convention at jumps in
Lemma~\ref{lem:png-general-hit-measures} justifies using the
right-continuous Markov state at the first hit.

For the plus composition \(\mathsf T_i^+\), the cancellation
\(\mathsf H_{c,-}^{-1}\mathsf H_c=I\) is first used for
\(0\leq c<1\) on a compatible exponentially weighted space.
For \(c\geq1\), the complete grouped identity is its analytic
continuation, justified by \cite[Theorem~1.1 and Appendix~A]{chenZhangPNG}.

Consequently, \(\delta^{-2}\mathsf T_i^-\) and
\(\delta^{-1}\mathsf T_i^+\) are integrals of the rows
in Lemma~\ref{lem:png-general-hit-row-bound} against
$\nu_{\delta,u}$.  On a compact set of hitting pairs, vague
convergence in Lemma~\ref{lem:png-general-hit-measures} and uniform
$L^2$ row convergence give convergence of these integrals in
the weighted row $L^2$ norm for almost every $u$.  A continuous
compact cutoff in $(s,b)$ justifies this Hilbert-space-valued
statement by uniform approximation by finite sums of scalar
continuous functions times fixed row vectors.  The cutoff can be
removed uniformly: the row bound tends uniformly to zero when
$s\to\infty$ or $b^-\to\infty$, while $b^+\leq C(1+s)$ on
the support.  The measures have total mass at most one.
Moreover, \eqref{eq:png-general-hit-measure-envelope} gives the
uniform row bound $Ke^{-\lambda|u|}$, with $\lambda$ as large
as required by first choosing $\mu>\lambda$.

Apply the transition formula, lattice embedding, and symplectic
conjugation above. For linear scattering
products, the rank-one trace norm is dominated by
$K(e^{a u}+e^{|\rho||u|})e^{-\lambda|u|}$ for $\lambda>\max(a,|\rho|)$.
For quadratic products it is dominated by
$K e^{|\rho|(|u|+|v|)}e^{-\lambda(|u|+|v|)}$ for $\lambda>|\rho|$.
These bounds are integrable on the full internal height line(s).
Dominated convergence therefore gives trace-norm convergence of
all scattering blocks along the given hypograph-convergent sequence;
the domination depends only on the common linear upper bound.
Together with the free-block estimates, this gives convergence of the
full kernel in $\|\cdot\|_{1|2}$.

Proposition~\ref{prop:png-kpz-kernel-identification} identifies the
limiting kernel. Mixed-norm Pfaffian continuity gives convergence of
the joint distribution functions to the KPZ transition law.
The same weighted bounds, with thresholds in a compact set, prove
continuity in the limiting thresholds, completing the weak convergence.
\end{proof}

For $0\leq a_1<\cdots<a_k$ and $b_1,\ldots,b_k\in\mathbb R$, set
\begin{equation}
 f(x)=\begin{cases}b_j,&x=a_j,\\-\infty,&x\notin\{a_1,\ldots,a_k\},\end{cases}
 \qquad
 h_{0,\delta}(y)=
 \begin{cases}\lfloor b_j/\delta\rfloor,&y=a_j\delta^{-2},\\
 -\infty,&\text{otherwise}.
 \end{cases}
 \label{eq:png-finite-wedge-scaled-data}
\end{equation}

The profiles in \eqref{eq:png-finite-wedge-scaled-data} converge locally
in hypograph with a common linear upper bound, so
Theorem~\ref{thm:png-general-data-kpz-limit} gives their finite-wedge
transition limit.

\begin{corollary}[Independent random initial data]
\label{cor:png-random-general-data-kpz-limit}
The conclusion of Theorem~\ref{thm:png-general-data-kpz-limit} also
holds when the initial profiles are random and independent of the PNG
environment, provided $f_\delta\Rightarrow f$ locally in hypograph,
$f\not\equiv-\infty$ almost surely, and their global linear upper bounds
have tight random constants. Explicitly, assume
\begin{equation*}
 \lim_{A\to\infty}\limsup_{\delta\downarrow0}
 \mathbb P\bigl(\exists s\geq0:
                    f_\delta(s)>A(1+s)\bigr)=0.
\end{equation*}
The limiting fixed point is started from the random $f$ independently
of its driving directed landscape.
\end{corollary}
\begin{proof}
Under a deterministic linear upper bound, the theorem gives sequential
convergence of the conditional transition probabilities. A Skorokhod
coupling of the initial profiles and bounded convergence therefore
give joint convergence with the terminal height vector.
For upper bounds with general tight random constants, use
\begin{equation*}
 f_\delta^A(s)=\min\{f_\delta(s),
            \delta\lfloor A(1+s)/\delta\rfloor\}.
\end{equation*}
At finite values this differs by at most $\delta$ from clipping by the
continuous ceiling, and both clips equal $-\infty$ where $f_\delta=-\infty$.
Thus it converges locally in hypograph to $\min\{f,A(1+\cdot)\}$.
The underlying integer cutoff is the floor of an increasing continuous
affine function, so it is upper-semicontinuous and has finitely many
jumps on compacts, preserving the allowed initial-data class.
The clipped profiles differ from the originals
with probability tending uniformly to zero as $A\to\infty$.
The same tight-envelope property passes to the limiting profile, so
its clipping error also tends to zero as $A\to\infty$.  Letting $A\to\infty$
proves the claim.
\end{proof}

\subsubsection{Subcritical finite-wedge transitions}
\label{subsubsec:png-subcritical-degeneration}

Fix $c_\star\in[0,1)$ before sending $\delta$ to zero. The limiting
scattering form contains a derivative and, when an initial wedge is
at zero, a jump contribution. We first treat finite wedges at positive
terminal positions; stationary comparison will include the terminal wall.
The boundary factors require the internal conjugation below, with
$d=1-c_\star$:
\begin{equation}
 C_{\delta,\star}
 =\operatorname{diag}\bigl(i d\delta^{-1/2},
                  -i\delta^{1/2}/d\bigr).
 \label{eq:png-subcritical-conjugation}
\end{equation}
Together with the lattice density factor $\delta^{-1}$, it multiplies
the four entries by $\delta^{-1}$, $-d^2\delta^{-2}$,
$-d^{-2}$, and $\delta^{-1}$, respectively.  Its determinant is one
and it preserves $J_0$.

\begin{proposition}[Subcritical finite-wedge limit away from the terminal wall]
\label{prop:png-subcritical-positive-endpoint-limit}
For the finite-wedge initial data in
\eqref{eq:png-finite-wedge-scaled-data}, any fixed $c_\star\in[0,1)$,
$t>0$, and $0<x_1<\cdots<x_m$, one has
\begin{equation*}
 \left(\delta\bigl[h^{(c_\star)}(t\delta^{-3},
      x_i\delta^{-2};h_{0,\delta})-2t\delta^{-3}\bigr]\right)_{i=1}^m
 \Longrightarrow
 \left(\mathfrak h^{-\infty}(t,x_i;f)\right)_{i=1}^m.
\end{equation*}
Here the initial wedges may include a wedge at the spatial origin.
\end{proposition}

For finite wedges the random-walk hitting operator can be reduced to
finitely many observation times.  Set $A_j=a_j\delta^{-2}$ and
$b_{j,\delta}^{\mathrm{lat}}=\lfloor b_j/\delta\rfloor$.  For all sufficiently small $\delta$,
every endpoint $T\pm X_i$ exceeds $A_k$.  With $A_0=0$, put
\begin{equation*}
 D_{j,\delta}
 =\left[\prod_{\ell=1}^{j-1}
 e^{(A_\ell-A_{\ell-1})\Delta}
 \mathbf1_{>b_{\ell,\delta}^{\mathrm{lat}}}\right]
 e^{(A_j-A_{j-1})\Delta}\mathbf1_{\leq b_{j,\delta}^{\mathrm{lat}}}.
\end{equation*}
The product is ordered by increasing $\ell$ from left to right;
an empty product is $I$.
The first-hitting decomposition is the finite identity
\begin{equation}
 \mathsf P_{0,q}^{h_{0,\delta}}e^{-q\Delta}
   =\sum_{j=1}^kD_{j,\delta}e^{-A_j\Delta},
 \qquad q\geq A_k.
 \label{eq:png-finite-wedge-first-hit-product}
\end{equation}
This conditions on the first successful height test, including time zero.

We interpret the limiting derivative operator as follows.
For row kernels $f,g$, write $(f\otimes g)(v,w)=f(v)g(w)$ for their
rank-one tensor kernel, with component indices included in $v,w$.
Let $F,G$ be vector-valued row functions with
exponentially integrable row norms and first derivatives.  They may
have one common jump at a height $b$.  Define the following trace-class bilinear expression
\begin{equation}
 \mathfrak q_D(F,G)
 =\frac12\int_{\mathbb R\setminus\{b\}}
          \bigl(F(u)\otimes G'(u)-F'(u)\otimes G(u)\bigr)\,du 
 +\frac12\bigl(F(b-)\otimes G(b+)
                    -F(b+)\otimes G(b-)\bigr).
 \label{eq:png-subcritical-skew-derivative-form}
\end{equation}
Omit the jump term for continuous rows. This midpoint convention
is needed when the first initial wedge is at zero.

\begin{lemma}[Subcritical operator limits]
\label{lem:png-subcritical-factor-estimates}
Use the height coordinates of
Lemma~\ref{lem:png-finite-wedge-airy-estimates}.  Then, locally uniformly,
\begin{align}
 \delta^{-1}\mathsf S_{T,x\delta^{-2}}&\longrightarrow\mathbf S^{t,x},\notag\\
 d^2\delta^{-2}\mathsf H_{c_\star}\mathsf S_{T,x\delta^{-2}}
       &\longrightarrow D\mathbf S^{t,x},\notag\\
 d^{-2}\mathsf H_{c_\star,-}^{-1}\mathsf S_{T,x\delta^{-2}}
       &\longrightarrow D^{-1}\mathbf S^{t,x}.
 \label{eq:png-subcritical-free-factor-limits}
\end{align}
Here $D^{-1}$ is the lower-triangular branch used above.  With terminal
weights $e^{\eta v}$ or $e^{-\theta v}$ on $v>r$, all three rows have
norm at most $C_a e^{au}$ for a sufficiently large fixed $a$; the
same statement holds with any fixed number of derivatives in either
height variable.  For embedded microscopic rows, these derivatives
mean scaled finite differences, not derivatives of the piecewise
constant embedding.  The corresponding scaled hitting rows converge to
\begin{equation*}
 \mathcal G_{i,+}^{\star}(u,v)
 =\mathbb E_u[\mathbf S^{t,x_i-\tau}(B(\tau),v)
                   \mathbf1_{\{\tau<\infty\}}],\qquad
 \mathcal G_{i,-}^{\star}(u,v)
 =\mathbb E_u[(D^{-1}\mathbf S)^{t,-x_i-\tau}(B(\tau),v)
                   \mathbf1_{\{\tau<\infty\}}],
\end{equation*}
where $\tau=\min\{a_j:B(a_j)\leq b_j\}$, with $\min\varnothing=\infty$,
using the microscopic normalizations $\delta^{-1}\mathsf T_i^+$ and
$d^{-2}\mathsf T_i^-$.  They and their first derivatives in the starting
height obey a bound
\begin{equation}
 Q_1(u)=C_{a,b}\bigl(e^{au}\mathbf1_{\{u<0\}}
                      +e^{-bu}\mathbf1_{\{u\geq0\}}\bigr),
 \label{eq:png-subcritical-derivative-envelope}
\end{equation}
for arbitrarily large fixed $a,b>0$, away from the sole possible jump
at $u=b_1$ when $a_1=0$.  At that jump both one-sided row norms are
bounded.  Discrete first differences divided by $\delta$ obey the
same estimate on each side of the jump.

If $F_\delta,G_\delta$ are any two such hitting rows, then
\begin{equation}
 \frac{d^2}{2}\sum_{n\in\mathbb Z}\sum_{k\geq1}c_\star^{k-1}
 \bigl[F_\delta(\delta n)\otimes G_\delta(\delta(n+k))
       -F_\delta(\delta(n+k))\otimes G_\delta(\delta n)\bigr]
 \longrightarrow\mathfrak q_D(F,G)
 \label{eq:png-subcritical-quadratic-limit}
\end{equation}
in trace norm.  At $c_\star=0$ only the term $k=1$ is retained.
\end{lemma}

\begin{proof}
On $z=e^{\delta(a+i\xi)}$ the centered Airy exponent has the same
uniform Gaussian estimate as in
Lemma~\ref{lem:png-finite-wedge-airy-estimates}.  For sufficiently small
$\delta$ the entire contour lies in the annulus
$c_\star<|z|<c_\star^{-1}$ when $0<c_\star<1$.  At $c_\star=0$
the Laurent contour is unrestricted.  In both cases no boundary pole
is crossed.  The symbols
satisfy
\begin{equation*}
 \frac{d^2}{\delta}h_{c_\star}(e^{\delta w})\longrightarrow w,
 \qquad
 \frac{\delta}{d^2}u_{c_\star}(e^{\delta w})\longrightarrow\frac1w.
\end{equation*}
On the whole contour the first normalized multiplier is bounded by
$C(1+|\xi|)$.  The second is bounded by a constant depending on $a$
and $c_\star$: near either of the possible zeros $z=\pm1$, use
$|1-z^{-2}|\geq1-e^{-2a\delta}$: the numerator is bounded and
the factor $\delta$ cancels this denominator bound.  Away from
these two arcs the denominator is bounded away from zero.  The
uniform Gaussian exponent bound also suppresses the arc near $-1$.  Multiplying by
powers of $w$ supplies the asserted derivative estimates.  Gaussian
domination proves the factor limits and weighted row convergence.

In the finite first-hitting identity
\eqref{eq:png-finite-wedge-first-hit-product}, the heat tests converge
by the local central limit theorem. Positive-time Brownian heights have
continuous densities; at time zero exclude $u=b_1$. The uniform
moment formula
\begin{equation*}
 \mathbb E e^{\eta\delta[\mathsf X(a_j\delta^{-2})-\mathsf X(0)]}
 =e^{2a_j\delta^{-2}(\cosh(\eta\delta)-1)}
\end{equation*}
and the preceding row bounds justify truncating hit heights and passing
to weighted row limits. For the starting-height
derivative, suppose first that $a_1>0$.  In each first-hitting term its
first factor is the heat kernel at time $a_1$.  Its embedded first
difference divided by $\delta$ converges to the derivative of the
variance-two Gaussian kernel and has uniform exponentially weighted
bounds.  One can see this directly from its Fourier integral by
inserting $(e^{\delta w}-1)/\delta$ and using the heat exponent's
Gaussian bound.  Integrating this estimate against the remaining
finite heat tests and the terminal Airy kernel gives
\eqref{eq:png-subcritical-derivative-envelope}.  Its positive-height
tail follows from exponential tilting together with the event that
one of the finitely many tested heights is at most $\max_jb_j$;
its negative-height tail follows from the Airy-kernel row estimate.  These
are the same exponential-moment bounds used for the undifferentiated
rows.  If $a_1=0$, split the starting height at $b_1$.  On the lower
side the row is the corresponding Airy row.  On the upper side the
first possible positive hit time is $a_2>0$, so the preceding heat
argument applies; if there is no second wedge, the upper row is zero.
The bond $n=\lfloor b_1/\delta\rfloor$ crossing the jump is
handled in the crossing sum below.
This proves the estimates and one-sided convergence.

For the last assertion, use skew-symmetry before taking absolute
values.  Away from the possible jump the tensor difference for an
increment $k\delta$ has trace norm bounded by $Ck\delta$ times the
product of the bounds in
\eqref{eq:png-subcritical-derivative-envelope}.  Summation over $n$
therefore leaves an integrable bound of order $Ck$, uniformly in
small $\delta$.  If the increment crosses the jump there are at most
$k$ starting lattice sites; their tensors have uniformly bounded row
norms, with at most a factor $e^{Ck\delta}$.  The sums in $k$ are
uniformly summable because $c_\star<1$.  For each fixed $k$, the
ordinary Riemann sum off the jump converges to
$\frac{k}{2}\int(F\otimes G'-F'\otimes G)$,
and the $k$ crossing terms converge to
$\frac{k}{2}(F(b-)\otimes G(b+)-F(b+)\otimes G(b-))$.
The coefficient identity
$d^2\sum_{k\geq1}kc_\star^{k-1}=1$ gives
\eqref{eq:png-subcritical-quadratic-limit}.  All these statements hold
in trace norm, since $\|f\otimes g\|_1=\|f\|_2\|g\|_2$.
\end{proof}

\begin{proof}[Proof of Proposition~\ref{prop:png-subcritical-positive-endpoint-limit}]
Apply the transition formula with
\eqref{eq:png-subcritical-conjugation}.  Choose
$0<\theta_1<\cdots<\theta_m$ and use the same reciprocal exponential
weights in the two components as in the finite-parameter proof.  Linear scattering products
converge in trace norm by
Lemma~\ref{lem:png-subcritical-factor-estimates}; the row norms of the Airy and hitting factors have integrable products.  For the quadratic scattering
products retain their antisymmetrized sums and use
\eqref{eq:png-subcritical-quadratic-limit}.  The diagonal heat-kernel
blocks converge in trace norm by the two-Hilbert--Schmidt-factor
argument already given.  Finally, the boundary block converges
in Hilbert--Schmidt norm to
$-D e^{(x_i+x_j)D^2}$.
Indeed, $x_i+x_j>0$ supplies a Gaussian Fourier bound, the normalized
boundary multiplier converges to $w$, and the terminal exponential
weights make the projected kernel square integrable.  Thus the whole
kernel converges in the mixed norm $\|\cdot\|_{1|2}$ to a kernel
$K^{\star}$ determined by the three limits in
\eqref{eq:png-subcritical-free-factor-limits} and the form
\eqref{eq:png-subcritical-skew-derivative-form}.

To identify the limiting Pfaffian, take the critical-window kernel
with $\rho=-R$, so its source parameter is \(\varrho=R\), and apply
the symplectic similarity multiplying its
$12$ entry by $R^2$ and its $21$ entry by $R^{-2}$.
Fix $0<\theta_1<\cdots<\theta_m$, independently of $R$,
and a sufficiently large fixed $a>\theta_m$. For $R>2a$, the boundary
convolution on $a+i\mathbb R$ has multiplier $w/(R^2-w^2)$,
without a residue term. Since
\begin{equation*}
 \left|\frac{R^2w}{R^2-w^2}\right|\leq2|w|,\qquad
 R^{-2}\mathbf H_R^{-1}\mathbf S
   =D^{-1}\mathbf S-R^{-2}D\mathbf S,
\end{equation*}
Gaussian contour domination gives
\begin{equation*}
 R^2(\mathbf H_R\mathbf S)\longrightarrow D\mathbf S,
 \qquad
 R^{-2}(\mathbf H_R^{-1}\mathbf S)
       \longrightarrow D^{-1}\mathbf S.
\end{equation*}
The same estimates and finite heat tests give the one-sided derivative
bounds above uniformly in large $R$. These fixed exponential conjugations preserve the Fredholm Pfaffians.
For the quadratic term the exact antisymmetrized expression is
\begin{equation*}
 \frac{R^2}{2}\int_0^\infty e^{-Rs}\int_{\mathbb R}
 \bigl[F(u)\otimes G(u+s)-F(u+s)\otimes G(u)\bigr]\,du\,ds.
\end{equation*}
Off the jump, divide the bracket by $s$ and use the derivative
bounds. The crossing interval has length $s$, giving the jump
term in \eqref{eq:png-subcritical-skew-derivative-form}. Since
$R^2\int_0^\infty se^{-Rs}\,ds=1$, dominated convergence gives
that form in trace norm. Positive terminal positions smooth
the boundary block. Hence the continuum kernels converge in
the same fixed mixed norm to $K^{\star}$.

The proof of \cite[Corollary~5.2]{dauvergne2026directed} gives
finite-dimensional convergence $\mathcal L^{-n}\to\mathcal L^{-\infty}$
on rational endpoints as $n\to\infty$. The uniform modulus in
\cite[Proposition~4.10(1)]{dauvergne2026directed} extends this convergence
to the fixed real endpoints used here. Taking maxima over the fixed finite initial support
identifies the weak limit of the corresponding fixed-point vectors as
$\mathfrak h^{-\infty}(t,x_i;f)$.  Mixed-norm Pfaffian continuity identifies its joint CDF with the Pfaffian of $K^{\star}$.
The weighted trace and Hilbert--Schmidt bounds give threshold
continuity, so Pfaffian convergence proves the microscopic limit at all
fixed thresholds.
\end{proof}

\begin{remark}[Wall observations]
\label{rem:png-subcritical-wall-obstruction}
At $x_i=x_j=0$ the boundary block tends to $-D$, a derivative
of a delta distribution, and has no individual Hilbert--Schmidt limit.
Lemma~\ref{lem:png-subcritical-wall-input} supplies wall observations
by stationary comparison and the one-point wall law.
\end{remark}

\subsubsection{Stationary comparison and wall observations}
\label{subsubsec:stationary-comparison}

We now include terminal observations at the wall and obtain spatial
comparison for tightness. We use exact stationarity, the one-parameter
functional limit of Section~\ref{subsec:stationary-horizon-convergence},
the critical transition limit, and Poisson passage bounds.
Write
\begin{equation*}
 G_c(z,s;y,t):=\mathcal L_{\mathrm{PNG}}^{(c)}(z,s;y,t),\qquad
 G_c(T):=G_c(0,0;0,T),
\end{equation*}
for the unscaled passage metric and its wall-to-wall specialization.

\begin{lemma}[A wall lower tail]
\label{lem:png-wall-lower-tail}
Let $G_c(T)$ be the half-space PNG passage value from $(0,0)$ to
$(0,T)$, with bulk nucleation intensity $2$ and boundary intensity $c$.
There are absolute $C,c_0>0$ such that, for $T\geq1$, $a\geq1$, and
every $c\geq0$,
\begin{equation*}
 \mathbb P\{G_c(T)-2T\leq-aT^{1/3}\}
 \leq C\exp(-c_0a^3).
\end{equation*}
\end{lemma}

\begin{proof}
It suffices to take $c=0$ by Poisson thinning.  The coordinates
$u=t+x$, $v=t-x$ turn the bulk process into intensity-one Poisson
points in $0\leq v<u\leq T$.  Reflecting these points about $u=v$
gives the Poissonized fixed-point-free involution model.  Any increasing
path can be folded into $u\geq v$ without changing its number of points:
each excursion below the diagonal is reflected separately, and a
point and its reflected partner are incomparable.  Thus its longest
increasing subsequence equals $G_0(T)$.

The triangle contains $\operatorname{Poisson}(T^2/2)$ points and
has no diagonal points, so $P_j(T;0):=\mathbb P(G_0(T)\leq j)$
is the CDF in \cite[(4.7)]{baik2001asymptotics}.
By \cite[Theorem~4.2, (4.21)]{baik2001asymptotics}, for $l\geq1$,
\begin{equation*}
 P_{2l}(T;0)=\frac12e^{-T^2/2}
       \bigl(D_l^{--}(T)+D_{l-1}^{++}(T)\bigr).
\end{equation*}
In that paper, $\pi_n(z;T)$ is the monic orthogonal polynomial on
the unit circle for $e^{2T\cos\theta}d\theta/(2\pi)$, and $N_n(T)$
is its squared norm for this measure \cite[(4.10)]{baik2001asymptotics}.
The determinant products in \cite[(4.17)--(4.18)]{baik2001asymptotics}
are
\begin{equation*}
 e^{-T^2/2}D_l^{--}(T)
 =\prod_{j\geq l}N_{2j+2}(T)^{-1}
   \prod_{j\geq l}(1+\pi_{2j+2}(0;T)),
\end{equation*}
and the corresponding product with the minus sign gives
$e^{-T^2/2}D_l^{++}(T)$.  Set $x=2(l-T)l^{-1/3}$.
The bounds \cite[Corollary~7.6, (7.24)--(7.26)]{baik2001asymptotics}
give, for $x\leq-M$,
\begin{equation*}
 e^{-T^2/2}D_l^{--}(T)
       \leq C\exp\{-c|x|^3+C|x|^{3/2}\}.
\end{equation*}
The same upper bound, with $l$ replaced by $l-1$, applies to
$e^{-T^2/2}D_{l-1}^{++}(T)$.  For $T/2\leq l\leq T$, rounding the threshold up to an even
integer gives $|x|\geq a/2$ for large $a$, hence the claimed cubic
bound; enlarged constants cover bounded $a,T$. For $l<T/2$, monotonicity bounds the probability by
$Ce^{-cT^2}$ at $l=\lceil T/2\rceil$. A nonempty event has
$a\leq2T^{2/3}$, so this also implies $Ce^{-c_0a^3}$.
\end{proof}

\begin{lemma}[Full-space compact and diagonal-block moments]
\label{lem:png-fullspace-diagonal-moments}
Let $L^{\mathrm{full}}$ be the full-space Poisson passage metric with bulk intensity
$2$ in PNG coordinates, and define its centered rescaling by
\begin{equation*}
 L_\delta^{\mathrm{full}}(x,s;y,t)
 :=\delta\bigl[L^{\mathrm{full}}(x\delta^{-2},s\delta^{-3};
                     y\delta^{-2},t\delta^{-3})
                    -2(t-s)\delta^{-3}\bigr].
\end{equation*}
For every compact set $K$ of endpoint pairs
with time separation bounded away from zero and every $p>0$,
\begin{equation*}
 \limsup_{\delta\downarrow0}
 \mathbb E[(\sup_K|L_\delta^{\mathrm{full}}|)^p]<\infty.
\end{equation*}
In particular there are $C_p,T_p$ such that
\begin{equation*}
 \mathbb E\left[
  \sup_{0\leq y\leq T^{2/3}}
  \left|\frac{L^{\mathrm{full}}(y,0;y,T)-2T}{T^{1/3}}\right|^p
       \right]\leq C_p,\qquad T\geq T_p.
\end{equation*}
Consequently, for $T\geq T_p$, $R\geq0$, and $a\geq1$,
\begin{align*}
 \sup_{y\in\mathbb R}\mathbb P\bigl(
     L^{\mathrm{full}}(y,0;y,T)-2T>aT^{1/3}\bigr)&\leq C_pa^{-p},\\
 \mathbb P\left\{\inf_{0\leq y\leq R}
     [L^{\mathrm{full}}(y,0;y,T)-2T]<-aT^{1/3}\right\}
 &\leq C_p(1+R/T^{2/3})a^{-p}.
\end{align*}
\end{lemma}

\begin{proof}
The proof of \cite[Theorem~13.8(3)]{dauvergne2021scaling}
establishes compact upper and lower exponential tightness for planar
Poisson passage. For the lower bound it applies Proposition~11.3 and
Corollary~11.4, as in Theorem~12.7, using the one-point input (13.12).
For the upper bound it transfers Corollary~10.8 to planar Poisson
passage by vertical compression and then follows Theorem~12.6.
These bounds imply all polynomial compact-supremum moments.

To check the normalization, the light-cone coordinates $(u,v)=(t+x,t-x)$
send the bulk intensity $2\,dt\,dx$ to $du\,dv$.
In \cite[Corollary~13.6]{dauvergne2021scaling},
the scale parameters are $\tau_T=2T^{2/3}$ and $\chi_T=T^{1/3}$,
with centering $h(u,v)=u+v$.  The cited coordinate map and ours have matrices
\begin{equation*}
 \mathsf M_T^{\mathrm{DV}}=\begin{pmatrix}2T^{2/3}&T\\0&T\end{pmatrix},\qquad
 \mathsf M_T^{\mathrm{PNG}}=\begin{pmatrix}T^{2/3}&T\\-T^{2/3}&T\end{pmatrix},
\end{equation*}
respectively, so
\begin{equation*}
 (\mathsf M_T^{\mathrm{DV}})^{-1}\mathsf M_T^{\mathrm{PNG}}
 =\left(\begin{smallmatrix}1&0\\-T^{-1/3}&1\end{smallmatrix}\right)
 \longrightarrow I.
\end{equation*}
The images of any fixed compact set therefore stay in a common compact
set with positive time separation, and $h\circ \mathsf M_T^{\mathrm{PNG}}(x,t)=2Tt$.
The cited estimates give the first assertion in our normalization.
Take $\delta=T^{-1/3}$ and restrict to
$\{(y,0;y,1):0\leq y\leq1\}$ for the diagonal-block assertion.
Translation invariance and Markov's inequality give the one-point
upper bound. A union bound over at most $1+\lfloor R/T^{2/3}\rfloor$
translated diagonal blocks gives the last bound, including $R=0$.
\end{proof}

\paragraph{\textbf{Wall passage tightness}}
For $R<\infty$ and $0<t_0<t_1<\infty$, the family
\begin{equation*}
 \delta\bigl[G_{\exp(\rho\delta)}(t\delta^{-3})
                 -2t\delta^{-3}\bigr],
 \qquad |\rho|\leq R,\quad t\in[t_0,t_1],\quad\delta\downarrow0,
\end{equation*}
is tight.  The same holds with boundary intensity zero.
It also holds with any fixed boundary intensity $c_\star\in[0,1)$.
Indeed, put $\eta=(t\delta^{-3})^{-1/3}$.  The first displayed random
variable is $t^{1/3}$ times
$\eta[G_{\exp(\kappa\eta)}(\eta^{-3})-2\eta^{-3}]$, where
$\kappa=\rho t^{1/3}$.  Couple by thinning and bound this quantity
between those with parameters $\kappa=\pm Rt_1^{1/3}$.  Theorem~\ref{thm:png-general-data-kpz-limit}, with one wedge at
the wall and time $1$, gives
tightness of both bounding variables. At intensity zero,
Lemma~\ref{lem:png-wall-lower-tail} controls the lower tail, while
$G_0\leq G_1$ and the critical limit control the upper tail.
For fixed $c_\star<1$, the thinning bounds
$G_0\leq G_{c_\star}\leq G_1$ give the fixed-subcritical claim.

\begin{lemma}[Comparison with full-space Poisson passage]
\label{lem:png-fullspace-profile-sandwich}
Fix a time slab $[s,t]$.  Couple $G_c$ and $G_0$ by retaining the
same bulk points, and extend the bulk environment independently to
the negative half-line.  Denote the resulting full-space passage
metric by $L^{\mathrm{full}}$.  Put $W_c=G_c(0,s;0,t)$ and
$W_0=G_0(0,s;0,t)$.  For every causally reachable $z,y\geq0$,
simultaneously,
\begin{equation*}
 L^{\mathrm{full}}(z,s;y,t)-L^{\mathrm{full}}(0,s;0,t)+W_0
 \leq G_0(z,s;y,t)
 \leq G_c(z,s;y,t)
 \leq L^{\mathrm{full}}(z,s;y,t)+W_c-W_0.
\end{equation*}
\end{lemma}

\begin{proof}
Take an optimizing wall-to-wall path in the half-space environment with no boundary atoms. Restrict first to the finite causal cone of
the two endpoint pairs, which contains finitely many Poisson atoms.
  Any path from $(z,s)$ to $(y,t)$ that visits the wall,
or leaves the half-space, must meet this path before and after each
such excursion.  Exchange the two paths between their first and
last intersections.  For comparison of $G_c$ and $G_0$, this puts
all modified boundary points into the wall-to-wall path, and gives
$G_0(0,s;0,t)+G_c(z,s;y,t)
 \leq G_c(0,s;0,t)+G_0(z,s;y,t)$.
For comparison of the restricted metric $G_0$ and $L^{\mathrm{full}}$, the same
exchange gives
$G_0(0,s;0,t)+L^{\mathrm{full}}(z,s;y,t)
 \leq L^{\mathrm{full}}(0,s;0,t)+G_0(z,s;y,t)$.
Finally $G_0\leq L^{\mathrm{full}}$ and $G_0\leq G_c$.  Path intersections can
be included as vertices of weight zero, so this is precisely the
finite directed-graph argument of
\cite[Proposition~3.10 and Corollary~3.11]{dauvergne2026directed}.
\end{proof}

\begin{lemma}[Maximizing starting points at large stationary slope]
\label{lem:png-positive-stationary-exit}
Fix $t>0$ and $b<\infty$.  Let $S_{\delta,\lambda}$ be the pinned
stationary PNG profile with spectral parameter $g=e^{\lambda\delta/2}$ and
boundary intensity $e^{\rho\delta}$, independent of the future
environment, and write
\begin{equation*}
 \widehat S_{\delta,\lambda}(z)
    =\delta S_{\delta,\lambda}(z\delta^{-2}),\qquad
 L_\delta(z;y)=\delta\bigl[
 G_{e^{\rho\delta}}(z\delta^{-2},0;y\delta^{-2},t\delta^{-3})
                  -2t\delta^{-3}\bigr].
\end{equation*}
Let $A_{\delta,\lambda}$ be the event that no maximizing starting point for
$\widehat S_{\delta,\lambda}(z)+L_\delta(z;0)$ lies in $[0,b]$.
Then, for $\lambda>\max(0,2\rho)$,
\begin{equation*}
 \lim_{\lambda\to\infty}\limsup_{\delta\downarrow0}
       \mathbb P(A_{\delta,\lambda}^{\mathrm c})=0.
\end{equation*}
\end{lemma}

\begin{proof}
Set $q_\lambda=\lambda t/2$. The sandwich in
Lemma~\ref{lem:png-fullspace-profile-sandwich} compares the metric
on $[0,b]$ and at $q_\lambda$ with full-space passage, up to tight
boundary-to-boundary passage terms. The compact full-space supremum is tight, and its
one-point limit at $q_\lambda$, plus $q_\lambda^2/t$, has a law
independent of $q_\lambda$.

The one-parameter stationary limit is the two-line Brownian passage
construction with drifts
$(\lambda,-\lambda)$ and bridge
$E_\lambda\sim\operatorname{Exp}(\lambda/2-\rho)$.  The elementary
Lipschitz bound for this construction gives, for $L<\infty$,
\begin{equation*}
 \sup_{0\leq z\leq L}
 |\mathcal H^\rho(\lambda,z)-\lambda z|
 \leq 2\sup_{0\leq z\leq L}|B_1(z)|
      +2\sup_{0\leq z\leq L}|B_2(z)|+2E_\lambda,
\end{equation*}
where the $B_i$ are centered Brownian motions of variance $2$.
At $L=q_\lambda\vee b$ this error is $O_{\mathbb P}(\sqrt\lambda)$.
Consequently, after first sending $\delta$ to zero, the difference
between the value at $q_\lambda$ and the supremum over $[0,b]$
is bounded from below by
\begin{equation*}
 \lambda^2t/4-\lambda b-O_{\mathbb P}(\sqrt\lambda).
\end{equation*}
It is positive with probability tending to one as
$\lambda\to\infty$.  Any maximizing starting point must then exceed $b$.
\end{proof}

\begin{lemma}[One-sided spatial comparison with stationarity]
\label{lem:png-one-sided-spatial-control}
Fix $s<t$, $b<\infty$, and a critical-window or fixed-subcritical
boundary regime. For each sufficiently large positive $\lambda$, one
can construct a height-scaled stationary output profile
$H_{\varepsilon,\lambda,t}$ and an event
$A_{\varepsilon,\lambda}$ such that
\begin{equation*}
 \lim_{\lambda\to\infty}\limsup_{\varepsilon\downarrow0}
 \mathbb P(A_{\varepsilon,\lambda}^{\mathrm c})=0,
\end{equation*}
and, on that event, simultaneously for $0\leq x\leq b$ and
$0\leq y_1\leq y_2\leq b$,
\begin{equation}
 L_\varepsilon^\rho(x,s;y_2,t)-L_\varepsilon^\rho(x,s;y_1,t)
 \leq H_{\varepsilon,\lambda,t}(y_2)
          -H_{\varepsilon,\lambda,t}(y_1).
 \label{eq:png-one-sided-spatial-comparison}
\end{equation}
For fixed $\lambda$, the output profiles are asymptotically
spatially equicontinuous. An analogous inequality controls positive
increments in $x$ after time reversal.
\end{lemma}

\begin{proof}
Time translation allows $s=0$. Independently of future atoms, take
the stationary initial law with
$g_\varepsilon=e^{\lambda\varepsilon^{1/2}/2}$ and the same boundary rate.
Let $H_{\varepsilon,\lambda,t}$ be its height-scaled output recentered
at the wall, and let $A_{\varepsilon,\lambda}$ be the event that
no maximizing starting point for the output at $0$ lies in $[0,b]$.
Lemma~\ref{lem:png-positive-stationary-exit} gives the probability
bound. Its proof applies also at a fixed subcritical boundary rate:
wall passages are bounded above by any critical-window passage and
below by the passage at boundary rate zero, and the stationary horizon limit
at each fixed positive slope is already proved in that regime.

Path exchange makes rightmost maximizing starting points nondecreasing
in the target. On $A_{\varepsilon,\lambda}$, a maximizing path for stationary initial data
to $y_1$ therefore starts to the right of $b$ and crosses an optimizer
from $x\leq b$ to $y_2\geq y_1$. Exchanging continuations proves
\eqref{eq:png-one-sided-spatial-comparison}.
The recentered output has the stationary initial law. Asymptotic equicontinuity of its uninterpolated profiles follows from the uniform coupling in the proof of
Proposition~\ref{prop:png-horizon-finite-slope-limit}: the
upper-semicontinuous profile and its right-continuous representative differ by
at most one, hence by at most $\varepsilon^{1/2}$ after scaling.
Reflecting time on a fixed deterministic slab and reversing passage
paths exchanges their endpoints and preserves the bulk and boundary
Poisson laws. The conventions $(s,t]$ and $[s,t)$ give the same scores
almost surely, since the deterministic endpoint slices contain no
atoms. This gives the comparison in the initial coordinate.
\end{proof}

\begin{lemma}[The subcritical finite-wedge limit at the wall]
\label{lem:png-subcritical-wall-input}
Fix $c_\star\in[0,1)$ and a nonempty finite initial support
$P\subset\mathbb R_{\geq0}$, with fixed initial heights $a_z$, $z\in P$.
The subcritical finite-wedge convergence proved at positive terminal
positions extends to finite terminal collections that include $0$.
\end{lemma}

\begin{proof}
Fix $t>0$ and $b>\max P$. Write
$X_{\varepsilon,z}=L_\varepsilon^{-\infty}(z,0;0,t)$ and
$Y_{\varepsilon,z}(h)=L_\varepsilon^{-\infty}(z,0;h,t)$.
Lemma~\ref{lem:png-one-sided-spatial-control} gives, simultaneously for
$z\in P$ and $0<h<b$,
\begin{equation*}
 X_{\varepsilon,z}\geq Y_{\varepsilon,z}(h)
             -\Omega_{\varepsilon,\lambda}(h)
\end{equation*}
on its comparison event, whose exceptional probability vanishes as
$\varepsilon\downarrow0$ and then $\lambda\to\infty$.
Here
\begin{equation*}
 \Omega_{\varepsilon,\lambda}(h)
 :=\sup_{u,v\in[0,h]}
 |H_{\varepsilon,\lambda,t}(u)-H_{\varepsilon,\lambda,t}(v)|
\end{equation*}
is the oscillation of the stationary output over $[0,h]$;
for fixed $\lambda$ it tends to zero in probability
as $\varepsilon\downarrow0$ and then $h\downarrow0$.

For each $n\in\mathbb N$, boundary thinning also gives
$X_{\varepsilon,z}\leq L_\varepsilon^{-n}(z,0;0,t)$ for small
$\varepsilon$, since $c_\star<e^{-n\varepsilon^{1/2}}$.
Theorem~\ref{thm:png-general-data-kpz-limit} identifies the limit of
this upper bound as $\mathcal L^{-n}(z,0;0,t)$.
The monotone degeneration in the proof of
\cite[Corollary~5.2]{dauvergne2026directed}, extended from rational
endpoints by its uniform modulus bound in Proposition~4.10(1), sends
this law to that of $\mathcal L^{-\infty}(z,0;0,t)$ as $n\to\infty$.
The lower bound has the same limiting law: by
Proposition~\ref{prop:png-subcritical-positive-endpoint-limit},
$Y_{\varepsilon,z}(h)\Rightarrow\mathcal L^{-\infty}(z,0;h,t)$,
and target continuity permits $h\downarrow0$.
Applying bounded increasing Lipschitz functions to these bounds, first
letting $\varepsilon\downarrow0$, then $h\downarrow0$ and
$\lambda\to\infty$ in the lower bound, and $n\to\infty$ in the
upper bound, proves
$X_{\varepsilon,z}\Rightarrow\mathcal L^{-\infty}(z,0;0,t)$,
including $z=0$.

These marginal limits and the one-sided comparison imply
\begin{equation}
 \lim_{h\downarrow0}\limsup_{\varepsilon\downarrow0}
 \mathbb P\!\left(
 \max_{z\in P}|X_{\varepsilon,z}-Y_{\varepsilon,z}(h)|>\eta\right)=0.
 \label{eq:png-subcritical-wall-approximation}
\end{equation}
Indeed, choose a bounded strictly increasing Lipschitz function.
The expected difference of its values at the two tight marginals tends
to zero as $\varepsilon\downarrow0$ and then $h\downarrow0$.
The negative part is bounded by the stationary modulus and the exit
exception. On a compact height interval, a positive difference exceeding
$\eta$ forces a fixed positive function difference. For each fixed height
cutoff, let $\lambda\to\infty$ and then remove the cutoff to obtain
absolute proximity;
a finite union gives \eqref{eq:png-subcritical-wall-approximation}.

Maxima over $P$ preserve this approximation. Apply the positive-terminal
finite-wedge limit jointly at $h$ and the other terminal positions, then
let $h\downarrow0$. Continuity of the limiting finite maximum
$\max_{z\in P}\{a_z+\mathcal L^{-\infty}(z,0;y,t)\}$ proves the claim.
\end{proof}

\subsubsection{Short-time bounds and localization}
\label{subsubsec:short-time-localization}

Short vertical-passage lower bounds and stationary upper bounds will
control nearby observation times and localize the variational supremum,
allowing the passage from fixed-time laws to space--time convergence.

\begin{lemma}[Lower bounds for short vertical passages]
\label{lem:png-uniform-time-buffers}
Fix $b\geq2$ and $\alpha\in(0,1/3)$, and choose
$p$ so that $\gamma:=p(1/3-\alpha)-5/3>0$. There is a constant $C_\alpha$
such that, with
\begin{equation*}
 \mathcal E_{\varepsilon}^{\mathrm{time}}(b,q)=\left\{
 \inf_{\substack{0\leq x\leq b,\ -b\leq s<t\leq b\\t-s\leq q}}
 L_\varepsilon^\rho(x,s;x,t)\geq-C_\alpha q^\alpha\right\},
\end{equation*}
one has
\begin{equation*}
 \lim_{q\downarrow0}\liminf_{\varepsilon\downarrow0}
 \mathbb P(\mathcal E_{\varepsilon}^{\mathrm{time}}(b,q))=1.
\end{equation*}
The conclusion is uniform over all nonnegative boundary rates.
\end{lemma}

\begin{proof}
First fix a duration $r$, put $T=r\varepsilon^{-3/2}$, and use the
full-space extension in Lemma~\ref{lem:png-fullspace-profile-sandwich}.
That comparison gives, with $Y=y\varepsilon^{-1}$ and
$S=s\varepsilon^{-3/2}$,
\begin{equation*}
 G_c(Y,S;Y,S+T)\geq L^{\mathrm{full}}(Y,S;Y,S+T)
                         -L^{\mathrm{full}}(0,S;0,S+T)
                         +G_0(0,S;0,S+T).
\end{equation*}
Split a lower deviation into three fixed fractions and apply
Lemmas~\ref{lem:png-fullspace-diagonal-moments}
and~\ref{lem:png-wall-lower-tail}. Since
$b\varepsilon^{-1}T^{-2/3}=br^{-2/3}$ and
$\varepsilon^{1/2}T^{1/3}=r^{1/3}$, for $T\geq T_p$ and $a\geq1$,
\begin{equation}
 \mathbb P\left(\inf_{0\leq y\leq b}
 L_\varepsilon^\rho(y,s;y,s+r)<-ar^{1/3}\right)
 \leq C_p(1+br^{-2/3})a^{-p}.
 \label{eq:png-short-vertical-lower}
\end{equation}
The bound is uniform in $s$ and the boundary rate. Since $G_c\geq0$,
the event is empty if $a\geq2T^{2/3}$.

Fix $0<q<1$ and use the nested dyadic grids with spacings $q2^{-j}$. At
level $j\geq0$, use the grid of intervals of length $r_j=q2^{-j}$
contained in $[-b-1,b+1]$. There are at most $C_b r_j^{-1}$ such
intervals. Apply \eqref{eq:png-short-vertical-lower} with threshold
$r_j^\alpha$, that is, $a_j=r_j^{\alpha-1/3}$. For all sufficiently
small $q$, $a_j\geq1$. The union bound gives
\begin{equation}
 \mathbb P\left(\exists j,I:\ 
  \inf_{0\leq x\leq b}L_\varepsilon^\rho(x,\inf I;x,\sup I)
       <-r_j^\alpha\right)
 \leq C_b\sum_{j\geq0}
   r_j^{-1}(1+b r_j^{-2/3})r_j^{p(1/3-\alpha)}
 \leq C_{b,p}q^\gamma.
 \label{eq:png-time-buffer-union}
\end{equation}
The estimate applies to every nontrivial event in this sum once
$\varepsilon$ is sufficiently small. Indeed, if
$2r_j\varepsilon^{-1}\leq r_j^\alpha$, the event is empty. Otherwise
\begin{equation*}
 r_j>(\varepsilon/2)^{1/(1-\alpha)},\qquad
 T_j=r_j\varepsilon^{-3/2}
    >2^{-1/(1-\alpha)}
      \varepsilon^{1/(1-\alpha)-3/2}\longrightarrow\infty,
\end{equation*}
where the exponent is negative because $\alpha<1/3$. Moreover
$a_j<2T_j^{2/3}$ in exactly this nontrivial range. Hence \eqref{eq:png-time-buffer-union} applies uniformly to all
nonempty events, and its right side tends to zero as $q\downarrow0$.

On the complementary event, decompose any interval $[s,t]\subset[-b,b]$
of length at most $q$ into maximal dyadic subintervals. There are at most
two of each length $q2^{-j}$. Centered superadditivity holds exactly with the half-open convention
$(s,t]$, including at random composition times. Finite compositions
therefore give the lower bound
$-2\sum_{j\geq0}(q2^{-j})^\alpha =-\frac{2q^\alpha}{1-2^{-\alpha}}$.
For arbitrary real endpoints, first use a finite inner dyadic
approximation. The two unfilled endpoint intervals have total length
at most $2q2^{-n}$; their centered passage values are bounded below by
$-4q2^{-n}\varepsilon^{-1}$. Sending $n\to\infty$ at fixed
$\varepsilon$ removes these remainders. This proves the assertion with
$C_\alpha=2/(1-2^{-\alpha})$.
\end{proof}

Write $U_\varepsilon=\mathfrak h_\varepsilon^\rho$ for the centered
height evolution. The functions $f_\varepsilon$ are rescalings
of PNG initial profiles in $\mathcal X_{\mathrm{adm}}$.

\begin{lemma}[Stationary upper bounds]
\label{lem:png-evolved-linear-bound}
Fix a critical-window or fixed-subcritical boundary regime, and suppose
$f_\varepsilon(z)\leq C(1+z)$ for one finite $C\geq0$.
For every fixed $t>0$,
\begin{equation}
 \lim_{A\to\infty}\limsup_{\varepsilon\downarrow0}
 \mathbb P\bigl(\exists z\geq0:
             U_\varepsilon(t,z)>A(1+z)\bigr)=0.
 \label{eq:png-evolved-linear-bound}
\end{equation}
Moreover, the centered evolution of a dominating stationary profile
has a tight upper supremum on every compact subset of
$(0,\infty)\times[0,\infty)$.
\end{lemma}

\begin{proof}
Put $\bar\rho=\rho$ in the critical window and $\bar\rho=0$ for a
fixed subcritical rate.  The boundary rate is then at most
$\bar c_\varepsilon=e^{\bar\rho\varepsilon^{1/2}}$ for small
$\varepsilon$.  Choose
$\lambda>\max\{C,0,2\bar\rho\}$ and let $H_\varepsilon$ be the scaled,
pinned stationary profile for $\bar c_\varepsilon$ and
$g_\varepsilon=e^{\lambda\varepsilon^{1/2}/2}$, independent of the future.
There are deterministic $C<d<\lambda$ and $D<\infty$, and tight
nonnegative random variables $K_\varepsilon$, such that
\begin{equation}
 dz-K_\varepsilon\leq H_\varepsilon(z)
                  \leq Dz+K_\varepsilon,\qquad z\geq0.
 \label{eq:png-stationary-global-linear-bounds}
\end{equation}
Indeed, write $\delta=\varepsilon^{1/2}$.  The two-line formula gives
$H_\varepsilon(z)\geq\delta S_1(z\delta^{-2})-\delta Y$, with
$Y\sim\operatorname{Geo}(e^{(\bar\rho-\lambda/2)\delta})$.
For a sufficiently small fixed $\beta>0$, the cumulant per unit distance
of $-\beta[\delta S_1(z\delta^{-2})-dz]$ equals
$-\beta(\lambda-d)+\beta^2+O(\delta)<0$.
Its exponential maximal inequality and the uniform exponential tail of
$\delta Y$ prove the lower bound.  Applying the same inequality with
both signs to each of the two walks bounds their absolute values by
a deterministic linear function plus tight constants.  Substitution in
the running-maximum formula proves the upper bound.

Let $Z_\varepsilon$ be the centered evolution from $H_\varepsilon$
with boundary rate $\bar c_\varepsilon$.  The one-parameter stationary limit,
\eqref{eq:png-stationary-global-linear-bounds}, and
Corollary~\ref{cor:png-random-general-data-kpz-limit} give tightness of
$Z_\varepsilon(t,0)$ for each fixed $t>0$.
Stationarity gives
\begin{equation*}
 Z_\varepsilon(t,\cdot)-Z_\varepsilon(t,0)
 \overset{\mathrm d}=H_\varepsilon.
\end{equation*}
Since $f_\varepsilon\leq H_\varepsilon+C+K_\varepsilon$,
the variational formula and boundary thinning imply
$U_\varepsilon(t,z) \leq Z_\varepsilon(t,z)+C+K_\varepsilon$.
The stationary upper bound and fixed-time tightness at zero prove
\eqref{eq:png-evolved-linear-bound}.

For the last assertion, cover a compact time interval by a finite time
grid of mesh $q$, extending it one step to the right.  On the event in
Lemma~\ref{lem:png-uniform-time-buffers}, if $t_j$ is the next grid
time after $t$, composition gives, simultaneously on a compact spatial
interval,
$Z_\varepsilon(t,y) \leq Z_\varepsilon(t_j,y)+C_\alpha q^\alpha$.
At the finitely many grid times, the height at zero is tight and the
stationary spatial profiles satisfy
\eqref{eq:png-stationary-global-linear-bounds}.  First choose $q$ so
that the exceptional probability is small, then bound these finitely
many variables.  This proves the compact upper bound.
\end{proof}

\begin{lemma}[Localization under a common linear upper bound]
\label{lem:png-linear-data-localization}
Fix a critical-window or fixed-subcritical regime.  Suppose
$f_\varepsilon\to f\not\equiv-\infty$ locally in hypograph and
$f_\varepsilon(z)\leq C(1+z)$.  For every compact
$\mathcal K\subset(0,\infty)\times[0,\infty)$,
\begin{equation*}
 \lim_{R\to\infty}\limsup_{\varepsilon\downarrow0}
 \mathbb P\left(
 \sup_{z>R}\{f_\varepsilon(z)+L_\varepsilon^\rho(z,0;y,t)\}
 >\sup_{0\leq z\leq R}\{f_\varepsilon(z)+L_\varepsilon^\rho(z,0;y,t)\}
 \text{ for some }(t,y)\in\mathcal K\right)=0.
\end{equation*}
\end{lemma}

\begin{proof}
Use the dominating critical-window stationary profile from the preceding
lemma.  Its lower bound gives, simultaneously for $z>R$,
$H_\varepsilon(z)-f_\varepsilon(z) \geq(d-C)R-C-K_\varepsilon$.
By the variational formula and boundary thinning, the tail contribution
is therefore at most
$Z_\varepsilon(t,y)-(d-C)R+C+K_\varepsilon$.
Its supremum over $\mathcal K$ tends to $-\infty$ in probability
as $R\to\infty$, uniformly for small $\varepsilon$.

Choose $z_0$ with $f(z_0)>-\infty$ and a hypograph recovery sequence
$z_\varepsilon\to z_0$ with
$f_\varepsilon(z_\varepsilon)\to f(z_0)$.
The comparison with full-space passage gives
\begin{equation*}
 L_\varepsilon^\rho(z_\varepsilon,0;y,t)
 \geq L_{\sqrt\varepsilon}^{\mathrm{full}}(z_\varepsilon,0;y,t)
       -L_{\sqrt\varepsilon}^{\mathrm{full}}(0,0;0,t)+W_{0,\varepsilon}(t),
\end{equation*}
where $W_{0,\varepsilon}(t)$ is the centered wall-to-wall passage
with no boundary atoms.  The two full-space terms have compact upper and
lower bounds by Lemma~\ref{lem:png-fullspace-diagonal-moments}.
The wall lower-tail estimate bounds $W_{0,\varepsilon}$ below at
finitely many positive grid times.  Between consecutive grid times,
composition and Lemma~\ref{lem:png-uniform-time-buffers}, applied at
boundary rate zero, give the same lower bound up to $C_\alpha q^\alpha$.  Thus this one-source contribution
has a tight lower infimum over $\mathcal K$.
For large $R$ the source lies in $[0,R]$, whereas the tail tends to
$-\infty$.  This proves the assertion.
\end{proof}

\begin{proposition}[Fixed-time transition limit in the subcritical regime]
\label{prop:png-subcritical-general-profile-limit}
The conclusion of Theorem~\ref{thm:png-general-data-kpz-limit} also
holds for a fixed boundary rate $c_\star\in[0,1)$, with limiting
parameter $\rho=-\infty$.
\end{proposition}

\begin{proof}
Fix $t>0$ and finitely many terminal positions $y_i\geq0$.
Write $L_\delta(z;y)=L_{\delta^2}^{-\infty}(z,0;y,t)$.
Lemma~\ref{lem:png-linear-data-localization} reduces the variational
supremum to $[0,R]$, with an error probability vanishing as
$R\to\infty$, uniformly for small $\delta$. Choose $R$ large enough
to contain a point where $f$ is finite, and to be a continuity point
of $f$ in the extended-real topology.

Take a partition $\pi:0=a_0<a_1<\cdots<a_N=R$ whose positive
endpoints are continuity points of $f$, and put
\begin{equation*}
 \begin{gathered}
 |\pi|:=\max_{0\leq j<N}(a_{j+1}-a_j),\\
 M_{j,\delta}:=\sup_{[a_j,a_{j+1}]}f_\delta,\qquad
 M_j:=\sup_{[a_j,a_{j+1}]}f,\qquad 0\leq j<N.
 \end{gathered}
\end{equation*}
Local hypograph convergence gives $M_{j,\delta}\to M_j$ in
$\mathbb R\cup\{-\infty\}$. The upper bound follows by compactness
and the hypograph limsup inequality. For the lower bound, approximate
an endpoint maximizer from inside the cell using continuity, and
then use a recovery sequence at an interior point. At the endpoint
$0$, a recovery sequence already lies in the first cell. Such
partitions can have arbitrarily small mesh, since the continuity
points of an upper-semicontinuous function are dense after a bounded
continuous change of its range.

Apply the initial-coordinate version of
Lemma~\ref{lem:png-one-sided-spatial-control} on a compact interval
containing $[0,R]$ and all $y_i$. Outside an event whose probability
vanishes as the auxiliary stationary parameter $\lambda\to\infty$,
it gives, simultaneously in $i$,
\begin{equation*}
 \max_{j<N}\{M_{j,\delta}+L_\delta(a_{j+1};y_i)\}
       -\Omega_{\delta,\lambda}(|\pi|)
 \leq\sup_{0\leq z\leq R}
             \{f_\delta(z)+L_\delta(z;y_i)\} 
 \leq\max_{j<N}\{M_{j,\delta}+L_\delta(a_j;y_i)\}
       +\Omega_{\delta,\lambda}(|\pi|).
\end{equation*}
Here $\Omega_{\delta,\lambda}(h)$ denotes the maximum of the
stationary spatial moduli supplied by the initial-coordinate comparison,
over the finitely many terminal positions on the chosen compact interval.
For each fixed $\lambda$, it satisfies
\begin{equation*}
 \lim_{h\downarrow0}\limsup_{\delta\downarrow0}
 \mathbb P\{\Omega_{\delta,\lambda}(h)>\eta\}=0
 \qquad(\eta>0).
\end{equation*}

For a fixed partition, both bounding vectors converge by the
finite-wedge limits, including terminal positions at zero. Varying
finite heights $M_{j,\delta}\to M_j$ do not change this conclusion:
the max-plus map is Lipschitz in the finitely many initial heights,
and replacing $M_j$ by $\delta\lfloor M_j/\delta\rfloor$ costs at
most $\delta$. Cells with $M_j=-\infty$ may be omitted, since their
heights tend to $-\infty$ and the finitely many single-source
passage values are tight.

Let $f^R=f$ on $[0,R]$ and $f^R=-\infty$ elsewhere. The limiting
left- and right-endpoint wedge profiles, with heights $M_j$, both
converge to $f^R$ in hypograph as $|\pi|\to0$. Indeed, each cell
maximizer is within $|\pi|$ of its assigned endpoint, while upper
semicontinuity supplies the limsup inequality. They have a common
linear upper bound, for example $2C(1+z)$ when $|\pi|\leq1$.
Continuity of the target fixed-point semigroup
\cite[Theorem~4.11]{dauvergne2026directed} therefore makes both
limiting vectors converge to
$(\mathfrak h^{-\infty}(t,y_i;f^R))_i$.
Apply the coordinatewise brackets at continuity points of its joint
distribution function, taking $\delta\downarrow0$, then
$|\pi|\downarrow0$ at fixed $\lambda$, followed by $\lambda\to\infty$.
This proves the transition limit for the truncated supremum.

Finally, $f^R\to f$ locally in hypograph with the same linear upper
bound. Target semigroup continuity and the prelimit localization
lemma allow $R\to\infty$, proving the assertion. 
\end{proof}

\subsubsection{Spatial tightness and the Markov limit}
\label{subsubsec:profile-markov-limit}

Stationary comparison and localization upgrade fixed-time vector laws
to spatial profile laws. Conditioning on the actual PNG profiles and
using independent future Poisson slabs then gives joint laws at
finitely many times. Short vertical-passage bounds control the
intervals between them and complete the positive-time field limit.

\begin{lemma}[Spatial profile convergence]
\label{lem:png-fixed-time-profile-convergence}
Under the hypotheses of Lemma~\ref{lem:png-linear-data-localization},
for each $t>0$,
\begin{equation*}
 U_\varepsilon(t,\cdot)\Longrightarrow
 \mathfrak h^\rho(t,\cdot;f)
 \quad\text{locally uniformly in space},
\end{equation*}
with linear spatial interpolation on a mesh tending to zero.
The interpolation error tends to zero in probability on compact intervals.
\end{lemma}

\begin{proof}
Theorem~\ref{thm:png-general-data-kpz-limit} and
Proposition~\ref{prop:png-subcritical-general-profile-limit} give the
finite-dimensional limits in the two boundary regimes.
Fix $R_{\mathrm{obs}}<\infty$.  Localization allows us, with arbitrarily high
probability, to restrict every maximizing source for $y\in[0,R_{\mathrm{obs}}]$
to $[0,R]$, for a fixed $R$.  Apply
Lemma~\ref{lem:png-one-sided-spatial-control} with
$b\geq R\vee R_{\mathrm{obs}}$.  Taking the supremum over the same initial sources
preserves its inequality: on its comparison event,
\begin{equation*}
 U_\varepsilon(t,y_2)-U_\varepsilon(t,y_1)
 \leq H_{\varepsilon,\lambda,t}(y_2)
              -H_{\varepsilon,\lambda,t}(y_1),
 \qquad 0\leq y_1\leq y_2\leq R_{\mathrm{obs}}.
\end{equation*}
Let
\begin{equation*}
 \Omega_{\varepsilon,\lambda}(h)
 :=\sup_{\substack{u,v\in[0,R_{\mathrm{obs}}+1]\\|u-v|\leq h}}
 |H_{\varepsilon,\lambda,t}(u)-H_{\varepsilon,\lambda,t}(v)|
\end{equation*}
be the spatial modulus of the stationary profile.  For neighboring mesh points
$y^-\leq y\leq y^+$ at distance at most $h$, this gives
\begin{equation*}
 U_\varepsilon(t,y^+)-\Omega_{\varepsilon,\lambda}(h)
 \leq U_\varepsilon(t,y)
 \leq U_\varepsilon(t,y^-)+\Omega_{\varepsilon,\lambda}(h).
\end{equation*}
For fixed $h$, the finite mesh vector converges to that of the
continuous fixed-point profile.  At fixed $\lambda$, the stationary
modulus tends to zero in probability as $h\downarrow0$, by the
one-parameter stationary limit.  First take $\varepsilon\downarrow0$, then
$h\downarrow0$, then $\lambda\to\infty$, and finally remove the
localization error by increasing $R$.  These brackets prove
asymptotic spatial equicontinuity, including negative increments.
Tightness identifies the spatial law, and the same modulus bounds
the interpolation error on any deterministic vanishing mesh.
\end{proof}

In $C_{\rm loc}$ statements, interpolate on expanding spatial intervals,
for example $[0,\varepsilon^{-1/4}]$, replacing any $-\infty$ mesh value
by zero, and extend constantly beyond the right endpoint. A recovered
finite initial source makes every fixed positive-time spatial compact
causally reachable for small $\varepsilon$, so these extensions have
no effect on the local limits.  Markov conditioning always uses the
actual upper-semicontinuous PNG profiles.

Under a common linear upper bound, the target transition is continuous
for local hypograph convergence. Indeed, the parabolic shape bound in
\cite[Proposition~4.10(2)]{dauvergne2026directed} makes the source
contribution tend to $-\infty$ uniformly on positive-time terminal
compacts, while a recovery sequence at one finite initial value gives
a uniform lower candidate. All maximizing sources therefore lie in
one random compact set. Hypograph convergence and continuity of the
landscape then give uniform convergence of the variational maxima.

\begin{proposition}[Convergence at finitely many observation times]
\label{prop:png-markov-profile-convergence}
Under the same initial-data hypotheses, for every
$0<t_1<\cdots<t_k$,
\begin{equation*}
 (U_\varepsilon(t_j,\cdot))_{j=1}^k
 \Longrightarrow
 (\mathfrak h^\rho(t_j,\cdot;f))_{j=1}^k
 \quad\text{in }C_{\rm loc}([0,\infty))^k,
\end{equation*}
with the spatial interpolation just described.
\end{proposition}

\begin{proof}
The preceding lemma gives sequential convergence of the spatial
transition laws whenever deterministic initial profiles converge
locally in hypograph under a common linear upper bound. For a random
current profile $g_\varepsilon$, put
$A_\varepsilon=\max\{0,\sup_{z\geq0}g_\varepsilon(z)/(1+z)\}$
in $[0,\infty]$; its sublevel sets are closed in local hypograph
topology. Suppose the current profile and the past tuple converge
jointly, the limiting profile is nontrivial almost surely, and the
coefficients $A_\varepsilon$ have tight tails. Pass to a subsequence
on which the profiles, past tuple, and coefficients converge jointly,
and take a Skorokhod coupling. The limiting coefficient is finite
almost surely, so the coupled coefficients are eventually bounded.
Deterministic sequential convergence therefore applies along each
coupled sample. Independence of the future slab from the entire past
identifies its conditional transition kernel. Bounded convergence,
after multiplication by bounded continuous functions of the past
tuple, proves joint convergence with the next profile.

Lemma~\ref{lem:png-evolved-linear-bound} gives tight global linear
bounds at each positive observation time, as required for the Markov
iteration in \cite[Lemma~3.10 and Theorem~3.13]{matetski2021kpz}.
Induction over $t_1,\ldots,t_k$, using the Markov property and independent
Poisson increments, gives the fixed-point transition from each current
profile and hence the claimed joint law.
At an intermediate time $s$, restart with the integer-valued profile
$h(s\varepsilon^{-3/2},\cdot)-\lfloor2s\varepsilon^{-3/2}\rfloor$.
Its rescaling differs from $U_\varepsilon(s,\cdot)$ by a single
constant of absolute value at most $\varepsilon^{1/2}$. Vertical
shift covariance makes the finitely many resulting errors vanish.
\end{proof}

\begin{proof}[Proof of Theorem~\ref{thm:intro-half-space-png-kpz-fixed-point-limit}]
Write $U_\varepsilon=\mathfrak h_\varepsilon^\rho$.
Fix a compact rectangle $[a,b]\times[0,R_{\mathrm{obs}}]$ with $a>0$ and a rational
time grid of mesh at most $q$, extending slightly beyond $[a,b]$.
On the event in Lemma~\ref{lem:png-uniform-time-buffers}, composition
gives, for consecutive grid times $t_j\leq t\leq t_{j+1}$,
\begin{equation*}
 U_\varepsilon(t_j,y)-C_\alpha q^\alpha
 \leq U_\varepsilon(t,y)
 \leq U_\varepsilon(t_{j+1},y)+C_\alpha q^\alpha,
 \qquad 0\leq y\leq R_{\mathrm{obs}}.
\end{equation*}
For fixed $q$, let $V_\varepsilon^q$ linearly interpolate in time the
spatially interpolated grid profiles. Let $D_\varepsilon(q)$ be the
largest uniform difference of neighboring raw grid profiles on
$[0,R_{\mathrm{obs}}]$, and $E_\varepsilon(q)$ their largest spatial
interpolation error. On the buffer event the brackets give
$\|U_\varepsilon-V_\varepsilon^q\|_{[a,b]\times[0,R_{\mathrm{obs}}]}
\leq D_\varepsilon(q)+C_\alpha q^\alpha+E_\varepsilon(q)$.
Proposition~\ref{prop:png-markov-profile-convergence} gives joint
convergence of the grid profiles and $E_\varepsilon(q)\to0$ in
probability at fixed $q$. The limiting neighboring differences tend
to zero as $q\downarrow0$ by the space--time continuity in
\cite[Theorem~4.11]{dauvergne2026directed}, and the limiting grid
interpolants converge uniformly to the target field. The buffer
exception vanishes in the same order, first $\varepsilon\downarrow0$
and then $q\downarrow0$. These approximations give space--time
tightness and identify the locally uniform limit.
For continuous representatives, use rectangular meshes of diameter at
most $\varepsilon$ on
$[\varepsilon^{1/8},\varepsilon^{-1/8}]\times[0,\varepsilon^{-1/4}]$,
replace any $-\infty$ nodal value by zero, interpolate bilinearly,
and extend by clamping each coordinate to this rectangle. The
interpolants are Borel $C_{\rm loc}$-valued random elements. A bounded
recovered source $z_\varepsilon$ reaches the whole rectangle for small
$\varepsilon$, since
$\varepsilon^{1/2}(\varepsilon^{-1/4}+z_\varepsilon)
<\varepsilon^{1/8}$ eventually. On each fixed compact, the
interpolation error is bounded by the raw oscillation modulus on a
slightly larger positive-time compact at the mesh diameter, and hence
tends to zero in probability. These rectangles exhaust the domain,
giving exactly the interpolation convention in the introduction.
\end{proof}

\subsection{Convergence to the half-space directed landscape}
\label{subsec:directed-landscape-convergence}

The fixed-point limit supplies finite-wedge laws on rational time
slabs. Together with superadditivity and independence of disjoint
Poisson slabs, these identify the rational passage field through
\cite[Theorem~5.1 and Corollary~5.2]{dauvergne2026directed}.
Stationary comparison in both spatial endpoints and short vertical
passages then give four-variable tightness and locally uniform
metric convergence. When $\rho=-\infty$, fix
$c_\star\in[0,1)$ in
\eqref{eq:intro-png-critical-boundary-scaling}.

\subsubsection{Identification on rational endpoints}
\label{subsubsec:rational-identification}

Set
\begin{equation*}
 D_{\mathbb Q}:=\mathbb H_\uparrow^2\cap\mathbb Q^4.
\end{equation*}
The raw metric $L_\varepsilon^\rho$ is $-\infty$ at noncausal
endpoint pairs. For the countable product space, use the finite-valued
completion
\begin{equation*}
 L_{\varepsilon,\mathbb Q}^\rho(x,s;y,t):=
 \begin{cases}
 L_\varepsilon^\rho(x,s;y,t),
       &\varepsilon^{1/2}|x-y|\leq t-s,\\
 0,&\text{otherwise},
 \end{cases}
 \qquad (x,s;y,t)\in D_{\mathbb Q}.
\end{equation*}
Every fixed compact subset of $\mathbb H_\uparrow^2$ is causally
reachable for all sufficiently small $\varepsilon$, so this
completion does not affect any of the limits below. We retain the
raw metric in all exact variational and composition identities.

For a nonempty finite $P\subset\mathbb Q_{\geq0}$ and
initial heights $a:P\to\mathbb Q$, put
\begin{equation*}
 f_{P,a}(x):=
 \begin{cases}
 a(x),&x\in P,\\
 -\infty,&x\notin P.
 \end{cases}
\end{equation*}
Theorem~\ref{thm:intro-half-space-png-kpz-fixed-point-limit},
applied to these finite-wedge data, gives, for rational $\tau>0$
and finite $Y\subset\mathbb Q_{\geq0}$,
\begin{equation}
 \left(\max_{x\in P}\{a(x)+L_\varepsilon^\rho(x,0;y,\tau)\}\right)_{y\in Y}
 \Longrightarrow
 \left(\mathfrak h^\rho(\tau,y;f_{P,a})\right)_{y\in Y}.
 \label{eq:png-finite-wedge-landscape-input}
\end{equation}
Rounding the microscopic initial heights costs at most
$\varepsilon^{1/2}$. Time translation gives the same conclusion on
every rational time slab.

\begin{theorem}[Identification on a countable dense set]
\label{thm:png-rational-landscape-identification}
Fix $\rho\in\mathbb R\cup\{-\infty\}$. Then
\begin{equation}
 L_{\varepsilon,\mathbb Q}^\rho
 \Longrightarrow \mathcal L^\rho|_{D_{\mathbb Q}}
 \label{eq:png-rational-landscape-convergence}
\end{equation}
in $\mathbb R^{D_{\mathbb Q}}$ with its countable product topology.
In the subcritical case this holds for every fixed $c_\star\in[0,1)$.
\end{theorem}

\begin{proof}
Fix an arbitrary sequence $\varepsilon_n\downarrow0$ and enumerate
$D_{\mathbb Q}=\{u_1,u_2,\ldots\}$. The one-source case of
\eqref{eq:png-finite-wedge-landscape-input} gives convergence of each
coordinate, hence tightness of its whole sequence of laws.
For $\eta>0$, choose $R_j<\infty$ with
\begin{equation*}
 \sup_n\mathbb P\bigl(
 |L_{\varepsilon_n,\mathbb Q}^\rho(u_j)|>R_j\bigr)
 \leq\eta2^{-j}.
\end{equation*}
The product $\prod_j[-R_j,R_j]$ is compact in the product topology
and has probability at least $1-\eta$. Thus the completed fields
are tight.

Let $M^\rho$ be any subsequential limit. We verify precisely the
three conditions in \cite[Definition~4.7]{dauvergne2026directed}.
First, centered superadditivity gives
\begin{equation*}
 L_\varepsilon^\rho(x,s;z,r)
 +L_\varepsilon^\rho(z,r;y,t)
 \leq L_\varepsilon^\rho(x,s;y,t),\qquad s<r<t.
\end{equation*}
Each rational triple is eventually causal, and the deterministic
centerings add. The inequality therefore passes to the limit
simultaneously for all rational triples.

Second, the fields on disjoint half-open Poisson time slabs are
independent. The deterministic noncausal completion preserves this
property. Independence passes to joint limits of finite collections
of coordinates; a monotone-class argument then gives independence
of the countable cylinder sigma-fields of the limiting slab fields.

Third, finite maxima are continuous functions of finitely many
metric entries. The continuous mapping theorem and
\eqref{eq:png-finite-wedge-landscape-input} give, for rational $s<t$,
\begin{equation*}
 \left(\max_{x\in P}\{a(x)+M^\rho(x,s;y,t)\}\right)_{y\in Y}
 \overset{\mathrm d}=
 \left(\mathfrak h^\rho(t-s,y;f_{P,a})\right)_{y\in Y}.
\end{equation*}
Thus $M^\rho$ satisfies the required finite-wedge tests and is a
half-space pre-landscape.
For finite $\rho$,
\cite[Theorem~5.1]{dauvergne2026directed} identifies its law as
$\mathcal L^\rho|_{D_{\mathbb Q}}$.

For $\rho=-\infty$,
\cite[Lemma~4.9]{dauvergne2026directed} first gives a continuous
extension of $M^{-\infty}$ to $\mathbb H_\uparrow^2$.
Fix $\rho'\in\mathbb R$. Boundary thinning couples the prelimit
fields so that
$L_{\varepsilon_n,\mathbb Q}^{-\infty}
\leq L_{\varepsilon_n,\mathbb Q}^{\rho'}$ for all sufficiently
large $n$, since $c_\star<e^{\rho'\varepsilon_n^{1/2}}$.
Pass to a further joint subsequence. The second field has limiting
law $\mathcal L^{\rho'}|_{D_{\mathbb Q}}$ by the finite-parameter
case, and the coordinatewise inequality passes to the limit.
The unique continuous extensions preserve this inequality by density.
Consequently the extension of $M^{-\infty}$ is stochastically
dominated by $\mathcal L^{\rho'}$ for each finite $\rho'$.
Separate couplings for different $\rho'$ suffice here.
The one-source tests already give the subcritical narrow-wedge
marginals at rational endpoint pairs. These are the two hypotheses
of \cite[Corollary~5.2]{dauvergne2026directed}, which identifies
the law as $\mathcal L^{-\infty}$.

The asserted law is therefore the unique subsequential limit in
both boundary regimes.
\end{proof}

\subsubsection{Local uniform metric convergence}
\label{subsubsec:metric-tightness}

For $b\geq2$, write
\begin{equation*}
 K_b=\{(x,s;y,t):0\leq x,y\leq b,\ |s|,|t|\leq b,\
                              t-s\geq b^{-1}\}.
\end{equation*}
These compact sets exhaust $\mathbb H_\uparrow^2$. For a
real-valued field $G$ on $K_b$, define
\begin{equation*}
 \omega_b(G,r):=
 \sup_{\substack{u,v\in K_b\\\|u-v\|_\infty\leq r}}|G(u)-G(v)|,
 \qquad \|u-v\|_\infty:=\max_{1\leq j\leq4}|u_j-v_j|.
\end{equation*}
For a field $F$ on a spatial square $[0,b]^2$, write instead
\begin{equation*}
 \omega_b^{\mathrm{sp}}(F,r):=
 \sup_{\substack{\boldsymbol x,\boldsymbol y\in[0,b]^2\\
                  \|\boldsymbol x-\boldsymbol y\|_\infty\leq r}}
 |F(\boldsymbol x)-F(\boldsymbol y)|.
\end{equation*}

\begin{lemma}[Spatial tightness on a rational time slab]
\label{lem:png-fixed-slab-spatial-tightness}
Fix rational $s<t$, $0\leq b<\infty$, and either boundary regime. For
$F_\varepsilon(x,y):=L_\varepsilon^\rho(x,s;y,t)$,
\begin{equation*}
 \lim_{r\downarrow0}\limsup_{\varepsilon\downarrow0}
 \mathbb P\bigl(\omega_b^{\mathrm{sp}}(F_\varepsilon,r)>\eta\bigr)=0,
 \qquad \eta>0.
\end{equation*}
Moreover, for every fixed finite collection of rational slabs,
their spatial fields converge jointly, locally uniformly in
$(x,y)$, to the corresponding restrictions of the same
$\mathcal L^\rho$, with continuous spatial interpolation.
\end{lemma}

\begin{proof}
Choose rational $R_{\mathrm{cmp}}>b$. Apply
Lemma~\ref{lem:png-one-sided-spatial-control} in both spatial
coordinates, using time reversal for the initial coordinate.
The two stationary comparison profiles
$S_{\varepsilon,\lambda}^{(1)}$ and
$S_{\varepsilon,\lambda}^{(2)}$ need not be independent.
Let $\Omega_{\varepsilon,\lambda}(r)$ be the larger of their
spatial moduli on $[0,R_{\mathrm{cmp}}]$. For fixed $\lambda$ it
vanishes in probability as $\varepsilon\downarrow0$ and then
$r\downarrow0$. The union of the two comparison exceptions has
probability tending to zero as $\varepsilon\downarrow0$ and then
$\lambda\to\infty$.

Take the finite mesh $r\mathbb N_0^2\cap[0,R_{\mathrm{cmp}}]^2$,
with positive rational spacing $r$ small enough that the
corners of every mesh cell meeting $[0,b]^2$ lie in
$[0,R_{\mathrm{cmp}}]^2$. On the comparison event, the two
one-sided increment bounds give
\begin{equation*}
 F_\varepsilon(x^+,y^+)-2\Omega_{\varepsilon,\lambda}(r)
 \leq F_\varepsilon(x,y)
 \leq F_\varepsilon(x^-,y^-)+2\Omega_{\varepsilon,\lambda}(r),
\end{equation*}
where $x^-\leq x\leq x^+$ and $y^-\leq y\leq y^+$ are the
neighboring mesh coordinates. Consequently
\begin{equation*}
 \omega_b^{\mathrm{sp}}(F_\varepsilon,r)
 \leq
 \max_{\substack{p,p'\in r\mathbb N_0^2\cap[0,R_{\mathrm{cmp}}]^2\\
                        \|p-p'\|_\infty\leq3r}}
 |F_\varepsilon(p)-F_\varepsilon(p')|
       +4\Omega_{\varepsilon,\lambda}(r).
\end{equation*}
For fixed $r$, Theorem~\ref{thm:png-rational-landscape-identification}
gives joint convergence of the finite mesh vector. Its limiting
mesh oscillation tends to zero as $r\downarrow0$ by continuity of
$\mathcal L^\rho$. Taking $\varepsilon\downarrow0$, then
$r\downarrow0$ at fixed $\lambda$, and finally
$\lambda\to\infty$ proves the modulus assertion.
Use $L_{\varepsilon,\mathbb Q}^\rho$ at rational mesh vertices for the
spatial interpolants; the completion is inactive on every fixed
spatial compact for small $\varepsilon$.
One-point tightness and this modulus give tightness of bilinear
spatial interpolations and a vanishing uniform interpolation error.

For finitely many rational slabs, marginal tightness gives product
tightness. Every product subsequential limit has the joint rational
spatial evaluations prescribed by
Theorem~\ref{thm:png-rational-landscape-identification}.
Evaluations on a countable dense set determine the joint law of
continuous fields. This proves joint convergence to the restrictions
of the same landscape.
\end{proof}

\begin{lemma}[Local uniform tightness of the directed metric]
\label{lem:png-metric-space-time-tightness}
For every $\rho\in\mathbb R\cup\{-\infty\}$, $b\geq2$, and $\eta>0$,
\begin{equation*}
 \lim_{r\downarrow0}\limsup_{\varepsilon\downarrow0}
 \mathbb P\bigl(\omega_b(L_\varepsilon^\rho,r)>\eta\bigr)=0.
\end{equation*}
There are Borel continuous mesh interpolations
$\widetilde L_\varepsilon^\rho$ that are tight in
$C_{\mathrm{loc}}(\mathbb H_\uparrow^2)$ and satisfy
\begin{equation*}
 \|\widetilde L_\varepsilon^\rho-L_\varepsilon^\rho\|_{K_b}
 \longrightarrow0\quad\text{in probability}.
\end{equation*}
\end{lemma}

\begin{proof}
Fix $b\geq2$ and a positive rational mesh
$q<(8b)^{-1}$. Let $\mathcal T_{b,q}$ be the finite set of
slabs $(a,d)$ with $a,d\in q\mathbb Z\cap[-b-1,b+1]$ and
$d-a\geq(2b)^{-1}$.
For $u=(x,s;y,t)$ and $v=(x',s';y',t')$ in $K_b$ with
$\|u-v\|_\infty\leq q$, choose common grid times
\begin{equation*}
 \begin{aligned}
 a^-&=q\lfloor\min(s,s')/q\rfloor,
 &a^+&=q\lceil\max(s,s')/q\rceil,\\
 d^-&=q\lfloor\min(t,t')/q\rfloor,
 &d^+&=q\lceil\max(t,t')/q\rceil.
 \end{aligned}
\end{equation*}
The added initial and terminal intervals have length at most $2q$.
Also $d^--a^+\geq b^{-1}-4q>(2b)^{-1}$, so both
$(a^+,d^-)$ and $(a^-,d^+)$ lie in $\mathcal T_{b,q}$.
Both resulting long-slab metric arguments lie in $K_{2b+2}$.

Put $e_q=2C_\alpha(2q)^\alpha$, with
$\alpha\in(0,1/3)$ as in
Lemma~\ref{lem:png-uniform-time-buffers}.
On $\mathcal E_\varepsilon^{\mathrm{time}}(b+1,2q)$,
centered superadditivity gives
\begin{equation*}
 L_\varepsilon^\rho(x,a^+;y,d^-)-e_q
 \leq L_\varepsilon^\rho(u)
 \leq L_\varepsilon^\rho(x,a^-;y,d^+)+e_q,
\end{equation*}
and the same bounds hold for $v$. Zero-duration vertical passages
are omitted. These exact comparisons use the raw metric.

Define the finite-slab oscillation
\begin{equation*}
 A_{\varepsilon,b}(q):=
 \max_{\substack{(a,d),(a',d')\in\mathcal T_{b,q}\\
                    |a-a'|\vee|d-d'|\leq2q}}
 \ \sup_{\substack{x,x',y,y'\in[0,b]\\
                    |x-x'|\vee|y-y'|\leq q}}
 |L_\varepsilon^\rho(x,a;y,d)
           -L_\varepsilon^\rho(x',a';y',d')|.
\end{equation*}
The preceding brackets imply
\begin{equation*}
 \omega_b(L_\varepsilon^\rho,q)
       \leq A_{\varepsilon,b}(q)+2e_q
\end{equation*}
on the time-buffer event. For fixed $q$,
Lemma~\ref{lem:png-fixed-slab-spatial-tightness} gives joint uniform
convergence on the finitely many slabs in $\mathcal T_{b,q}$.
Thus $A_{\varepsilon,b}(q)$ converges in law to the same
functional of $\mathcal L^\rho$. Its limit is bounded by
$\omega_{2b+2}(\mathcal L^\rho,2q)$, which tends to zero almost
surely as $q\downarrow0$. The buffer exceptional probability tends
to zero by Lemma~\ref{lem:png-uniform-time-buffers}. This proves the
asserted modulus estimate.

Choose deterministic rational spacings
$\eta_\varepsilon\downarrow0$ on the rectangular grid in
$(x,s,y,t)$, with $x,y\geq0$. At grid vertices use
$L_\varepsilon^\rho$ when $t>s$ and the endpoints are causal,
and zero otherwise; interpolate multilinearly between vertices.
The resulting $\widetilde L_\varepsilon^\rho$ is continuous and
is Borel measurable as a $C_{\mathrm{loc}}$-valued random element,
since it is determined by countably many measurable nodal values.
For each $K_b$, all cells meeting it have causal vertices for
small $\varepsilon$ and lie in $K_{2b+2}$.
The uniform interpolation error is then at most
$\omega_{2b+2}(L_\varepsilon^\rho,\eta_\varepsilon)$ and tends to zero
in probability. One-point tightness, the modulus estimate, and
Arzel\`a--Ascoli give tightness on every $K_b$; a diagonal argument
gives tightness in $C_{\mathrm{loc}}(\mathbb H_\uparrow^2)$.
\end{proof}

\begin{proof}[Proof of Theorem~\ref{thm:intro-half-space-png-directed-landscape-limit}]
Lemma~\ref{lem:png-metric-space-time-tightness} gives tight
continuous interpolations with vanishing local uniform error.
Theorem~\ref{thm:png-rational-landscape-identification} identifies
the joint rational evaluations of every continuous subsequential
limit as those of $\mathcal L^\rho$. These evaluations determine
the Borel law on $C_{\mathrm{loc}}(\mathbb H_\uparrow^2)$.
Hence the interpolations converge to $\mathcal L^\rho$, and their
vanishing error proves the asserted convergence of the raw metric.
Both arguments apply to every fixed $c_\star\in[0,1)$ when
$\rho=-\infty$.
\end{proof}

\subsubsection{Joint convergence from several initial profiles}
\label{subsubsec:joint-profile-limit}

Joint tightness retains the coupling between the metric and finitely
many height evolutions. Localization and the exact variational formula
identify every limiting height as a function of the same landscape.

\begin{corollary}[Joint convergence from several initial profiles]
\label{cor:png-general-initial-transfer}
Fix $\rho\in\mathbb R\cup\{-\infty\}$ and $k\in\mathbb N$.
For $1\leq j\leq k$, take deterministic PNG initial profiles
$h_{0,\varepsilon,j}\in\mathcal X_{\mathrm{adm}}$ and put
\begin{equation*}
 f_{\varepsilon,j}(z):=\varepsilon^{1/2}
 h_{0,\varepsilon,j}(\varepsilon^{-1}z).
\end{equation*}
Suppose $f_{\varepsilon,j}\to f_j\in\mathrm{UC}_+^{\mathbb R}$
locally in hypograph and
$f_{\varepsilon,j}(z)\leq C(1+z)$ for one finite $C$.
For the same PNG directed metric, put
\begin{equation*}
 U_{\varepsilon,j}(t,y):=
 \sup_{z\geq0}\{f_{\varepsilon,j}(z)
                       +L_\varepsilon^\rho(z,0;y,t)\}.
\end{equation*}
Then
\begin{equation*}
 (L_\varepsilon^\rho,U_{\varepsilon,1},\ldots,U_{\varepsilon,k})
 \Longrightarrow
 (\mathcal L^\rho,\mathfrak h^\rho(\cdot,\cdot;f_1),\ldots,
                   \mathfrak h^\rho(\cdot,\cdot;f_k))
\end{equation*}
in
\begin{equation*}
 C_{\mathrm{loc}}(\mathbb H_\uparrow^2)
 \times C_{\mathrm{loc}}((0,\infty)\times[0,\infty))^k,
\end{equation*}
using continuous interpolations with vanishing local uniform errors.
All limiting evolutions use the same landscape through
\eqref{eq:intro-half-space-kpz-variational}.
For $\rho=-\infty$, the assertion holds for each fixed
$c_\star\in[0,1)$.
\end{corollary}

\begin{proof}
Theorems~\ref{thm:intro-half-space-png-kpz-fixed-point-limit}
and~\ref{thm:intro-half-space-png-directed-landscape-limit}
give joint tightness of the continuous metric and height
interpolations $\widetilde L_\varepsilon^\rho$ and
$\widetilde U_{\varepsilon,j}$ in this finite product. Consider any joint
subsequential limit $(\mathcal L^\rho,V_1,\ldots,V_k)$ and use
Skorokhod representation for these continuous fields. Include
their interpolation errors on a countable compact exhaustion in
the coupling; these errors converge to zero in probability and,
on a further subsequence, almost surely.

Fix $j$, $(t,y)\in\mathbb Q_{>0}\times\mathbb Q_{\geq0}$,
and an integer $R\geq1$.
Uniform metric convergence and the hypograph limsup
inequality give
\begin{equation*}
 \limsup_{\varepsilon\downarrow0}
 \sup_{0\leq z\leq R}
 \{f_{\varepsilon,j}(z)+\widetilde L_\varepsilon^\rho(z,0;y,t)\}
 \leq
 \sup_{0\leq z\leq R}
 \{f_j(z)+\mathcal L^\rho(z,0;y,t)\}.
\end{equation*}
Indeed, approximate maximizers have a convergent subsequence in
$[0,R]$. For the reverse bound, every $z\in[0,R]$ with
$f_j(z)>-\infty$ has a hypograph recovery sequence
$z_\varepsilon\to z$, contained in $[0,2R]$ for small
$\varepsilon$. Hence
\begin{equation*}
 \liminf_{\varepsilon\downarrow0}
 \sup_{0\leq z\leq2R}
 \{f_{\varepsilon,j}(z)+\widetilde L_\varepsilon^\rho(z,0;y,t)\}
 \geq
 \sup_{0\leq z\leq R}
 \{f_j(z)+\mathcal L^\rho(z,0;y,t)\}.
\end{equation*}
Choose a finite sum $e_\varepsilon(R)$ of the coupled exhaustion
errors that bounds the sum of the height interpolation error at
$(t,y)$ and the metric interpolation error on
$\{(z,0;y,t):0\leq z\leq2R\}$. This sum tends to zero almost
surely. Only the tail
where this compact is causal is relevant. The exact prelimit
variational formula gives
\begin{equation*}
 \widetilde U_{\varepsilon,j}(t,y)
 \geq \sup_{0\leq z\leq2R}
       \{f_{\varepsilon,j}(z)+\widetilde L_\varepsilon^\rho(z,0;y,t)\}
       -e_\varepsilon(R).
\end{equation*}
This inequality concerns only the represented fields and error
variables, so it remains valid in the coupling. The second
compact-source bound, followed by $R\to\infty$, implies
$V_j(t,y)\geq\mathfrak h^\rho(t,y;f_j)$ almost surely.

For the upper bound, the same variational formula gives, on the
original probability space,
\begin{equation*}
 \mathbb P\left(\widetilde U_{\varepsilon,j}(t,y)>
    \sup_{0\leq z\leq R}
       \{f_{\varepsilon,j}(z)+\widetilde L_\varepsilon^\rho(z,0;y,t)\}
       +e_\varepsilon(R)\right) 
 \leq
 \mathbb P(\text{localization to }[0,R]\text{ fails}).
\end{equation*}
The event on the left is measurable in the represented variables,
so its probability is unchanged by the coupling.
The first compact-source bound implies that
$\{V_j(t,y)>\mathfrak h^\rho(t,y;f_j)+\eta\}$ is contained in
the limit inferior of these comparison events. Fatou's lemma gives
\begin{equation*}
 \mathbb P\bigl(V_j(t,y)>
                    \mathfrak h^\rho(t,y;f_j)+\eta\bigr)
 \leq
 \limsup_{\varepsilon\downarrow0}
 \mathbb P(\text{localization to }[0,R]\text{ fails}).
\end{equation*}
Lemma~\ref{lem:png-linear-data-localization} makes the right-hand
side tend to zero as $R\to\infty$. Then let $\eta\downarrow0$.
Thus equality holds
almost surely at every rational terminal point, simultaneously
for the finitely many $j$. Continuity of both limiting height
fields extends the equality everywhere.
Every joint subsequential limit is consequently the claimed
common-landscape coupling, which proves the corollary.
\end{proof}

\end{document}